\documentclass[11pt]{article} % use larger type; default would be 10pt
\usepackage[utf8]{inputenc} % set input encoding (not needed with XeLaTeX)
\usepackage{amscd, amsmath, amssymb, amsfonts, epic, amsthm}
\usepackage{graphicx}
\usepackage{tikz}
\usepackage{tikz-cd} %for chain maps
\usetikzlibrary{patterns}
\usepackage{a4wide}
\usepackage{paralist}
\usepackage{booktabs}
\usepackage{multirow}
\usepackage{geometry} % to change the page dimensions
\usepackage{setspace}
\usepackage{graphicx} % support the \includegraphics command and options

\usepackage{booktabs} % for much better looking tables
\usepackage{array} % for better arrays (eg matrices) in maths
\usepackage{paralist} % very flexible & customisable lists (eg. enumerate/itemize, etc.) 
\usepackage{verbatim} % adds environment for commenting out blocks of text & for better verbatim
\usepackage{subfig} % make it possible to include more than one captioned figure/table in a single float
\newtheoremstyle{exampstyle}
{0em} % Space above
{0em} % Space below
{} % Body font
{} % Indent amount
{\bfseries} % Theorem head font
{.} % Punctuation after theorem head
{0.5em} % Space after theorem head
{} % Theorem head spec (can be left empty, meaning `normal')

\theoremstyle{exampstyle} \newtheorem{theorem}{Theorem}[section]

\newtheorem{lem}[theorem]{Lemma}
\newtheorem{exm}[theorem]{Example}
\newtheorem{thm}[theorem]{Theorem}
\newtheorem{cor}[theorem]{Corollary}

\newtheorem{defn}[theorem]{Definition}
\newtheorem{remark}[theorem]{Remark}

\newtheorem{setup}[theorem]{Setup}

\usepackage{amsopn}

\usepackage[mathcal]{euscript}

\newcommand{\IR}{{\mathbb{R}}}
\newcommand{\IZ}{{\mathbb{Z}}}

\newcommand{\IN}{{\mathbb{N}}}

\newcommand{\kh}{{\mathcal H}}
\newcommand{\ki}{{\mathcal I}}
\newcommand{\ka}{{\mathcal A}}
\newcommand{\kb}{{\mathcal B}}
\newcommand{\kC}{{\mathcal C}}
\newcommand{\kl}{{\mathcal L}}

\newcommand{\kn}{{\mathcal N}}
\newcommand{\kg}{{\mathcal G}}

\newcommand{\kS}{{\mathcal S}}
\newcommand{\kt}{{\mathcal T}}
\newcommand{\ku}{{\mathcal U}}
\newcommand{\kv}{{\mathcal V}}

\newcommand{\inv}{^{-1}}

\DeclareMathOperator{\Hom}{\mathsf{Hom}}
\renewcommand{\setminus}{\smallsetminus}

\newcommand{\sC}{\mathsf{C}}
\newcommand{\sD}{\mathsf{D}}

\newcommand{\sK}{\mathsf{K}}
\newcommand{\sk}{\mathsf{k}}

\newcommand{\sT}{\mathsf{T}}
\newcommand{\sU}{\mathsf{U}}

\DeclareMathAlphabet{\mathpzc}{OT1}{pzc}{m}{it}

\usepackage{mathrsfs}

\usepackage[utf8]{inputenc}
\usepackage[english]{babel}
\usepackage{authblk}
\newcommand{\Addresses}{{% additional braces for segregating \footnotesize
  \bigskip
  \footnotesize

  \textsc{School of Engineering, Computing and Mathematics, University of Plymouth,
    Drake Circus, Plymouth, PL4 8AA}\newline\nopagebreak
 Email: \texttt{callum.page@plymouth.ac.uk}}}

\newcommand{\bib}[6]{{\bibitem{#2} #3: {\emph{#4},} #5#6.}}
\usepackage{mathrsfs}
\usetikzlibrary{arrows}
\usepackage{blindtext}

\usepackage{faktor}\usepackage{amsmath}\usepackage{amssymb}

\DeclareMathOperator*{\dprime}{\prime \prime}
\makeatletter
\DeclareRobustCommand*{\mfaktor}[3][]
{
   { \mathpalette{\mfaktor@impl@}{{#1}{#2}{#3}} }
}
\newcommand*{\mfaktor@impl@}[2]{\mfaktor@impl#1#2}
\newcommand*{\mfaktor@impl}[4]{
   \settoheight{\faktor@zaehlerhoehe}{\ensuremath{#1#2{#3}}}%
   \settoheight{\faktor@nennerhoehe}{\ensuremath{#1#2{#4}}}%
      \raisebox{-0.5\faktor@zaehlerhoehe}{\ensuremath{#1#2{#3}}}%
      \mkern-4mu\diagdown\mkern-5mu%
      \raisebox{0.5\faktor@nennerhoehe}{\ensuremath{#1#2{#4}}}%
}
\makeatother
\newcommand{\pf}{\textbf{Proof: }}

\newcommand{\simgen}{\sim_{\text{gen}}}

\newcommand{\dhomd}{\hom_{\sD/\Sigma}}
\newcommand{\td}{\text{thick}_\sD}
\newcommand{\tds}{\text{thick-st}(\sD)}
\newcommand{\leqgen}{\leq_{\text{gen}}}

\newcommand{\acs}{\text{Arc}({\Sigma})}
\usepackage{parskip}
\usepackage{pgfplots}
\pgfplotsset{compat=1.15}
\usepackage{mathrsfs}
\usetikzlibrary{arrows}
\usepackage{stackengine}
\newcommand\xleftrightarrow[1]{%
    \mathrel{{\stackon[4pt]{$\longleftrightarrow$}{$\scriptscriptstyle#1$}}}}
\newcommand\restr[2]{{% we make the whole thing an ordinary symbol
  \left.\kern-\nulldelimiterspace % automatically resize the bar with \right
  #1 % the function
  \littletaller % pretend it's a little taller at normal size
  \right|_{#2} % this is the delimiter
  }}

\newcommand{\littletaller}{\mathchoice{\vphantom{\big|}}{}{}{}}
\newcommand{\pac}{pointed arc-collection}
\usepackage{float}
\begin{document}
\begin{center}
    \textbf{THICK SUBCATEGORIES OF DERIVED CATEGORIES OF GENTLE ALGEBRAS}
   
    \small CALLUM T PAGE
\end{center}
\begin{abstract}
We study thick subcategories of derived categories of gentle algebras. Any thick subcategory of a derived category of a gentle algebra is generated by a set of string objects or a set of band objects. We show the thick subcategories generated by string objects are in bijection with sets of non-crossing paths on the geometric model of the derived category.
\end{abstract}
\section{Introduction}
In this paper, we will produce a classification of certain thick subcategories of a bounded derived category of a gentle algebra. We call a subcategory of a triangulated category \textit{thick} when it is full triangulated subcategory which is closed under taking direct summands. Classification of thick subcategories is, in general, a hard problem. However, classifications have been produced in multiple settings. The work of Hopkins and Smith in \cite{Devinatz, Hopkins} considered endomorphisms to classify thick subcategories of a homotopy category with finite spectra. Ingalls and Thomas in \cite{Ingalls} proved that the finitely generated thick subcategories of representations of a quiver without oriented cycles are in bijection with cluster tilting objects in the associated cluster category. This was then built on in \cite{Krause, Antieau} %11, 15
by authors who considered a relationships between the posets of thick subcategories and non-crossing partitions. More recently, in \cite{Elagin1, Elagin}, the work of Elagin and Lunts classified thick subcategories in the derived category of a weighted projective curve by considering them as derived categories of nilpotent representations of some quivers or the orthogonal to an exceptional collection of torsion sheaves. There are several other authors that have considered classifications of thick subcategories (see \cite{Neeman, Thomason, Garkusha, Bruning, Takahashi}).%92, 97, 06, 07, 09

The algebras that became known as \textit{gentle} were introduced by Assem and Skowroński, so they could classify special biserial algebras in \cite{Assem}. In general derived categories are complex objects to study, but an understanding of derived categories of gentle algebras has been developed by various authors.  Bekkert and Merklen described their indecomposable objects in terms of homotopy strings and bands in \cite{Bekkert}. Arnesen, Laking and Pauksztello used the paths in the associated quiver to describe all maps between indecomposable objects in these derived categories \cite{Arnesen}. The cones of these morphisms can be computed using the work of Çanakçı, Pauksztello and Schroll in \cite{Canakci}.

We will also consider geometric models of derived categories of gentle algebras. These were first treated by Haiden, Katzarkov and Kontsevich who considered the Fukaya category of a surface \cite{Haiden}. Lekili and Polishchuk expanded on this work and constructed a graded surface for any homologically smooth (graded) gentle algebra \cite{Lekili}. Independent work of Opper, Plamondon and Schroll gave a very explict formulation of the geometric model for all derived categories of gentle algebras in \cite{Opper}. This is the model that will be used in this paper.

A motivation for our work is Broomhead's classification of the thick subcategories of discrete derived categories (DDCs). The paper \cite{Broomhead} provides an isomorphism between the posets of thick subcategories of a DDC and the configurations of non-crossing arcs on the derived category's geometric model. These derived categories are a class of derived categories of gentle algebras, so we build on this work to achieve our aim. However, there are properties of DDCs that were used by Broomhead in his classification which do not hold in general for the derived category of a gentle algebra. For example, the geometric model corresponding to a DDC is always a cylinder and the hom spaces between any two objects in a DDC are zero- or one-dimenisonal. %Neither of these facts generalise to all derived categories of gentle algebras. 

In this paper, we consider thick subcategories generated by string objects in the derived category. Our classification of these concurs with the complete classification in the DDC case \cite{Broomhead}, because there are no band objects in a DDC. 

We begin, in Section 1, by recalling some background on derived categories of gentle algebras and the classification of their indecomposable objects in terms of homotopy strings and bands. We call an indecomposable object a string (respectively band) object if it corresponds to a homotopy string (respectively band). We then discuss some details of the geometric model from \cite{Opper}. In section 2, we discuss some preliminaries, first we consider curves on the geometric model that we can associate to objects in the unbounded homotopy category and then we look at morphisms between two thick subcategories. In Section 3, we introduce the class of finitely generated Ext-connected thick subcategories which we will go on to classify. (We consider the Ext-connected, finitely generated thick subcategories generated by string and band objects.)

\thm (Theorem~\ref{cstobj}) Let $\sT$ be a Ext-connected thick subcategory containing a string object. Then $\sT$ is generated by string objects. Furthermore, if $\sT$ is finitely generated then it is generated by finitely many  string objects.

This means we have two distinct classes of Ext-connected thick subcategories. Namely, an Ext-connected, finitely generated thick subcategory is either generated by only string objects or only band objects. For the remainder of the paper we consider the former class, which we denote as $\tds$. In other words, we classify the Ext-connected, finitely generated thick subcategories generated by string objects.

\begin{remark}\label{posnlat}
    The set of thick subcategories $\tds$ is a sub-poset of all the thick subcategories of $\sD$ but not a sub-lattice, (see Example~\ref{exm1}).
\end{remark}

In Section 4, we introduce arc-collections, which are non-crossing configurations of paths on the geometric model of the derived category. We discuss concatenations of paths in arc-collection and prove some technical results that will be used later in the paper. Section 5 gives us a further simplification, when considering Ext-connected thick subcategories containing a string object.

\thm (Theorem~\ref{ctwsgen}) Let $\sT$ be a Ext-connected thick subcategory containing a string object. Then $\sT$ is generated by exceptional and spherelike string objects. Furthermore, if $\sT\in\tds$ then it is generated by finitely many exceptional and spherelike string objects.

\begin{remark}\label{nsphexcbds}
    An analogous result does not hold for all thick subcategories generated by band objects, (see Example~\ref{exm2}).
\end{remark} 

Section 6 considers the thick subcategories generated by finitely many exceptional and spherelike string objects. If a thick subcategory is generated by a set of objects corresponding to a set of arcs which form an arc-collection then we say the thick subcategory corresponds to an arc-collection.

\thm (Theorem~\ref{strESc}) $\sT\in \tds$ if and only if $\sT$ corresponds to an connected finite arc-collection.

We then introduce a binary relation $\leqgen$ on the set of arc-collections $\acs$ and consider the induced equivalence relation. The binary relation induces a well-defined partial order on the set of equivalence classes. We conclude the paper by showing the following result.

\thm (Theorem~\ref{mt}) There is an isomorphism of posets:\[\faktor{\acs}{\simgen}\xleftrightarrow{1-1}\tds\]

\subsection*{Acknowledgements} I would like to thank Nathan Broomhead, my PhD supervisor, for all the lengthly chats and for his support at every step in the process of producing this paper. His invaluable advice has greatly improved this work.

\section{Background}
We collect, for the convenience of the reader, details of gentle algebras, their derived categories and the geometric models of such derived categories. Most of these can be found in \cite{Arnesen, Opper} and a reader familiar with such objects can proceed directly to Section 2.
\subsection{The derived category of a gentle algebra}
%If $\sC$ is an abelian category then we denote the bounded derived category of $\sC$ by $\sD^b(\sC)$. The homotopy category $\sK^b(\sC)$ has a full subcategory $\sN$ such that $X\in\sN$ if $H_n(X)=0$ for every degree $n\in\IZ$. The definition of derived categories says \[\sD^b(\sC)\cong\faktor{\sK^b(\sC)}{\sN},\] so quasi-isomorphisms in $\sK^b(\sC)$ become true isomorphisms in $\sD^b(\sC)$. 

Our main object of study will be the bounded derived category of finitely generated modules over a gentle algebra $\Lambda$. %These algebras are path algebras of a quiver modulo relations. We say a quiver $Q$ is equal to the quadruple $(Q_0,Q_1,s,t)$, where $Q_0$ is the set of vertices in $Q$, $Q_1$ is the set of arrows in $Q$, $s:Q_1\rightarrow Q_0$ is a function which sends each arrow to its source and $t:Q_1\rightarrow Q_0$ is a function which sends each arrow to its target. 
We now recall the definition of a gentle algebra.
\begin{defn}
    \cite[Definition 1.7]{Opper} Let $\sk Q$ be the path algebra of a quiver $Q$ over an algebraically closed field $\sk$ and $I$  be an ideal of $\sk Q$. The algebra $\Lambda$ is \textit{gentle} if $\Lambda\cong\sk Q/I$ and \begin{itemize}
        \item $Q=(Q_0, Q_1,s,t)$ is a finite connected quiver;
        \item for any vertex $v\in Q_0$ the number of arrows with a source (respectively target) of $v$ is at most two; 
        \item for $a\in Q_1$ there is at most one arrow $b$ such that $ab\in I$ and at most one arrow $c$ such that $ac\notin I$;
        \item for $a\in Q_1$ there is at most one arrow $b$ such that $ba\in I$ and at most one arrow $c$ such that $ca\notin I$;
        \item $I$ is generated by paths of length 2;
        \item $I$ is an admissible ideal.
    \end{itemize}
\end{defn}

Now we fix some notation for our main objects of study. For an additive category $\sC$, we denote the bounded derived category by $\sD^b(\sC)$ and the bounded homotopy category by $\sK^b(\sC)$. 

Let $\Lambda$ be a gentle algebra. The bounded derived category $\sD^b(\text{mod } \Lambda)$ can be identified with the homotopy category of projective resolutions $\sK^{-}(\text{proj }\Lambda)$ containing right bounded complexes of finitely generated projective modules \cite[Section 1]{Arnesen}. The category $\sK^{-}(\text{proj }\Lambda)$ is equivalent to the bounded homotopy category $\sK^b=\sK^b(\text{proj }\Lambda)$ if every $\Lambda$-module has a finite projective dimension. This happens when the algebra is homologically smooth. It is known that a gentle algebra $\sk Q/I$ is homologically smooth if there are no cycles $c_1c_2...c_n$ ($c_i\in Q_1$) of $Q$ with $c_jc_{j+1}\in I$ for all $j\in\{1,2,...,n-1\}$ and $c_nc_1\in I$ \cite[Lemma 3.1.3]{Lekili}.

For the rest of the paper we will consider $\Lambda$ to be isomorphic to the homologically smooth gentle algebra $\sk Q/I$ and $\sD=\sD^b(\text{mod }\Lambda)$. We can think of objects in $\sD$ as bounded chain complexes of projective modules and the non-zero morphisms as chain maps which are not homotopic to zero.

\defn If there exists a non-zero morphism between $A$ and $B[n]$ for some $n\in\IZ$ then we say $A$ and $B$ have a non-zero morphism between them upto shift.

We now define the dimension of the hom space from $A$ to $B$ in $\sD$ upto shift:\[
 \dhomd(A,B):=\sum_i \hom_\sD(A,B[i])
\]

For a triangulated category $\sC$, if there are two sets $\kS$ and $\kS^\prime$ of objects in $\sC$ such that any object in $\kS$ has no morphism to a shift of an object in $\kS^\prime$ then we denote this by \[
\hom_{\sC/\Sigma}(\kS,\kS^\prime)=0.
\]

\subsection{Homotopy strings and bands}
In this section, we see how homotopy strings and bands describe the objects in the derived category of a gentle algebra $\Lambda=\sk Q/I$. The details of this section can be found in \cite{Arnesen}. We recall some definitions and notation that will be used throughout this paper.

A path $p$ of a bound quiver $(Q, I)$ is a sequence $a_1a_2...a_n$ of arrows in $Q$ with $t(a_i)=s(a_{i+1})$. We say that $p$ is permitted if $p\notin I$, otherwise $p$ is forbidden. For any permitted path $p=a_1a_2...a_n$ there is a direct homotopy letter $\boldsymbol{p}$ with source $s(\boldsymbol{p}):=s(a_1)$ and target $t(\boldsymbol{p}):=t(a_n)$. We define an inverse homotopy letter $\overline{\boldsymbol{p}}$ to be the formal inverse of a direct homotpy letter $\boldsymbol{p}$ with source $s(\overline{\boldsymbol{p}})=t(\boldsymbol{p})$ and target $t(\overline{\boldsymbol{p}})=s(\boldsymbol{p})$. The inverse of an inverse homotopy letter $\overline{\boldsymbol{p}}$ is $\boldsymbol{p}$.

A sequence of homotopy letters $\boldsymbol{\sigma}=\boldsymbol{\sigma}_1\boldsymbol{\sigma}_2...\boldsymbol{\sigma}_n$ is called an \textit{homotopy string} if it satisfies the conditions given in \cite[Section 2]{Arnesen} %for any consecutive letters $\boldsymbol{\sigma}_i$ and $\boldsymbol{\sigma}_{i+1}$ in the sequence $t(\boldsymbol{\sigma}_i)=s(\boldsymbol{\sigma}_{i+1})$ and the following is satisfied:\begin{itemize}
    %\item if $\boldsymbol{\sigma}_i=\boldsymbol{p}$ and $\boldsymbol{\sigma}_{i+1}=\boldsymbol{q}$ are both direct then $pq\in I$;
    %\item if $\boldsymbol{\sigma}_i=\overline{\boldsymbol{p}}$ and $\boldsymbol{\sigma}_{i+1}=\overline{\boldsymbol{q}}$ are both inverse then $qp\in I$;
    %\item if $\boldsymbol{\sigma}_i=\boldsymbol{p}$ is direct and $\boldsymbol{\sigma}_{i+1}=\overline{\boldsymbol{q}}$ is inverse then $p$ and $q$ do not end with the same arrow;
    %\item if $\boldsymbol{\sigma}_i=\overline{\boldsymbol{p}}$ is inverse and $\boldsymbol{\sigma}_{i+1}=\boldsymbol{q}$ is direct then $p$ and $q$ do not start with the same arrow.
%\end{itemize}

\defn An \textit{homotopy band} is an homotopy string $\boldsymbol{\sigma}=\boldsymbol{\sigma}_1\boldsymbol{\sigma}_2...\boldsymbol{\sigma}_n$ such that $\boldsymbol{\sigma}$ has the same number of direct and inverse of homotopy letters, $\boldsymbol{\sigma}_n\boldsymbol{\sigma}_1$ is a homotopy string and $\boldsymbol{\sigma}\ne\boldsymbol{\mu}^m$ for some homotopy string $\boldsymbol{\mu}$ and integer $m>1$.

There exists a grading $\kg_{\boldsymbol{\sigma}}$ of homotopy string or band $\boldsymbol{\sigma}$ which is a sequence of integers satisfying the conditions given in \cite[Definition 2.3]{Opper}.
%\defn The pair is $(\boldsymbol{\sigma}, \kg_{\boldsymbol{\sigma}})$ a graded homotopy string (respectively band) if $\boldsymbol{\sigma}=\boldsymbol{\sigma}_1\boldsymbol{\sigma}_2...\boldsymbol{\sigma}_n$ and for some integer $i_0$ and $\kg_{\boldsymbol{\sigma}}=(i_0,i_1,...,i_n)\in\IN^{n+1}$ (respectively $\kg_{\boldsymbol{\sigma}}=(i_0,i_1,...,i_{n-1})\in\IN^n$) such that\[
%i_{j+1}=\begin{cases}
%    i_j+1 & \text{$\boldsymbol{\sigma}_{j+1}$ is a direct homotopy letter}\\
%    i_j-1 & \text{$\boldsymbol{\sigma}_{j+1}$ is a inverse homotopy letter}\end{cases}\]

%We call the pair $(\boldsymbol{\sigma}, \kg_{\boldsymbol{\sigma}})$ where $\boldsymbol{\sigma}$ is a homotopy string (respectively band) a graded homotopy string (respectively band). 
Let $(\boldsymbol{\sigma},\kg_{\boldsymbol{\sigma}})$ be a graded homotopy string, where $\boldsymbol{\sigma}=\boldsymbol{\sigma}_1\boldsymbol{\sigma}_2...\boldsymbol{\sigma}_n$ and $\kg_{\boldsymbol{\sigma}}=(b_0,b_1,...,b_n)$. The inverse $(\boldsymbol{\sigma}^{-1},\kg_{\boldsymbol{\sigma}^{-1}})$ of $(\boldsymbol{\sigma}, \kg_{\boldsymbol{\sigma}})$ is a graded homotopy string where $\boldsymbol{\sigma}^{-1}=\overline{\boldsymbol{\sigma}_n}\overline{\boldsymbol{\sigma}_{n-1}}...\overline{\boldsymbol{\sigma}_1}$ and $\kg_{\boldsymbol{\sigma}^{-1}}=(b_n,b_{n-1},...,b_0)$. If $(\boldsymbol{\sigma},\kg_{\boldsymbol{\sigma}})$ a homotopy band then a cyclic rotation is of the form $\boldsymbol{\sigma}_i\boldsymbol{\sigma}_{i+1}...\boldsymbol{\sigma}_n\boldsymbol{\sigma}_1\boldsymbol{\sigma}_2...\boldsymbol{\sigma}_{i-1}$ with the grading $(b_{i-1},b_i,...,b_n,b_0,b_1,...,b_{i-2})$.

We need some equivalence relations for graded homotopy strings and bands, these are given in \cite[Section 2.4]{Arnesen}. 

\defn If $(\boldsymbol{\sigma},\kg_{\sigma})$ and $(\boldsymbol{\tau},\kg_{\boldsymbol{\tau}})$ are graded homotopy strings and $(\boldsymbol{\tau},\kg_{\boldsymbol{\tau}})$ is the inverse of $(\boldsymbol{\sigma},\kg_{\boldsymbol{\sigma}})$ then $(\boldsymbol{\sigma},\kg_{\boldsymbol{\sigma}})\sim^{str}(\boldsymbol{\tau},\kg_{\boldsymbol{\tau}})$.

\defn If $(\boldsymbol{\sigma},\kg_{\boldsymbol{\sigma}})$ and $(\boldsymbol{\tau},\kg_{\boldsymbol{\tau}})$ are graded homotopy bands and $(\boldsymbol{\tau},\kg_{\boldsymbol{\tau}})$ or $(\boldsymbol{\tau}^{-1},\kg_{\boldsymbol{\tau}^{-1}})$ is a cyclic rotation of $(\boldsymbol{\sigma},\kg_{\boldsymbol{\sigma}})$ then $(\boldsymbol{\sigma},\kg_{\boldsymbol{\sigma}})\sim^{bd}(\boldsymbol{\tau}, \kg_{\boldsymbol{\tau}})$.

As $\sD\cong\sK^b$ then \cite[Theorem 2.9]{Arnesen} gives the following classification of the indecomposable objects of $\sD$.

\thm There exists the isomorphism\[
\{\text{The indecomposable objects in } \sD\}\xleftrightarrow{1-1} \text{St}\sqcup\text{Ba},\] where $\text{St}=\{\text{Graded strings upto the equivalence } \sim^{str}\}$
and\\ $\text{Ba}= \{\text{Graded bands upto the equivalence } \sim^{bd}\}\times \sk^*\times \IN$.

By \cite[Proposition 1.1]{Arnesen} there exists a bijection between the set of permitted paths from a vertex $v$ to a vertex $u$ and a basis for $\Hom_\Lambda(P_u, P_v)$ where $P_u$ and $P_v$ are the projective $\Lambda$-modules corresponding to the vertices $u$ and $v$ respectively of $Q$. We shall identify a basis morphism in $\Hom_\Lambda(P_u, P_v)$ with a corresponding path from $v$ to $u$. % Let $\boldsymbol{\sigma}$ be a homotopy letter with source $s$ and target $t$. If $\boldsymbol{\sigma}$ is direct then $\sigma$ is a path from $s$ to $t$ and we have the map $\sigma:P_t\rightarrow P_s$, where $P_s$ is the projective $\Lambda$-module corresponding to the vertex $s$ in $Q$ and $P_t$ is the projective $\Lambda$-module corresponding to the vertex $t$ in $Q$. Otherwise, $\boldsymbol{\sigma}$ is inverse, so $\sigma$ is a path from $t$ to $s$ and we have the map $\sigma:P_s\rightarrow P_t$, where $P_s$ and $P_t$ are defined as above.

We now recall how an object in the derived category $\sD$ arises from a homotopy string or band.

\defn \cite[Definition 2.4]{Opper} We define the object $X\in\sD$ corresponding to a homotopy string/band $\boldsymbol{\sigma}_X=\boldsymbol{\sigma}_1\boldsymbol{\sigma}_2...\boldsymbol{\sigma}_n$ with a grading $\kg$. We denote $X$ by $(X^i, M^i)$ where $X^j$ is the object in the $j$th degree of $X$ and $M^k$ is the differential from $X^k$ to $X^{k-1}$. \begin{itemize}
    \item Let $\boldsymbol{\sigma}_X$ be a homotopy string with the grading $\kg=(b_0,b_1,...,b_n)$ and define $c_0=s(\boldsymbol{\sigma}_1)$ and $c_m=t(\boldsymbol{\sigma}_m)$ for $1\leq m\leq n$ \begin{itemize}
        \item The object in the $i$th degree is \[
        X^i=\bigoplus_{b_l=i}P_{c_l}
        \] where $P_{c_l}$ is the projective module corresponding to the vertex $c_l$.
        \item The differential from $X^{i}$ to $X^{i-1}$ is $M^i=(m^i)$ where \[
        m^i =\begin{cases}
            \sigma_u & \boldsymbol{\sigma}_u \text{ is direct and $m^i$ is the map from $P_{c_u}$ to $P_{c_{u-1}}$}\\
            \sigma_u & \boldsymbol{\sigma}_u \text{ is inverse and $m^i $ is the map from $ P_{c_{u-1}}$ to $P_{c_{u}}$}\\
            0 & \text{Otherwise}
            \end{cases}
        \]
    \end{itemize}
    \item Let $\boldsymbol{\sigma}_X$ be a homotopy band with the grading $\kg=(b_0,b_1,...,b_{n-1})$ and $\lambda$ be a scalar in $\sk^*$. Define $c_0=s(\boldsymbol{\sigma}_1)$ and $c_m=t(\boldsymbol{\sigma}_m)$ for $1\leq m\leq n-1$. \begin{itemize}
        \item The object $X^l$ is defined in the same way as for string objects.%\[
      %  X^l=\bigoplus_{b_t=l}P_{c_t} \] where $P_{c_t}$ is the projective module corresponding to the vertex $c_t$.
        \item For $n>2$ the morphism $M^i$ is defined in the same way as for string objects apart from if $m^i:P_{c_0}\rightarrow P_{c_{n-1}}$ and $\boldsymbol{\sigma}_n$ is direct then $m^i=\lambda\sigma_n$ or if $m^i:P_{c_{n-1}}\rightarrow P_{c_0}$ and $\boldsymbol{\sigma}_n$ is inverse then $m^i=\lambda\sigma_n$. For $n=2$ the only non-zero differential in the complex is $\lambda\sigma_1+\sigma_2$.
    \end{itemize} 
\end{itemize}

For homotopy bands there exists an object in $\sD$ for each scalar $\lambda\in\sk^*$ and natural number $r$. The $r$-dimensional band object $B_{\boldsymbol{\sigma}, \lambda, r}$, corresponding to the band $\boldsymbol{\sigma}$ and scalar $\lambda$ is defined in \cite[Section 5.1]{Arnesen}. %\begin{itemize}
%    \item The object in the $j$th degree is \[
%    B_r^j=(B^j)^{\oplus r}.\]
%    \item The differential from $B_r^k$ to $B_r^{k+1}$ is given by \[\begin{pmatrix}        D^k&A^k&&0\\0&D^k&&0\\&&\ddots&\\0&0&&A^k\\0&0&&D^k    \end{pmatrix}\] where $A^k=(a^k_{lm})$ and $a_{lm}^k$ is $\sigma_n$ if $d^k_{lm}=\lambda\sigma_n$ and zero otherwise.\end{itemize}

\subsection{The geometric model}

We can construct a ribbon graph from a gentle algebra $\Lambda$ from \cite[pp. 7]{Opper} and there is a ribbon surface $\Sigma_\Lambda$ defined in \cite[Definition 1.10]{Opper} corresponding to this ribbon graph. There exist a set of non-intersecting curves called a lamination $\kl_\Lambda$ of $\Sigma_\Lambda$ satisfying the conditions of \cite[Definition 1.14 \& Proposition 1.16]{Opper}. We define the geometric model constructed in \cite{Opper}. 

\defn The geometric model of $\sD^b(\text{mod }\Lambda)$ is provided by the ribbon surface $\Sigma_\Lambda$ and lamination $\kl_\Lambda$.

Objects in the derived category $\sD$ correspond to particular paths on $\Sigma_\Lambda$ and morphisms between two objects correspond to intersection points between their corresponding paths. We now consider some details of this.

\defn Let $\gamma:[0,1]\rightarrow\Sigma_\Lambda$ be a  path on $\Sigma_\Lambda$. \begin{itemize}
    \item The image of $\gamma$ is called an arc if $\gamma(0)$ and $\gamma(1)$ are marked points on the boundary of $\Sigma_\Lambda$ and these are the only points on the boundary of $\Sigma_\Lambda$.
    \item The image of $\gamma$ is called a closed curve if $\gamma(t)$ is not on the boundary of $\Sigma_\Lambda$ for any $t\in[0,1]$ and $\gamma(0)=\gamma(1)$.
\end{itemize}

For a path $\gamma$ we define its inverse $\gamma^{-1}:[0,1]\rightarrow\Sigma_\Lambda$ where $\gamma\inv(t)=\gamma(1-t)$. We use the notation $\gamma^\prime$ to denote a path homotopic to either $\gamma$ or its inverse.

A point $p=\gamma(t_1)=\psi(t_2)$ for some $t_1, t_2\in[0,1]$ is an intersection point of paths $\gamma$ and $\psi$.  If $\gamma=\psi$ then we additionally require $t_1\not=t_2$ and we call $p$ a self-intersection point. If $p$ is not a marked point then we call $p$ an internal intersection point. 

A collection of paths on $\Sigma_\Lambda$ is in minimal position if the number of intersection points between any two (possibly non-distinct) paths in the collection is minimal with respect to their homotopy classes. By \cite[Corollary 4]{Thurston}, we can assume a collection of paths can always be put into minimal position.

From \cite[Lemma 2.8]{Opper}, the image of a path is completely determined by the sequence of laminates it crosses in minimal position. It follows from \cite[Lemma 2.14]{Opper} that we can attach an homotopy letter $\boldsymbol{\sigma}(\ell, \ell^\prime)$ to each part of a path's image between two consecutive laminates $\ell$ and $\ell^\prime$ that it crosses. We attach a sequence of homotopy letters for a path $\gamma$ whose image crosses a sequence of laminates $(\ell_1, \ell_2,..,\ell_n)$\begin{itemize}
    \item If the image of $\gamma$ is an arc then it corresponds to the string $\boldsymbol{\sigma}_1\boldsymbol{\sigma}_2...\boldsymbol{\sigma}_{n-1}$, where the homotopy letter $\boldsymbol{\sigma}_i=\boldsymbol{\sigma}(\ell_i,\ell_{i+1})$ for all $1\leq i<n$.
    \item If the image of $\gamma$ is a closed curve then it corresponds to the string $\boldsymbol{\sigma}_1\boldsymbol{\sigma}_2...\boldsymbol{\sigma}_{n}$, where $\boldsymbol{\sigma}_i=\boldsymbol{\sigma}(\ell_i,\ell_{i+1})$ for $1\leq i<n$ and $\boldsymbol{\sigma}_n=\boldsymbol{\sigma}(\ell_n, \ell_1)$.
\end{itemize}

\defn The degree of a closed curve corresponding to a string $\boldsymbol{\sigma}$ is the number of direct homotopy letters in $\boldsymbol{\sigma}$ minus the number of inverse homotopy letters in $\boldsymbol{\sigma}$.

By \cite[Theorem 2.12]{Opper}, a string object corresponds to an arc on the geometric model. We call a closed curve graded if it has degree zero. A 1-dimensional band object corresponds to a graded closed curve. For a object $A\in\sD$ we denote its corresponding path as $\gamma_A$. 

\defn A subpath of a path $\gamma$ is a restriction of $\gamma$ to an interval $[t_1, t_2]\subseteq[0, 1]$.

Oriented intersection points between two paths are defined in \cite[Theorem 3.7]{Opper}. An internal intersection point between $\gamma_A$ and $\gamma_B$ is an oriented intersection point in each direction between $\gamma_A$ and $\gamma_B$. Under the correspondence in \cite[Theorem 3.3]{Opper}, such an intersection point between $\gamma_A$ and $\gamma_B$ corresponds to one morphism is each direction between $A$ and $B$ upto shift. Alternatively, an intersection point between $\gamma_A$ 
and $\gamma_B$ that is a marked point corresponds to one morphism in just one direction between $A$ and $B$ upto shift. Therefore, the geometric model gives a description of the objects and morphisms in $\sD^b(\text{mod }\Lambda)$.

In our work we will consider generators of the thick subcategories of $\sD$. If a thick subcategory has a string (respectively one-dimensional band) object generator then this corresponds to arc (respectively closed curve) on $\Sigma_\Lambda$. So the only generators that do not correspond paths on $\Sigma_\Lambda$ are $n$-dimensional band objects, where $n>1$. The next lemma implies we do not have to consider such generators.

\begin{lem} \label{bddim} Any thick subcategory generated by a $n$-dimensional band object $B_{w, \lambda, n}$ is generated by the one-dimensional band object $B_{w, \lambda, 1}$.\end{lem}

\pf If we let $B_j=B_{w, \lambda, j}$, from the band tubes mentioned in \cite[Section 5]{Arnesen} we have the triangles: \[
B_1\rightarrow B_2\rightarrow B_1\rightarrow B_1[1]\]
\[B_i\rightarrow B_{i-1}\oplus B_{i+1}\rightarrow B_i\rightarrow B_i[1]\]
for $i>1$. From these triangles and by the definition of thick subcategories we have the two filtrations:
\[\cdots\subseteq\td(B_n)\subseteq\cdots\subseteq\td(B_2)\subseteq\td(B_1)\]
\[\td(B_1)\subseteq\cdots\subseteq\td(B_{n-1})\subseteq\td(B_n)\subseteq\cdots\]
These filtrations give $\td(B_1)=\td(B_n)$.\qed 

By the definition of thick subcategories, we may assume (without loss of generality) that band objects in the generating set of a thick subcategory are one-dimensional.

\section{Preliminaries}
Let $\sD$ be a derived category of a gentle algebra $\Lambda$ corresponding to the geometric model $\Sigma$.
\subsection{Objects in the unbounded homotopy category}
In the previous section we recalled the paths on $\Sigma$ that correspond to objects in $\sD$. We would like to describe some objects in the unbounded homotopy category by paths on $\Sigma$.
First, we need some details on exceptional and spherelike objects in $\sD$. Let us start by recalling their definition.
 \defn \label{excsphobj} For an object $A\in\sD$ \begin{itemize}
    \item $A$ is exceptional if and only if $\dhomd(A,A)=1$;
    \item $A$ is spherelike if and only if $\dhomd(A,A)=2$.
\end{itemize}

The next result tells us about arcs on $\Sigma$ that correspond to spherelike and exceptional string objects in $\sD$.

\lem \label{excsph} Let $A$ be a string object in $\sD$.\begin{itemize}\item $A$ is exceptional if and only if $\gamma_A$ is an arc which is not closed and has no internal self-intersection points. \item $A$ is spherelike if and only if $\gamma_A$ is a closed arc with no internal self-intersection points. \end{itemize}

\pf Firstly, if $\gamma_A$ has an internal self-intersection point then by \cite[Theorems 3.3 \& 3.7]{Opper} this intersection point corresponds to two morphisms from $A$ to a shift of itself. Neither of these are the identity as from \cite[Remark 3.12]{Opper} the identity morphism in $\Hom_\sD(A,A)$ corresponds to an intersection point occurring at one of the marked points on $\Sigma$, so $\dhomd(A,A)\geq 3$.

Let $\gamma_A$ be an arc with no internal self-intersection points, any intersection between two lifts is a common end point, otherwise it would be the lift of an internal intersection point. If $\gamma_A$ is not closed, so $\gamma_A(0) \neq \gamma_A(1)$, then any common end-point of two lifts is of the form $\widehat{\gamma}_1(t) = \widehat{{\gamma}}_2(t)$ for some $t \in \{ 0,1\}$. However the unique lifting property \cite[Proposition 1.34]{Hatcher} then implies that the lifts $\widehat{\gamma}_1$ and $\widehat{{\gamma}}_2$ are the same (so distinct lifts do not intersect). Given two copies of the same lift, by \cite[Remark 3.12]{Opper} we can choose representations of them so they only intersect at both end-points. Precisely one of these will be graded intersection point oriented from $\widehat{\gamma}_1$ to $\widehat{{\gamma}}_2$. We apply \cite[Theorem 3.3]{Opper} to deduce that there is only one morphism in any degree from $A$ to $A$, that is, $\dhomd(A, A)=1$. 

If $\gamma_A$ was closed then $\gamma_A(0)=\gamma_A(1)=q$. Let $\hat\gamma_1$ be any lift of $\gamma_A$ with $\hat\gamma_1(0)=\hat q_1$ and $\hat \gamma_1(1)=\hat q_2$. There are precisely two lifts of $\gamma_A$, which are not homotopic to $\hat{\gamma}_1$ but which have a common point with $\hat\gamma_1$. We denote these by $\hat\gamma_2$ and $\hat\gamma_3$ such that $\hat\gamma_2(1)=\hat q_1$ and $\hat\gamma_3(0)=\hat q_2$. It follows from \cite[Theorem 3.7]{Opper} that one of these common end-points $\hat q_1$ or $\hat q_2$ is a graded intersection point oriented from $\hat\gamma_1$ but not both. Hence, we have a graded intersection from $\hat\gamma_1$ to $\hat{\gamma}_2$ or $\hat{\gamma}_3$ which does not correspond to the identity morphism 
$\text{id}_A$ by \cite[Remark 3.12]{Opper}. Since the only other lift that intersects $\hat\gamma_1$ is itself (upto homotopy) and as in the non-closed case we know the intersection points between two homotopic lifts give the identity morphism between the corresponding object and itself then $\dhomd(A,A)=2$.

Therefore $\dhomd(A, A)\geq 3$ if and only if $\gamma_A$ has an internal self-intersection point. If $A$ is an exceptional object then it follows from \cite[Theorem 3.3]{Opper} that $\gamma_A$ does not have any self-intersection points, so $A$ is exceptional if and only if $\gamma_A$ is a non-closed arc with no internal self-intersection points. We deduce that, $A$ is spherelike if and only if $\gamma_A$ is a closed arc with no internal intersection points. \qed

We now consider neighbourhoods around a point which we will use throughout this paper. For any point $p$ on $\Sigma$, there exists a open neighbourhood $\sU$ of $p$ and a continuous map $f:\sU\rightarrow\IR^2$. If $f(p)=(a,b)$ then we define a open neighbourhood $\kb_p$. If $p$ is a marked point on the boundary of $\Sigma$ then we can take $\kb_p$ to be homeomorphic to the half-ball
\[ \{(x,y)\in\IR^2\mid (x-a)^2+(y-b)^2<r\text{ and }y\geq b\}\] for some $r>0$ such that $\kb_p\subset f(\sU)$. Otherwise, $p$ is in the interior of $\Sigma$ and we can take $\kb_p$ to be homeomorphic to the ball \[\{(x,y)\in\IR^2\mid(x-a)^2+(y-b)^2<r\}\] for some $r>0$ such that $\kb_p\subset f(\sU)$. We let $\kn_p$ be the inverse image of $\kb_p$ under $f$ and $\overline{\kn_p}$ be the inverse image of the closure of $\kb_p$ under $f$. We abuse notation for the remainder of this paper by identifying $\kn_p$ and $\overline{\kn_p}$ with their images under $f$.

We also consider tubular neighbourhoods of paths, so we discuss some technicalities relating to these here. We begin by considering when we can take a tubular neighbourhood of a path. A path $\gamma:[t_1,t_2]\rightarrow\Sigma$ with no self-intersection points is an injective map. Therefore, the image of $\gamma$ is a submanifold (with boundary) $X$ of $\Sigma$ with an inclusion $i:X\hookrightarrow \Sigma$. By \cite[Theorem 6.5]{Cannas} there exists a tubular neighbourhood of $\gamma$. By \cite[Corollary 6.9]{Bott} any vector bundle over a contractible manifold is trivial, so the normal bundle to $\gamma$ is isomorphic to $[t_1, t_2]\times \IR$. The tubular neighbourhood of $\gamma$ is diffeomorphic to a subset of the normal bundle consisting of vectors of a particular length \cite[pp. 66]{Bott}, so we can take the 
tubular neighbourhood $\kt$ of $\gamma$ that contains all normal vectors to $\gamma$ with length $2\varepsilon$. This tubular neighbourhood is homeomorphic to $[t_1,t_2]\times[-\varepsilon,\varepsilon]$. A section in $\kt$ is a map from $[t_1, t_2]$ to $[t_1, t_2]\times [-\varepsilon, \varepsilon]$. By definition, $\gamma$ corresponds to the zero section given by $x\mapsto(x, 0)$.

We now consider paths, namely ungraded closed curves, that do not correspond to any object in $\sD$. We construct such paths from spherelike objects.

For a spherelike object $S$, $\gamma_S$ is a closed arc with no internal self-intersection points and end point $v$. There exists two subpaths $\restr{\gamma_S}{[0,t_1]}$ and $\restr{\gamma_S}{[t_2,1]}$ in $\overline{\kn_v}$ such that $\gamma_S(t_1)$ and $\gamma_S(t_2)$ are on the boundary of $\overline{\kn_v}$. When $\gamma_S$ is in minimal position, then we may choose the radius of $\overline{\kn_v}$ to be small enough such that $\gamma_S\cap\overline{\kn_v}=\restr{\gamma_S}{[0,t_1]}\cup\restr{\gamma_S}{[t_2,1]}$. For any $t_3<t_1<t_2<t_4$ there exists a tubular neighbourhood of $\restr{\gamma_S}{[t_3,t_4]}$ containing vectors with a length of $2\varepsilon$. There are two sections $x\mapsto(x, \lambda(x-t_3)(x-t_4))$ and $x\mapsto(x, \lambda(x-t_3)(t_4-x))$ for some $0\leq\lambda\leq\varepsilon$ which intersect the boundary of $\overline{\kn_v}$ at two points. We can choose $\sigma$ to be such a section that intersects the chord $\chi$ of $\kn_p$ between $\gamma_S(t_1)$ and $\gamma_S(t_2)$. As $\Sigma$ is an oriented surface, the paths $\sigma$ and $\chi$ intersect at two points $\sigma(u_1)=\chi(w_1)$ and $\sigma(u_2)=\chi(w_2)$ for some $0<w_1<w_2<1$ and $t_3<u_2<u_1<t_4$. We see this section in the diagram below.

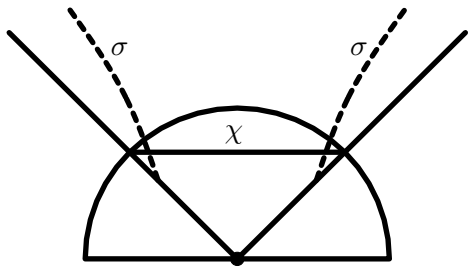
\begin{figure}[ht]
    \centering    
\definecolor{sqsqsq}{rgb}{0,0,0}
\begin{tikzpicture}[line cap=round,line join=round,>=triangle 45,x=2cm,y=2cm]
\draw [shift={(0,0)},line width=2pt,color=sqsqsq]  (0,0) --  plot[domain=0:3.141592653589793,variable=\t]({1*1*cos(\t r)+0*1*sin(\t r)},{0*1*cos(\t r)+1*1*sin(\t r)}) -- cycle ;
\draw [line width=2pt] (0,0)-- (-1.5,1.5);
\draw [line width=2pt] (0,0)-- (1.5,1.5);
\draw[line width=2pt,dashed, smooth,samples=100,domain=0:0.4630591222013562,variable=\t] plot({-0.24052067925304674*\t^3-0.40768783751385745*\t-0.5211349194290228},{-0.12113745452125312*\t^3+1.2150779399405016*\t+0.5211349194290227});
\draw[line width=2pt,dashed, smooth,samples=100,domain=0.4630591222013562:1,variable=\t] plot({0.20742562023366035*\t^3-0.6222768607009812*\t^2-0.11953686063144549*\t-0.5656118989012338},{0.10446923614065697*\t^3-0.31340770842197097*\t^2+1.360204238293518*\t+0.4987342339877959});
\draw[line width=2pt,dashed, smooth,samples=100,domain=0:0.46113380827557565, variable=\t] plot({0.2507163679155906*\t^3+0.40208584183361684*\t+0.52},{-0.1261872154388609*\t^3+1.2195453790784152*\t+0.52});
\draw[line width=2pt,dashed, smooth,samples=100,domain=0.46113380827557565:1,variable=\t] plot({-0.21455009668348507*\t^3+0.6436502900504553*\t^2+0.10527693238497161*\t+0.5656228742480581},{0.10798449059274134*\t^3-0.32395347177822403*\t^2+1.3689312772236018*\t+0.49703770396188096});
\draw (-0.9,1.5) node[anchor=north west] {$\sigma$};
\draw (0.675,1.5) node[anchor=north west] {$\sigma$};
\draw [line width=2pt] (-0.7071067811865476,0.7071067811865476)-- (0.7071067811865476,0.7071067811865476);
\draw (-0.15,0.95) node[anchor=north west] {$\chi$};
\begin{scriptsize}
\draw [fill=black] (0,0) circle (2.5pt);
\end{scriptsize}
\end{tikzpicture}
\end{figure}
We define the closed curve $\lambda_S$ homotopic to $\gamma_S$ as the path $\restr{\chi}{[w_1,w_2]}\restr{\sigma}{[u_2,u_1]}$. This path does not intersect $\gamma_S$ since $\chi$ and $\sigma$ only intersect $\gamma_S$ at their end points, which are not points on $\restr{\chi}{[w_1,w_2]}\restr{\sigma}{[u_2,u_1]}$.

Consider the closed curve $\lambda_S$. When it is in minimal position with respect to the lamination on $\Sigma$ it crosses a finite sequence of the laminates which we will denote by $(\ell_1, \ell_2,...,\ell_n)$. From \cite[Lemma 2.14]{Opper} there is a homotopy letter attached to each path starting at one laminate and ending at another, so we can attach the sequence of homotopy letters $\boldsymbol{\sigma}=\boldsymbol{\sigma}_1\boldsymbol{\sigma}_2...\boldsymbol{\sigma}_n$ to the closed curve $\lambda_S$. If $\boldsymbol{\sigma}$ is a band then there is a 1-dimensional band object defined in \cite[Definition 2.4(3)]{Opper} corresponding to it. Otherwise, we give $\boldsymbol{\sigma}$ the grading $\kg_{\boldsymbol{\sigma}}=(g_0,g_1,...,g_n)$ following \cite[Definition 2.3]{Opper}. We then turn $\boldsymbol{\sigma}$ into the infinite homotopy string $...\boldsymbol{\sigma}_n\boldsymbol{\sigma}_1\boldsymbol{\sigma}_2...\boldsymbol{\sigma}_n\boldsymbol{\sigma}_1...$ with the grading $(...,g_{-1},g_0,g_1,...,g_n,g_{n+1},...)$ satisflying the conditions of \cite[Definition 2.3]{Opper}. Now we follow \cite[Definition 2.4(1)]{Opper} to define a complex $P_S^\infty$ from the graded infinite homotopy string that we attach to $\lambda_S$. The object $P_S^\infty$ is in the unbounded homotopy category.

We now assume $\lambda_S$ does not correspond to a band object. As $S$ is spherelike then $\lambda_S$ does not has any self-intersection points. This implies $\lambda_S\not\simeq\omega^m$ for a closed curve $\omega$ and integer $m>1$. Therefore, if $\lambda_S$ does not correspond to a band object then it is ungraded, i.e. it has a non-zero degree. If $\lambda_S$ crosses the laminates $(\ell_1, \ell_2,...,\ell_j)$ then (a power of the arc corresponding to the spherelike object $S$) $\gamma_S^i$ crosses the same sequence of laminates $i-1$ times without crossing any other laminate in between. This implies that the string corresponding to $\gamma_S^i$ is of the form $\boldsymbol{\beta}_1\boldsymbol{\sigma}^{i-1}\boldsymbol{\beta}_2$, where $\boldsymbol{\beta}_1$, $\boldsymbol{\beta}_2$ are strings. By constuction, there exists $b>0$ such that the complex $P_S^\infty$ is periodic with the fundamental domain $C_1\rightarrow C_2\cdots \rightarrow C_{b-1}\rightarrow C_b$, since $\lambda_S$ has a non-zero degree. Define $S_k$ to be the object corresponding to $\gamma_S^k$.  Therefore, there exists $i_1>0$ such that there is one copy of the fundamental domain of $P_S^\infty$ in the complex of $S_{i_1}$ i.e. $S_{i_1}$ has the form
\[\cdots 0\rightarrow L_1\rightarrow L_2\rightarrow\cdots\rightarrow L_{l}\rightarrow C_1\rightarrow C_2\rightarrow\cdots\rightarrow C_b\rightarrow R_{1}\rightarrow R_2\rightarrow\cdots\rightarrow R_m\rightarrow 0\cdots\]

Consider the unfolded diagram of $S_{i_1+i_2}$ for $i_2\in\IN$

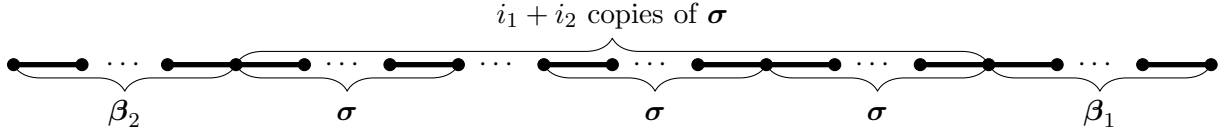
\begin{center}
    \begin{tikzpicture}[line cap=round,line join=round,>=triangle 45,x=1cm,y=1cm, scale=0.9]
\draw [line width=2pt] (0,0)-- (1,0);
\draw [line width=2pt] (2.25,0)-- (4.25,0);
    \draw [decorate,decoration={brace,amplitude=10pt},xshift=0pt,yshift=0pt]
    (3.25,0) -- (0,0) node [black,midway,xshift=0cm,yshift=-0.65cm] 
    {$\boldsymbol{\beta}_2$};
    \draw [decorate,decoration={brace,amplitude=10pt},xshift=0pt,yshift=0pt]
    (6.5,0) -- (3.25,0) node [black,midway,xshift=0cm,yshift=-0.65cm] 
    {$\boldsymbol{\sigma}$};
    \draw [decorate,decoration={brace,amplitude=10pt},xshift=0pt,yshift=0pt]
    (11,0) -- (7.75,0) node [black,midway,xshift=0cm,yshift=-0.65cm] 
    {$\boldsymbol{\sigma}$};
    \draw [decorate,decoration={brace,amplitude=10pt},xshift=0pt,yshift=0pt]
    (14.25,0) -- (11,0) node [black,midway,xshift=0cm,yshift=-0.65cm] 
    {$\boldsymbol{\sigma}$};
    \draw [decorate,decoration={brace,amplitude=10pt},xshift=0pt,yshift=0pt]
    (17.5,0) -- (14.25,0) node [black,midway,xshift=0cm,yshift=-0.65cm] 
    {$\boldsymbol{\beta}_1$};
    \draw [decorate,decoration={brace,amplitude=10pt},xshift=0pt,yshift=0pt]
    (3.25,0) -- (14.25,0) node [black,midway,xshift=0cm,yshift=0.65cm] 
    {$i_1+i_2$ copies of $\boldsymbol{\sigma}$};
\draw (1.2,0.25) node[anchor=north west] {$\cdots$};
\draw (4.4,0.25) node[anchor=north west] {$\cdots$};
\draw (6.65,0.25) node[anchor=north west] {$\cdots$};
\draw (8.9,0.25) node[anchor=north west] {$\cdots$};
\draw (12.15,0.25) node[anchor=north west] {$\cdots$};
\draw (15.4,0.25) node[anchor=north west] {$\cdots$};
\draw [line width=2pt] (5.5,0)-- (6.5,0);
\draw [line width=2pt] (7.75,0)-- (8.75,0);
\draw [line width=2pt] (10,0)-- (12,0);
\draw [line width=2pt] (13.25,0)-- (15.25,0);
\draw [line width=2pt] (16.5,0)-- (17.5,0);
\begin{scriptsize}
\draw [fill=black] (0,0) circle (2.5pt);
\draw [fill=black] (1,0) circle (2.5pt);
\draw [fill=black] (2.25,0) circle (2.5pt);
\draw [fill=black] (3.25,0) circle (2.5pt);
\draw [fill=black] (4.25,0) circle (2.5pt);
\draw [fill=black] (5.5,0) circle (2.5pt);
\draw [fill=black] (6.5,0) circle (2.5pt);
\draw [fill=black] (7.75,0) circle (2.5pt);
\draw [fill=black] (8.75,0) circle (2.5pt);
\draw [fill=black] (11,0) circle (2.5pt);
\draw [fill=black] (12,0) circle (2.5pt);
\draw [fill=black] (10,0) circle (2.5pt);
\draw [fill=black] (13.25,0) circle (2.5pt);
\draw [fill=black] (14.25,0) circle (2.5pt);
\draw [fill=black] (15.25,0) circle (2.5pt);
\draw [fill=black] (16.5,0) circle (2.5pt);
\draw [fill=black] (17.5,0) circle (2.5pt);
\end{scriptsize}
\end{tikzpicture}
\end{center}

We label the copies of $\boldsymbol{\sigma}$ from the right in this unfolded diagram. We take the shift of $S_{i_1}$ so the last $C_b$ before $R_{1}$ in the complex $S_{i_1}$ is in degree zero. This means that the $1$st to the $i_1$th copies of $\boldsymbol{\sigma}$ gives a copy of the fundamental domain of $P^\infty_S$ with the last $C_b$ in degree zero. If $d$ is the degree of $\boldsymbol{\sigma}$ then for $c\in\IN$ the $c$th to $(c+i_1)$th copies of $\boldsymbol{\sigma}$ give a copy of the fundamental domain of $P_S^\infty$ with the last $C_b$ in the $((c-1)d)$th degree. This proves the following lemma.
\begin{lem}
    \label{remAjobj}
    There exists $i_1\in\IN$ such that for any $i_2\in\IN$, the complex $S_{i_1+i_2}$ has the form \begin{equation}\label{Si}
        \cdots 0\rightarrow L_1\rightarrow L_2\rightarrow\cdots \rightarrow L_{l}\rightarrow \underbrace{C_1\rightarrow C_2\rightarrow\cdots \rightarrow C_b}_{i_2+1 \text{ copies}}\rightarrow R_{1}\rightarrow R_2\rightarrow\cdots\rightarrow R_m\rightarrow 0\cdots.
    \end{equation}
\end{lem} 

Let $\sK$ be the unbounded homotopy category of projective modules over the gentle algebra $\Lambda$. For two objects $A,B\in\sK$ we define $\hom_{\sK/\Sigma}(A,B)$ to be $\infty$ if  $A$ or $B$ is unbounded and $\hom_\sK(A,B[n])\not=0$ for some $n\in\IZ$. Otherwise, \[
\hom_{\sK/\Sigma}(A,B):=\sum_{i\in\IZ} \hom_\sK(A,B[i]).\]

\begin{lem}\label{hominfAi}
    Assume $S$ is a spherelike object, $P$ is a string object and $\hom_{\sK/\Sigma}(P_S^\infty,P)=~0$ then $\dhomd(S_i,P)=\dhomd(S_{i+1},P)$ for all $i\gg1$ 
\end{lem}

\pf $P$ is a indecomposable bounded complex, so it only has non-zero modules in a finite number of degrees and these are all in consecutive degrees, i.e. \[
P:\cdots 0\rightarrow B_1\rightarrow B_2\rightarrow\cdots\rightarrow B_{a-1}\rightarrow B_a\rightarrow0\cdots
\] for some $a\in\IN$.

By Lemma~\ref{remAjobj}, we know there exists $i_1>1$ such that for any $i_2>0$, $S_{i_1+i_2}$ takes the form of (\ref{Si}) from Lemma~\ref{remAjobj}. %\[\cdots 0\rightarrow R_1\rightarrow R_2\cdots R_{m}\rightarrow C_{i-1}\rightarrow \underbrace{C_i\cdots C_b\rightarrow C_1\cdots C_{i-1}}_{j_2+1 \text{ copies}}\rightarrow C_i\rightarrow R_{m+1}\cdots R_l\rightarrow 0\cdots    \]

We can let $i_2$ be the minimum natural number such that $i_2b>a$. Any non-zero morphism from $S_{i_1+i_3}$ (for any $i_3>i_2$) to a shift of $P$ is a chain map between their complexes in $\sK^b=\sK^b(\text{proj }\Lambda)$ (we have identified $\sD$ with $\sK^b$ \cite[Section 1]{Arnesen}).

Now assume $\hom_{\sK/\Sigma}(P_S^\infty,P)=0$. As $S_{i_1+i_3}$ and $P$ are bounded complexes, $\hom_\sD(S_{i_1+i_3},P[u])$ can only be non-zero if the degrees where the complexes are supported overlap, and this overlap is not in contained in the periodic part of $S_{i_1+i_3}$ as this would induce a morphism from $P_S^\infty$ to the shift $P[u]$ of $P$. Also, by the construction of $S_{i_1+i_3}$, the supported part of $P[u]$ cannot overlap $S_{i_1+i_3}$ in the degrees of both $L_l$ and $R_1$. If $B_1$ in $P$ and $L_1$ in $S_{i_1+i_3}$ are both in degree zero, then this happens when $u\in\ki_1\cup\ki_2$ where $\ki_1=\{-l,-l+1,...,0,1,...,a-1\}$ and $\ki_2=\{a-l-(i_3+1)b,a-l-(i_3+1)b-1...,1-m-l-(i_3+1)b\}$. 

Let $x=i_1+i_3$. There exists a map \[
\Phi_{u^\prime}:\Hom_\sD(S_x,P[u^\prime])\rightarrow\Hom_\sD(S_{x+1},P[u^*])
\] where \[
u^*=\begin{cases}
    u^\prime -b& \text{if } u^\prime\in\ki_2,\\
    u^\prime & \text{otherwise.}
\end{cases} \] A map $f\in\Hom_\sD(S_x,P[u^\prime])$ is mapped under $\Phi_{u^\prime}$ to the map in $\Hom_\sD(S_{x+1}, P[u^*])$ which is a chain map with the same sequence of non-zero homomorphisms as $f$. This is a well-defined map as any homotopy between $f$ and $g$ in $\Hom_\sD(S_x,P[u^\prime])$ induces an homotopy between $\Phi_{u^\prime}(f)$ and $\Phi_{u^\prime}(g)$. Now define the map \[
\Phi:\bigoplus_{\widetilde{u}\in\IZ}\Hom_\sD(S_x,P[\widetilde{u}])\rightarrow\bigoplus_{\overline{u}\in\IZ}\Hom_\sD(S_{x+1},P[\overline{u}])
\]where for $h\in\Hom_\sD(S_x,P[\widetilde{u}])$, $\Phi(h)=\Phi_{\widetilde{u}}(h)$.

When doing the following calculations, we consider morphisms in indecomposable summands of $\bigoplus_{\widetilde{u}\in\IZ}\Hom_\sD(S_x,P[\widetilde{u}])$ and $\bigoplus_{\overline{u}\in\IZ}\Hom_\sD(S_{x+1},P[\overline{u}])$ because for any $v_1\in\IZ$, if the image of one map in $\Hom_\sD(S_x,P[v_1])$ under $\Phi$ is in $\Hom_\sD(S_{x+1},P[v_2])$ for some $v_2\in\IZ$ then all maps in $\Hom_\sD(S_x,P[v_1])$ have an image under $\Phi$ in $\Hom_\sD(S_{x+1},P[v_2])$.

\textit{$\Phi$ is well defined:} Let $f,g\in\Hom_D(S_x,P[u_1])$ for some $u_1\in\IZ$ be two homotopic maps. As $\Phi_{u_1}$ is well-defined then $\Phi(f)\simeq\Phi(g)$ and $\Phi$ is a well-defined map.

\textit{$\Phi$ is injective:} Let $\Phi(f)$ and $\Phi(g)$ be homotopic maps, so they are both in $\Hom_\sD(S_{x+1},P[u_2])$ for some $u_2\in\IZ$. If $\Phi(f)$ and $\Phi(g)$ are homotopic to zero then $f$ and $g$ are homotopic to zero. Therefore we need to check if $f$ and $g$ can be homotopic when $\Phi(f)$ and $\Phi(g)$ are not homotopic to zero. If $f\in\Hom_\sD(S_x,P[u_3])$ and $g\in\Hom_\sD(S_x,P[u_4])$ then $u_3,u_4\in\{u_2,u_2+b\}$, otherwise they would not map to morphisms in $\Hom_\sD(S_{x+1},P[u_2])$ by the definition of $\Phi$. We have the following cases, of which only one can happen.

\textit{Case 1:} $u_3=u_4$ and the chain homotopy between $\Phi(f)$ and $\Phi(g)$ induces a chain homotopy between $f$ and $g$.

\textit{Case 2:} $u_3=u_2$ and $u_4=u_2+b$. Since if $f\simeq0$, it follows that $\Phi(f)\simeq0$, so $f\not\simeq0$ and from the definition of $\Phi$, $u_3\in\ki_1$. The definition of $\Phi$ also says that $u_4\in\ki_2$. This implies $f$ and $\Phi(f)$ have a non-zero component to some $L_{l^\prime}$ for $1\leq l^\prime\leq l$ or the first $C_1$ from the left but no map to any $R_{m^\prime}$ for $1\leq m^\prime\leq m$. Similarly, $g$ and $\Phi(g)$ have a non-zero component to some $R_{m^{\dprime}}$ for $1\leq m^{\dprime}\leq m$ or the last $C_b$ from the left but no map to any $L_{l^{\dprime}}$ for $1\leq l^{\dprime}\leq l$. Therefore, $\Phi(f)$ and $\Phi(g)$ must be maps to different shifts of $P$, which is a contradiction that they are homotopic.

\textit{Case 3:} $u_3=u_2+b$ and $u_4=u_2$. By a similar argument to Case 2, we cannot have this case either.

Therefore, if $\Phi(f)\simeq\Phi(g)$ then $f\simeq g$ and $\Phi$ is injective

\textit{$\Phi$ is surjective:} Let $h$ be a morphism in $\Hom_\sD(S_{x+1},P[u_5])$. If $u_5\notin\ki_1\cup\ki_2$ then $h\simeq0$ and there exists zero map in $\Hom_\sD(S_x,P[u_5])$ which maps to $h$ under $\Phi$. If $u_5\in\ki_1$ then $h$ induces a morphism $h^\prime\in\Hom_\sD(S_x,P[u_5])$ and $\Phi(h^\prime)=h$, If $u_5\in\ki_2$ then $h$ induces the morphism $h^{\dprime}\in\Hom_\sD(S_x,P[u_5+b])$ and $\Phi(h^{\dprime})=h$. Hence, $\Phi$ is surjective.

Therefore, $\Phi$ is an well-defined bijection. We can write the above argument for any $i_3>i_2$ and the result follows.\qed

For a spherelike object $S$ consider $\gamma_S$ and $\lambda_S$. They are homotopic closed paths that do not self-intersect or intersect each other and so by \cite[Lemma 2.4]{Thurston} there exists a cylinder $\kC$ bounded by them.

\lem \label{cycsa} Assume $S$ is a spherelike object and $P$ is an exceptional or spherelike string object which does not equal $S$. Suppose that $\gamma_P$ does not internally intersect $\gamma_S$ when in minimal position. Then the intersection of $\gamma_P$ and the cylinder $\kC$ bounded by $\gamma_S$ and $\lambda_S$ is either empty or a union of subpaths of $\gamma_P$ which each have an end point on $\lambda_S$ and another which is the end point $q$ of $\gamma_S$.

\pf Assume $\gamma_S$, $\gamma_P$ and $\lambda_S$ are in minimal position with respect to each other. If $\gamma_P$ does not have a subpath on $\kC$ then $\kC\cap\gamma_P =\emptyset$. We assume otherwise and let $\chi$ be a subpath of $\gamma_P$ on $\kC$. We have the following possible cases since we have assumed $\gamma_P$ does not intersect $\gamma_S$ when in minimal position. The third of these cases is the only one that can be true.

\textit{Case 1:} $\chi$ has both end points at $q$, so $\gamma_P$ is entirely on $\kC$. This would imply that $\gamma_P$ is homotopic to $\gamma^{\pm1}_S$, since $\gamma_P$ has no internal self-intersection points and the fundamental group of a cylinder is isomorphic to $\IZ$. This contradicts our assumption that $P\not=S$ and therefore this case is not possible.

\textit{Case 2:} $\chi$ has both end points on $\lambda_S$ which implies $\gamma_P$ forms a bigon with $\lambda_S$. This implies $\gamma_P$ and $\lambda_S$ are not in minimal position with respect to each other \cite[Lemma 10]{Thurston}, which would contradict our assumption. 

\textit{Case 3:} $\chi$ has end point of $q$ and another on $\lambda_S$. This is the only possible case for $\gamma_P$ to have a subpath on $\kC$. \qed

If $\boldsymbol{\sigma}$ is a homotopy band associated to a closed curve $\lambda_S$ corresponds to a band object $B_{\boldsymbol{\sigma}}$ then by \cite[Theorem 3.3]{Opper} $\lambda_S$ intersects an arc $\gamma_P$ if and only if there exists morphisms between $P$ and $B_{\boldsymbol{\sigma}}$. We now prove an analogous result when $\lambda_S$ corresponds to a unbounded complex.

\begin{lem} \label{hominf}
   Let $S$ be a spherelike string object in $\sD$ such that $\lambda_S$ does not correspond to a band. Let $P\not=S$ be a non-zero spherelike or exceptional object in $\sD$. Suppose $\gamma_S$ and $\gamma_P$ do not internally intersect. Then the following are equivalent: \begin{enumerate}[i)]
%\item $\gamma_P$ internally intersects $\alpha^k$ for some $k>1$ (in minimal position)
\item $\hom_{\sK/\Sigma}(P_S^\infty, P)\not=0$
\item $\gamma_P$ intersects $\lambda_S$ (in minimal position).\end{enumerate}\end{lem}

\pf (i) $\implies$ (ii) Assume that $\hom_{\sK/\Sigma}(P_S^\infty,P)$ is non-zero.

We have the objects \[
P_S^\infty:\cdots C_b\rightarrow C_1\rightarrow C_2\cdots \rightarrow C_{b-1}\rightarrow C_b \rightarrow C_1\cdots
\]\[
P:\cdots 0\rightarrow Q_1\rightarrow Q_2\rightarrow\cdots\rightarrow Q_{a-1}\rightarrow Q_a\rightarrow0\cdots
\]

There exist a morphism $f$ from $P^\infty_S$ to a shift of $P$ given by the chain map which is not homotopic to zero\[
\begin{tikzcd}[cells={nodes={minimum height=2em}}]
P^\infty_S:\cdots D_m\arrow[r] & D_1 \arrow[r] \arrow[d, "f_1"] &D_2 \arrow[r] \arrow[d, "f_2"] &  \cdots\arrow[r] &  D_a \arrow[r] \arrow[d, "f_a"] &D_{a+1}\arrow[r]\cdots& D_m \cdots\\
P[u_1]:\cdots 0\arrow[r] & Q_1 \arrow[r]  &Q_2 \arrow[r]   & \cdots \arrow[r]&  Q_a \arrow[r]  & 0 \cdots\arrow[r]&0\cdots
\end{tikzcd}
\]
where $m=sb$ such that $s\in\IN$ is the minimal natural number with $sb>a$. The object $D_1=C_i$ for some $1\leq i\leq b$. By Lemma~\ref{remAjobj} there exist $j>1$ such that $S_{j+5s-1}$ is a complex with $5s$ copies of the fundamental domain of $P_S^\infty$ \[
\cdots 0\rightarrow R_1\rightarrow R_2\cdots R_{n}\rightarrow C_{i-1}\rightarrow \underbrace{C_i\cdots C_b\rightarrow C_1\cdots C_{i-1}}_{5s \text{ copies}}\rightarrow C_i\rightarrow R_{n+1}\cdots R_l\rightarrow 0\cdots
\]

which is equivalent to \[
\cdots 0\rightarrow R_1\rightarrow R_2\cdots R_{n}\rightarrow C_{i-1}\rightarrow \underbrace{D_1\rightarrow D_2\cdots D_{m}}_{5 \text{ copies}}\rightarrow C_i\rightarrow R_{n+1}\cdots R_l\rightarrow 0\cdots\]
Let $Q_1$ in $P[u_1]$ and $R_1$ in $S_{j+5s-1}$ both be in degree zero. The family of maps $\{f_c\mid 1\leq c\leq a\}$ induce a chain map from $S_{j+5s-1}$ to $P[u_1-u_2]$ for each $u_2\in\{mt+n+1\mid0\leq t\leq 4\}$. We see the form of each of these maps below: 
\[\begin{tikzcd}[cells={nodes={minimum height=2em}}]
S_{j+5s-1}:\cdots D_m\arrow[r] & D_1 \arrow[r] \arrow[d, "f_1"] &D_2 \arrow[r] \arrow[d, "f_2"] &  \cdots\arrow[r] &  D_a \arrow[r] \arrow[d, "f_a"] &D_{a+1}\cdots\\
P[u_1-u_2]:\cdots 0\arrow[r] & Q_1 \arrow[r]  &Q_2 \arrow[r]   & \cdots \arrow[r]&  Q_a \arrow[r]  & 0 \cdots
\end{tikzcd}
\]

These maps are not homotopic to each other as they are to different shifts 
of $P$. A homotopy between any of them and the zero map would induce a homotopy between $f$ and the zero map. This would be a contradiction of $f$ not being homotopic to zero. Hence, $\dhomd(S_{j+5s-1}, P)\geq5$ and $\gamma_P$ must internally intersect $\gamma_S^{j+5s-1}$. Otherwise the two arcs can only have intersections at common end points, as $\gamma_S$ has both end points at $q$ and $\gamma_P$ has at most two end points at $q$ then this would imply $\dhomd(S_{j+5s-1}, P)\leq4$, which is a contradiction.

%(i) $\implies$ (ii) and (iii) Assume $\gamma_P$ internally intersects $\gamma_S^k$ for some $k>1$. 

Now consider $\gamma_S$ and $\lambda_S$. We have the cylinder $\kC$ bounded by them. We may assume by \cite[Corollary 4]{Epstein} that $\gamma_P$ is in minimal position with respect to $\gamma_S$ and $\lambda_S$. For each $k>1$, there is a representative of $\gamma_S^k$ which is contained in $\kC$. Since $\gamma_P$ internally intersects $\gamma_S^{j+5s-1}$ then $\gamma_P$ must have a subpath on $\kC$. By Lemma~\ref{cycsa}, this subpath of $\gamma_P$ on $\kC$ has a end point on $\lambda_S$, so (ii) holds.

(ii) $\implies$ (i) Suppose $\gamma_P$ intersects $\lambda_S$ when in minimal position. $\gamma_P$ must have a subpath $\phi$ on the cylinder $\kC$ bounded by $\lambda_S$ and $\gamma_S$.

By Lemma~\ref{cycsa}, there exists a subpath $\phi$ of $\gamma_P$ with the end points $r$ on $\lambda_S$ and $q$ the end point of $\gamma_S$. Consider the universal cover $\widehat{\kC}$ of $\kC$ which is homeomorphic to $\IR\times[0,1]$. For any $n_1\in\IZ$ we denote the lift of a point $c$ in $[n_1, n_1+1)\times[0,1]$ by $(c, n_1)$. The lifts of $\gamma_S^k$ have the end points $(q, n_2)$ and $(q, n_2+k)$ for some $n_2\in\IZ$ and the lifts of $\phi$ have the end points $(q, n_3)$ and $(r, n_4)$ for some $n_3,n_4\in\IZ$, we denote this lift of  $\phi$ by $\widehat{\phi}_{n_3}$. By \cite[Lemma 3.8]{Broomhead} the lift $\widehat{\gamma_S^k}$  of $\gamma_S^k$ with the end points $(q,0)$ and $(q,k)$ is intersected once by the each of the lifts $\widehat{\phi}_0, \widehat{\phi}_1,...,\widehat{\phi}_k$. The lifts $\widehat{\phi}_0$ and $\widehat{\phi}_k$ only intersect $\widehat{\gamma_S^k}$ at their common end point with $\widehat{\gamma_S^k}$, which are both lifts of the same marked point $q$. It follows from \cite[Corollary 3.6]{Broomhead} that $\gamma_S^k$ has $k$ intersection points with $\phi$. Only one of these intersection points is a end point of $\gamma_S^k$, therefore $\gamma_S^k$ and $\phi$ have $k-1$ internal intersection points.

%and only one of these is a graded oriented intersection point from $\widehat{\gamma_S^k}$. The other lifts must internally intersect $\widehat{\gamma_S^k}$, since they do not share any common end points with this lift of $\gamma_S^k$. Therefore there are $k$ lifts of $\phi$ that have a oriented intersection point from $\widehat{\gamma_S^k}$.

The arc $\gamma_P$ can only have at most two subpaths on $\kC$ when in minimal position with respect to $\gamma_S$ and $\lambda_S$, since any such subpath must have an end point which is the end point $q$ of $\gamma_P$. We have the two following cases:

\textit{Case 1:} There is only one subpath $\phi_1$ of $\gamma_P\simeq\phi_1^\prime\psi^\prime$ on $\kC$. In minimal position $\phi_1$ has $k-1$ graded oriented internal intersection points from $\gamma_S^k$. Looking locally around the point $q$ we see that the number of morphisms from $S_k$ to $P$ corresponding to $q$ is $\dhomd(S_1,P)$. Hence $\dhomd(S_k, P)=k-1+\dhomd(S_1,P)$. Otherwise, $\psi$ would be a subpath that formed part of a bigon between $\gamma_P$ and $\gamma_S^k$. This would imply that $\gamma_P$ is homotopic to a path entirely on $\kC$. As $\gamma_P$ has no internal intersection points (when in minimal position) and the fundamental group of $\kC$ is isomorphic to $\IZ$, then this means that either $P=0$ or $P=S_1$, which is a contradiction of our assumptions about $P$.

\textit{Case 2:} There are two subpaths $\phi_1$ and $\phi_2$ of $\gamma_P$ on $\kC$ with end points $r_1$ and $r_2$ respectively on $\lambda_S$. The arc $\gamma_P$ is homotopic to $\phi_1\psi\phi_2$, where $\psi$ is the subpath of $\gamma_P$ outside of $\kC$. We can set the paths $\phi_1$, $\phi_2$ and $\gamma_S^k$ to be in minimal position on $\kC$, so $\phi_1$ and $\phi_2$ each intersect $\gamma_S^k$ $k$ times. If there are no bigons between $\gamma_S^k$ and $\gamma_P$ then they are in minimal position by \cite[Lemma 10]{Thurston}. If such a bigon between $\gamma_S^k$ and $\gamma_P$ did exist then it would be formed by a subpath of $\gamma_S^k$ and a subpath $\alpha$ of $\gamma_P$ that had one end $\gamma_P(w_1)$ on $\phi_1$ and the other $\gamma_P(w_2)$ on $\phi_2$. If $r_1=\gamma_P(w_3)$ and $r_2=\gamma_P(w_4)$ then (without loss of generality) we could assume $\gamma_P$ is oriented such that $w_1<w_3<w_4<w_2$. Hence, $\psi$ would be a subpath of $\alpha$, which would imply $\psi$ forms a bigon with $\lambda_S$. This contradicts the fact that the two simple closed curves $\gamma_P$ and $\lambda_S$ are in minimal position with respect to each other by \cite[Proposition 1.7]{Farb}.  Therefore, the paths $\gamma_S^k$ and $\gamma_P$ form no bigons between them and are in minimal position. Hence, $\gamma_S^k$ has $2k$ oriented intersection points to $\gamma_P$ when they are in minimal position and $\dhomd(S_k, P)=2k$.

The arc $\gamma_S^y$ and a subpath $\phi_3$ of $\gamma_P$ on $\kC$ have $y$ intersection points between them if and only if there are $y+1$ intersection points between $\gamma_S^{y+1}$ and $\phi_3$. Hence, if $\gamma_P$ has one subpath on $\kC$ then $\dhomd(S_y,P)=y-1+\dhomd(S_1,P)$ and $\dhomd(S_{y+1},P)=y+\dhomd(S_1,P)$. Otherwise, $\gamma_P$ has two subpaths on $\kC$, so $\dhomd(S_y,P)=2y$ and $\dhomd(S_{y+1},P)=2(y+1)$. Therefore,  $\dhomd(S_y,P)<\dhomd(S_{y+1},P)$ for any $y\geq 1$. By the contrapositive of Lemma~\ref{hominfAi} this implies $\hom_{\sK/\Sigma}(P^\infty_S,P)~\neq~0$. \qed

\subsection{Morphisms between thick subcategories}

Now we see when thick subcategories are semi-orthogonal or fully orthogonal by considering morphisms between their generating sets.

First we introduce some notation. The objects in a thick subcategory are generated in a finite number of steps by taking shifts, direct summands and cones of morphisms between objects generated in fewer steps. Note that, the objects generated in zero steps are precisely the ones in generating set of the thick subcategory.

\defn For a thick subcategory $\sT$, we denote the set of all shifts of direct summands of objects generated in $n$ or less steps by $\sT_{\leq n}$.

\begin{lem} \label{zhom} Let $\sD$ be a full triangulated subcategory of a triangulated category $\sC$ and $\sT$ be a thick subcategory generated by a set of objects $\kS$ in $\sD$. If $X\in\sC$ is such that $\hom_{\sC/\Sigma}(\kS,X)=0$ then $\hom_{\sC/\Sigma}(\sT,X)=0$. Furthermore, if $\sT^\prime$ is another thick subcategory of $\sD$ generated by a set of objects $\kS^\prime$ with $\dhomd(\kS,\kS^\prime)=0$ then $\dhomd(\sT,\sT^\prime)=0$. \end{lem}

\pf  Since $\sD$ is closed under taking shifts, finite extensions and direct summands in $\sC$ then so is $\sT$. Any object in a thick subcategory $\sT^*$ of $\sD$ is in $\sT^*_{\leq l}$ for some $l\in \IZ$. 

%We define $\sT_{\leq n}$ to be the collection of objects generated in $n$ steps or fewer, where the objects in $\kS$ are precisely the objects in $\sT_{\leq0}$.

Let $X\in\sC$ be an object such that $\hom_{\sC/\Sigma}(\sT_{\leq0},X)=0$.

We do a proof by induction. Let $i\geq 0$ and suppose $\hom_{\sC/\Sigma}(\sT_{\leq i},X)=0$.

We prove that $\hom_{\sC/\Sigma}(\sT_{\leq i+1},X)=0$. If $Y\in\sT_{\leq i+1}$ is a shift of a direct summand of an object in $\sT_{\leq i}$ then  by our assumption $\hom_{\sC/\Sigma}(Y,X)=0$.
%$Y=Z[k_1]$ and $Z\in\sT_{\leq i}$ then a morphism from $Y$ to $X[k_2]$ would induce a morphism from $Z$ to $X[k_2-k_1]$. Our assumption says $\hom_{\sC/\Sigma}(Z,X)=0$, so $\hom_{\sC/\Sigma}(Y,X)=0$.

%If $Y$ is a direct summand of an object $Z\in\sT_{\leq i}$ then any morphism from $Y$ to $X[k_3]$ induces a morphism from $Z$ to $X[k_3]$. As $\hom_{\sC/\Sigma}(Z,X)=0$ then there can be no morphisms from $Y$ to any shift of $X$.

If $Y$ is the cone of a morphism from $Z_1$ to $Z_2$, where $Z_1$ and $Z_2$ are in $\sT_{\leq i}$ then by the definition of cones it fits into the triangle $Z_1\rightarrow Z_2\rightarrow Y\rightarrow Z_1[1]$. For any $k\in\IZ$ there exist part of a long exact sequence \[
\cdots\rightarrow\Hom_\sC(Z_1, X[k-1])\rightarrow\Hom_\sC(Y,X[k])\rightarrow\Hom_\sC (Z_2, X[k])\rightarrow\cdots
\]
Since both $Z_1$ and $Z_2$ have no morphism to any shift of $X$ and this is a long exact sequence then $\hom_\sC(Y, X[k])=0$ for any $k\in\IZ$. Since any morphism from 
a direct summand of $Y$ to $X[k]$ would induce a morphism from $Y$ to $X[k]$ and $\hom_\sC(Y,X[k])=0$ then for any direct summand $Y^\prime$ of $Y$, $\hom_{\sC/\Sigma}(Y^\prime,X)=0$. Therefore $\hom_{\sC/\Sigma}(\sT_{\leq i+1},X)=0$. It follows by induction that $\hom_{\sC/\Sigma}(\sT,X)=0$. 

%Now assume $\dhomd(F,G)=0$ for all $F\in\kS$ and $G\in\kS^\prime$. From the argument above $\dhomd(H,G)=0$ for all $H\in\sT$ and $G\in\kS^\prime$. So we assume that any $H\in\sT$ has no morphisms to any object in $\sT^\prime_{\leq j}$. Let $W$ be an object in $\sT^\prime_{\leq j+1}$. If $W$ is a shift or direct summand of an object in $\sT^\prime_{\leq j}$ then by our assumption none of it's direct summands have any non-zero morphisms from an $H$. If $W$ is cone of a morphism between objects $U_1$ and $U_2$ in $\sT^\prime_{\leq j}$ then it fits into the triangle \[U_1\rightarrow U_2\rightarrow W\rightarrow U_1[1].\] For any shift of $H$ we have the following part of a long exact sequence \[\cdots\rightarrow\Hom_\sD( H[k_5],U_2)\rightarrow\Hom_\sD(H[k_5],W)\rightarrow\Hom_\sD(H[k_5],U_1[1])\rightarrow\cdots\] for any $k_5\in\IZ$.

%By our assumption, there is no morphism from $H[k_5]$ to $U_2$ or $U_1[1]$ for any $k_5\in\IZ$. Since the sequence above is exact then $\hom_\sD(H[k_5], W)=0$. As any non-zero morphism from $H[k_5]$ to a direct summand of $W$ would induce a non-zero morphism from $H[k_5]$ to $W$ then there are no morphisms from any shift of $H$ to any direct summand of $W$. Hence, there are no objects in $\sT^\prime_{\leq j+1}$ with morphisms from an object in $\sT$ and the result follows by induction.

Now assume $\dhomd(\kS,\kS^\prime)=0$, it follows from the above argument that $\dhomd(\sT,\sT^\prime_{\leq0})=0$. An similar argument to the one above gives $\dhomd(\sT,\sT^\prime)=0$ \qed

\begin{remark}\label{rzhom} If $\sT=\td(\kS)$ and $Z$ is a non-zero object in $\sD$ such that $\dhomd(Z,\kS)=0$,  then $Z$ cannot be in $\sT$, since $\dhomd(Z,Z)\not=0$.\end{remark} 

\section{Thick subcategories generated by string objects}
In this section, we provide the details of the thick subcategories that we later give a classification for. Firstly, the thick subcategories we will consider are finitely generated.

\defn\label{tsfg} A thick subcategory generated by a set of objects $\kS$ in $\sD$ is \textit{finitely generated} if $\left|\kS\right|<\infty$. 
\subsection{Ext-connected thick subcategories}
A general category can be written as a disjoint union of connected categories. We consider a similar concept for finitely generated thick subcategories.

\defn \label{conts} A subcategory $\sU$ of $\sD$ which is closed under shift is \textit{Ext-connected} if for any two indecomposable objects $X,Y\in\sU$, there exists a sequence of indecomposable objects in $\sU$, \[X=X_0,X_1,...,X_m=Y\] for some finite $m\in\IN$ such that  $\dhomd(X_i,X_{i+1})\not=0$ or $\dhomd(X_{i+1},X_i)\not=0$ (for each $i=0,...,m-1$).

\lem \label{tsecgsec} A thick subcategory $\sT=\td(\kS)$ is Ext-connected if and only if $\kS$ is Ext-connected. 

\pf Assume $\kS$ is not Ext-connected. It follows from Definition~\ref{conts} that there exists $\kS_1$ and $\kS_2$ such that $\kS=\kS_1\cup\kS_2$ and $\dhomd(\kS_1,\kS_2)=0=\dhomd(\kS_2,\kS_1)$. By Lemma~\ref{zhom}, $\dhomd(\td(\kS_1),\td(\kS_2))=0=\dhomd(\td(\kS_2),\td(\kS_1))$ As $\sT=\td(\kS_1,\kS_2)$, it follows that $\sT$ is not Ext-connected.

Now we let $\kS$ be Ext-connected and proceed inductively. Assume $\sT_{\leq b}$ is Ext-connected for some $b\geq0$.

We need to show $\sT_{\leq b+1}$ is Ext-connected. Let $U$ be an object generated in $b+1$ steps. We first want to show that every indecomposable direct summand of $U$ has a morphism to or from an object in $\sT_{\leq b}$. If $U$ is direct summand of an object in $\sT_{\leq b}$ then each of its direct summands has an identity morphism to an direct summand of a object in $\sT_{\leq b}$. If $U$ is the cone of a morphism between $A_1,A_2\in\sT_{\leq b}$ then by definition $U$ fits into the triangle $
A_1\rightarrow A_2\xrightarrow{f}U\rightarrow A_1[1].$ 

Let $U^\prime$ be any indecomposable direct summand of $U$. Either $f$ induces a non-zero morphism from a direct summand of $A_2$ to $U^\prime$ or $f$ induces a zero morphism from $A_2$ to $U^\prime$, and so $U^\prime$ is a direct summand of $A_1[1]$. In both cases, every direct summand $U^\prime$ has a morphism from a direct summand of an object in $\sT_{\leq b}$. This means that for any indecomposable object $U^\prime$ in $\sT_{\leq b+1}$, there exists an indecomposable object in $\sT_{\leq b}$ with a morphism to $U^\prime$. Let $V_1$ and $V_2$ be indecomposable direct summands of objects generated in $b+1$ steps. There exists indecomposable objects $W_1$ and $W_2$ in $\sT_{\leq b}$ with morphisms to $V_1$ and $V_2$ respectively. By our assumption that $\sT_{\leq b}$ is Ext-connected and we have the sequence of indecomposable objects in $\sT_{\leq b}$ \[
W_1=Z_0, Z_1,...,Z_c=W_2
\] with morphisms in some direction between shifts of consecutive objects. There exists the sequence of objects in $\sT_{\leq b+1}$\[
V_1,Z_0,Z_1,...,Z_c,V_2
\] with morphisms in some direction between shifts of consecutive objects. Hence, $\sT_{\leq b+1}$ is Ext-connected and we have proved by induction that $\sT$ is Ext-connected.\qed

\cor Every finitely generated thick subcategory is equivalent to a direct sum of Ext-connected thick subcategories. 

\pf Any finite set of objects can be written as $\kS=\kS_1\sqcup\kS_2\sqcup...\sqcup\kS_n$ (for a finite $n\in\IZ$) such that $\dhomd(\kS_i,\kS_j)=0=\dhomd(\kS_j,\kS_i)$ for all $i\not=j$ and each $\kS_k$ is Ext-connected. From Lemma~\ref{zhom} and Remark~\ref{rzhom}, $\td(\kS_i)\cap\td(\kS_j)=0$. Therefore, $\td(\kS)=\bigoplus_{l=1}^n\td(\kS_l)$ and by Lemma~\ref{tsecgsec} each $\td(\kS_k)$ is Ext-connected.\qed

\subsection{Thick subcategories containing a string object}
We would like to know which thick subcategories contain a string object, the next result tells us such thick subcategories must contain a string object in any generating set.

\lem \label{bdgnstr}Let $\sT=\td(B_i\mid i\in I)$ such that each $B_j$ is a band object. Then there does not exist any string object in $\sT$.

\pf Let $\sT$ be a thick subcategory generated by a set of band objects $\kS$ in $\sD$. By Lemma~\ref{bddim}, we may assume all of the band objects in $\kS$ are 1-dimensional. 

Let $X\in\sD$ be a string object. The arc $\gamma_X$ has an end point $p$ on a boundary component $\kb$. Starting at $p$ and travelling anti-clockwise along $\kb$, the boundary component passes through the end points of a sequence of laminates $\kl=(\ell_1,\ell_2,...,\ell_n)$ before passing through $p$ again. By \cite[Lemma 2.8]{Opper} there exists a unique arc $\gamma_Y$ upto homotopy equivalence that crosses the sequence of laminates $\kl$. Note that, $\gamma_Y$ is a closed arc with end point $p$. As $\kb$ could also be thought of as a closed path crossing the set of laminates $\kl$ which has an end point $p$ then we see that $\gamma_Y\simeq\kb$. A closed curve only has points within the interior of $\Sigma$ and no points on any boundary component. This implies, as $\gamma_Y$ is homotopic to a boundary component, that it does not intersect any closed curve on $\Sigma$ when in minimal position. It follows from \cite[Theorem 3.3]{Opper} and Lemma~\ref{zhom} that \[\dhomd(Y,\sT)=\dhomd(\sT,Y)=0.\] However, \cite[Theorem 3.3]{Opper} also implies that at least one of $\dhomd(Y,X)$ or $\dhomd(X,Y)$ is non-zero. Hence $X$ cannot be in $\sT$. The result follows.\qed

A direct consequence of the above is the following.
\cor A thick subcategory $\sT$ contains a string object if and only if every generating set of $\sT$ contains a string object. 

Now we consider Ext-connected thick subcategories and aim to show exactly which ones are generated by string objects. The next lemma is the key step towards this aim.
\lem \label{bdtostr} Let $\sT=\td(\kS)$ be a Ext-connected thick subcategory with a string object in $\kS$ and $B$ be a band object in $\kS$. Then $\sT=\td(\kS\backslash\{B\}, A_1, A_2)$ where $A_1$ and $A_2$ are string objects.

\pf Let $\sT=\td(\kS)$ be Ext-connected with a string object $S\in\kS$ and $B$ be a band object in $\kS$. We would like to find an indecomposable string object $Y_B$ with a morphism to $B$. As $\sT$ is Ext-connected then there exists a finite sequence of indecomposable objects in $\sT$ \[S=X_0,X_1,...,X_n=B\]
such that there is a non-zero morphism in at least one direction between $X_i$ and a shift of $X_{i+1}$ for $0\leq i<n$. Let $X_j$ be the string object such that for every $j^*>j$, then $X_{j^*}$ is a band object. %So we have the sequence \[X_j,X_{j+1},...,X_n=B\] with morphisms in some direction between $X_k$ and a shift of $X_{k+1}$ for $j\leq k<n$. 
(Note that, if there exist a morphism between $X_k$ and $X_{k+1}$ then there exists a basis morphism between $X_k$ and $X_{k+1}$.) We may assume $X_{j+1},X_{j+2},...,X_n$ are one-dimensional band objects as any morphism between a band object $B^*$ and any object $Z\in \sD$ is induced by a morphism between $Z$ and the one-dimensional band object corresponding to the same homotopy band and scalar as $B^*$ \cite[Lemma 5.4 and Theorem 5.10]{Arnesen}. This implies that the basis morphism between $X_{k^\prime}$ and $X_{k^\prime+1}$ for $j\leq k^\prime\leq n-1$ corresponds to an internal intersection point between either an arc and closed curve or two closed curves. By \cite[Theorem 3.7]{Opper} an internal intersection point between two paths corresponds to an morphism in each direction between shifts of the paths' corresponding objects. By fixing $B=X_n$ and if necessary shifting the any other $X_{k^\prime}$, we may assume without loss of generality, that there exists morphism from $X_{k^\prime}$ to $X_{k^\prime+1}$ which is in degree zero. i.e. we have the sequence of non-zero morphisms\[
X_j\xrightarrow[]{h_j}X_{j+1}\xrightarrow{h_{j+1}}X_{j+2}\rightarrow\cdots\rightarrow X_{n-1}\xrightarrow{h_{n-1}}X_n=B
\] If $j=n-1$ then $Y_B=X_{n-1}$. Otherwise, there exists the exact triangle \[
X_j\xrightarrow{h_j} X_{j+1} \rightarrow C(h_j)\rightarrow X_j[1]
\] and by \cite[Propositions 2.9, 3.4, 4.2 \& 5.2(2)]{Canakci} $C(h_j)$ is an indecomposable string object. By the definition of thick subcategories $X_j$ and $X_{j+1}$ are in $\sT$, as is $C(h_j)$. Applying $\Hom(-,X_{j+2})$ to the triangle, we have the long exact sequence \[
\dots\rightarrow\Hom(C(h_i),X_{j+2})\xrightarrow{m}\Hom(X_{j+1},X_{j+2})\rightarrow\Hom(X_j,X_{j+2})\rightarrow\dots
\] There are now two cases. If $\hom(X_j,X_{j+2})=0$ then $m$ is surjective. As $h_{j+1}$ is a non-zero morphism in $\Hom(X_{j+1},X_{j+2})$  then there exists a non-zero morphism in $\Hom(C(h_j),X_{j+2})$ whose image under $m$ is $h_{j+1}$. There exists a sequence of non-zero morphisms of length $n-j-1$ \[
C(h_j)\xrightarrow[]{}X_{j+2}\rightarrow X_{j+3}\rightarrow\cdots\rightarrow X_{n-1}\xrightarrow{h_{n-1}}X_n=B
\] from a string object in $\sT$ to $B$. However, if there exists a non-zero morphism in $\Hom(X_j,X_{j+2})$ then we again have a sequence of non-zero morphisms of length $n-j-1$ \[
X_j\xrightarrow[]{}X_{j+2}\rightarrow X_{j+3}\rightarrow\cdots\rightarrow X_{n-1}\xrightarrow{h_{n-1}}X_n=B\] from a string object in $\sT$ to $B$.

Working inductively, we find an indecomposable string object $Y_B$ in $\sT$ with a non-zero basis morphism to $t:Y_B\rightarrow B$. This morphism fits into the triangle: \[
Y_B\xrightarrow{t}B\rightarrow C(t)\rightarrow Y_B[1]
\] and by \cite[Propositions 2.9, 3.4, 4.2 \& 5.2(2)]{Canakci} $C(t)$ is an indecomposable string object. Since $Y_B$ and $B$ are in $\sT$ then by the definition of thick subcategories $C(t)\in \sT$. Therefore, the result follows as $\sT=\td(\kS\backslash\{B\},Y_B,C(t))$. \qed

We conclude this section with a result that says in order to achieve the aim of this paper we only need to classify Ext-connected thick subcategories generated by a finite set of string objects in $\sD$.
\thm\label{cstobj} Let $\sT$ be a Ext-connected thick subcategory containing a string object. Then $\sT$ is generated by string objects. Furthermore, if $\sT$ is finitely generated then it is generated by finitely many  string objects.

\pf From Lemma~\ref{bdgnstr} there must exist a string object in the generating set of $\sT$. By Lemma~\ref{bdtostr} each band object in the generating set of $\sT$ can be replaced by two string objects. % By Corollary~\ref{corstrgen} each string object can be replaced by a finite number of exceptional and spherelike string objects. 
The result follows.\qed

This section has told us that an Ext-connected thick subcategory either contains a string object and is generated by string objects or it is generated by band objects.

\defn We denote the poset of all Ext-connected thick subcategories generated by finitely many string objects as $\tds$, where the elements are partially ordered by inclusion.

\section{Concatenations of paths}
The geometric arguments in this section aid us to prove results about generating sets of thick subcategories in $\tds$ later in the paper.

Throughout this section we will considered piecewise smooth paths with a finite number of self-intersection points that do not coincide with any singular points. The following lemma allows us to assume (without loss of generality) that such a path is homotopic to a smooth one.
 
\lem \label{cosm} If a piecewise smooth path $\gamma$ is smooth in a neighbourhood of all self-intersection points then it is homotopic to a smooth path which has the same number of intersection points with any path on $\Sigma$ as $\gamma$, when they are both in minimal position.

\pf If a point $\gamma(t)$ on the path is not smooth then there exists a neighbourhood $\ku=\kn_{\gamma(t)}$ around $\gamma(t)$ such that $\psi=\restr{\gamma}{[t_1,t_2]}$ is the only subpath of $\gamma$ in $\ku$, $t\in(t_1,t_2)$ and $\gamma(t)$ is the only singular point in the closure $\overline{\ku}$ of $\ku$. There are two smooth subpaths $\psi_1=\restr{\gamma}{[t_1,t]}$ and $\psi_2=\restr{\gamma}{[t,t_2]}$ of $\gamma$ such that $\psi=\psi_1\psi_2$. We have the figure below.

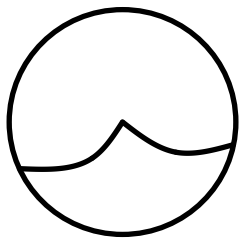
\begin{center}
\begin{tikzpicture}[line cap=round,line join=round,>=triangle 45,x=1cm,y=1cm]
\draw [line width=2pt] (6.62,-0.34) circle (1.490503270710937cm);
\draw[line width=2pt, smooth,samples=100,domain=0:0.6072751397528033, variable=\t] plot({-0.4578512155679058*\t^3)+1.7434558855355122*\t+5.26377964781218},{0.8188455661034301*\t^3-0.07427578342985591*\t-0.9582769252620937});
\draw[line width=2pt, smooth,samples=100,domain=0.6072751397528033:1,variable=\t] plot({0.7079808004641731*\t^3-2.1239424013925197*\t^2+3.0332733041680595*\t+5.0026882967602875},{-1.266190674250735*\t^3+3.7985720227522055*\t^2-2.3810541394077904*\t-0.49132720909368});
\draw[line width=2pt, smooth,samples=100,domain=0:0.5126130633294061,variable=\t] plot({0.18901547712157066*\t^3+1.315884398328176*\t+6.619999999999999},{0.9612430109031274*\t^3-1.0329035536071094*\t-0.33999999999999997});
\draw[line width=2pt, smooth,samples=100,domain=0.5126130633294061:1,variable=\t] plot({-0.19879852218822003*\t^3+0.5963955665646601*\t^2+1.010164239995389*\t+6.672238715628173},{-1.0109949351311027*\t^3+3.0329848053933084*\t^2-2.587651185731316*\t-0.07433868453088999});
\end{tikzpicture}
\end{center}
 
%If $\phi_1=\restr{\gamma}{[t^\prime_1, t]}$ and $\phi_2=\restr{\gamma}{[t, t_2^\prime]}$ are subpaths of $\gamma$ such that $t_1<t^\prime_1$ and $t_2^\prime<t_2$ then there exists a sub-region $\ku_2$ such that $\phi_1\phi_2$ is the only subpath of $\gamma$ contained in $\ku_2$. 
(We identify $\ku$ with an open subset in $\IR^2$.) By \cite[Lemma 1]{de Boor}, there exists a cubic spline $\mu$ in $\ku$ such that $\mu\simeq\psi$ and $\restr{\gamma}{[0,t_1]}\mu\restr{\gamma}{[t_2,1]}$ has no singular points in $\overline{\ku}$, since we can stipulate the value of the cubic spline $\mu$ and its first derivative at both end points. The path $\mu$ has no self-intersection points and it does not intersect $\gamma$ outside of $\ku$ apart from at its common end points with $\restr{\gamma}{[0,t_1]}$ and $\restr{\gamma}{[t_2,1]}$, neither of which causes an self-intersection point of $\restr{\gamma}{[0,t_1]}\mu\restr{\gamma}{[t_2,1]}$. Therefore, $\gamma$ is homotopic to a path $\overline{\gamma}$ which has no singular point in $\overline{\ku}$ and the same number of intersections as $\gamma$.

If we repeat this process for every singular point of $\gamma$ then $\gamma$ is homotopic to a smooth path with the same number of self-intersection points. By \cite[Corollary 4]{Thurston}, we can move the smooth path into minimal position with respect to any arc that intersects $\gamma$ and the result follows.\qed

Now we will discuss simple concatenations, we begin by defining them in a general setting.

\defn \label{dsimcc} Let $\psi$ and $\phi$ be two paths such that $\psi(1)=\phi(0)=p$ and the subpaths $\restr{\psi}{(t_1,1]}$ and $\restr{\phi}{[0,t_2)}$ be in $\kn_p$ with $\psi(t_1)$ and $\psi(t_2)$ on the boundary of $\overline{\kn_p}$. The concatenation $\psi\phi$ is \textit{simple} if there exists a path in $\overline{\kn_p}$ which is equivalent to $\restr{\psi}{[t_1,1]}\restr{\phi}{[0,t_2]}$ up to end point fixing homotopy and which only intersects $\psi$ and $\phi$ at its end points.

Before we prove our first result which relates to simple concatenations, we need to discuss the following.

As we have done before, we identify the tubular neighbourhood $\kt$ of a path with $[t_1,t_2]\times[-\varepsilon,\varepsilon]$ throughout this section. Define two subspaces of $\kt$ by $\kt^+=\{(y_1,y_2)\in\kt\mid y_2\geq0\}$ and $\kt^-=\{(y_1,y_2)\in\kt\mid y_2\leq0\}$. %The result below says if we can take a tubular neighbourhood of a path then we can construct two other paths such that all three paths are homotopic and they only intersect each other at their end points.
We now prove a result relating to these subsets.
\lem\label{tbnhnsip} Let $\psi$ be a path with no self-intersection points and $\phi$ be a path such that $\kt\cap\phi$ is either a subset of $\kt^{+}$ or $\kt^{-}$. Then there exists two paths $\psi_{+}$ and $\psi_{-}$ in a tubular neighbourhood $\kt$ of $\psi$ such that $\psi_{+}\cap\phi\subseteq\{\psi_{+}(0),\psi_{+}(1)\}$ or $\psi_{-}\cap\phi\subseteq\{\psi_{-}(0),\psi_{-}(1)\}$.

\pf A tubular neighbourhood $\kt$ of a non-self-intersecting path $\psi:[t_1,t_2]\rightarrow\Sigma$ consists of normal vectors to $\psi$ of a length less than $\varepsilon_1$ for some $\varepsilon_1>0$. There are two sections $\psi_{+}$ and $\psi_{-}$ given by $x\mapsto(x, \frac{2\varepsilon_1}{3t_1^2+2t_1t_2+t_2^2}(x-t_1)(x-t_2))$ and $x\mapsto(x, \frac{2\varepsilon_1}{3t_1^2+t_1t_2+t_2^2}(x-t_1)(t_2-x))$ respectively which only intersect the zero section when $x\in\{t_1,t_2\}$. Note that $\psi_{+}$ (respectively $\psi_{-}$) is a subset of $\kt^{+}$ (respectively $\kt^{-}$). We see this below.
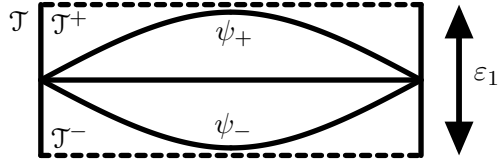
\begin{figure}[ht]
    \centering
    \begin{tikzpicture}[line cap=round,line join=round,>=triangle 45,x=1cm,y=1cm]
\draw [line width=2pt] (0,2)-- (0,0);
\draw [line width=2pt] (5,2)-- (5,0);
\draw [line width=2pt, dashed] (5,2)-- (0,2);
\draw [line width=2pt,dashed] (0,0)-- (5,0);
\draw[line width=2pt, smooth,samples=100,domain=0:0.5, variable=\t] plot({3.5527136788005006E-16*\t^3+5*\t,-3.599999999999998*\t^3+2.6999999999999975*\t+1.0000000000000009});
\draw[line width=2pt, smooth,samples=1000,domain=0.5:1, variable=\t] plot({-3.5527136788005006E-16*\t^3+1.0658141036401502E-15*\t^2+4.999999999999999*\t,3.599999999999998*\t^3-10.799999999999994*\t^2+8.099999999999994*\t+0.10000000000000127});
\draw [line width=2pt] (0,1)-- (5,1);
\draw[line width=2pt, smooth,samples=1000,domain=0:0.5, variable=\t] plot({3.5527136788005006E-16*\t^3+5*\t,3.5999999999999983*\t^3-2.6999999999999984*\t+0.9999999999999994});
\draw[line width=2pt, smooth,samples=100,domain=0.5:1, variable=\t] plot({-3.5527136788005006E-16*\t^3+1.0658141036401502E-15*\t^2+4.999999999999999*\t,-3.5999999999999983*\t^3+10.799999999999995*\t^2-8.099999999999996*\t+1.8999999999999988});
\draw (2.15,0.7) node[anchor=north west] {$\psi_{-}$};
\draw (2.16,1.94) node[anchor=north west] {$\psi_{+}$};
\draw (0,2.05) node[anchor=north west] {$\mathcal{T}^{+}$};
\draw (-0.55,2.05) node[anchor=north west] {$\mathcal{T}$};
\draw (0,0.55) node[anchor=north west] {$\mathcal{T}^{-}$};
\draw [->,line width=2pt] (5.5,0.1) -- (5.5,2);
\draw [->,line width=2pt] (5.5,1.9) -- (5.5,0);
\draw (5.56,1.3) node[anchor=north west] {$\varepsilon_1$};
\end{tikzpicture}
\end{figure}

Let $\phi$ be a path such that the intersection of $\phi$ and $\kt$ is either a subset of $\kt^+$ or $\kt^-$. If all the points in $\phi\cap\kt$ are in $\kt^+$ then $\phi$ can only possibly intersect $\psi_{-}$ at the end points of $\psi_{-}$. Otherwise, $\phi$ can only possibly intersect $\psi_{+}$ at the end points of $\psi_{+}$. \qed

\lem \label{simcc} Let $\psi\phi$ be a simple concatenation. Suppose $\psi$ has no internal self-intersection points and $\psi$ has no internal intersection points with $\phi$. Then path $\psi\phi$ only intersects $\psi$ at common end points when they are in minimal position. Furthermore, the number of internal self-intersection points of $\phi$ is an upper bound for the number of internal self-intersection points of $\psi\phi$ in minimal position.

\pf Let $\psi\phi$ be a simple concatenation at the point $p=\psi(1)=\phi(0)$. We have the subpaths $\restr{\psi}{[t_1,1]}$ and $\restr{\phi}{[0,t_2]}$ of $\psi$ and $\phi$ respectively which are homotopic to radii of $\kn_p$. Denote the path between $\psi(t_1)$ and $\phi(t_2)$ in $\kn_p$, by $\chi$. 

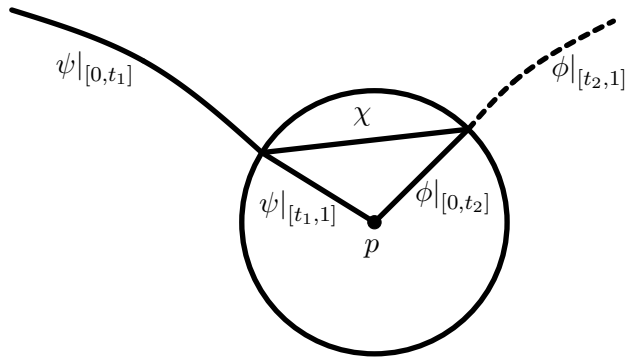
\begin{figure}[ht]
    \centering
    \begin{tikzpicture}[line cap=round,line join=round,>=triangle 45,x=1cm,y=1cm]

\draw [line width=2pt] (21.16636570058514,3.7411279409123743) circle (1.7446878802355925cm);
\draw [line width=2pt] (21.16636570058514,3.7411279409123743)-- (19.686289571089553,4.664875985748373);
\draw [line width=2pt] (21.16636570058514,3.7411279409123743)-- (22.40004633175385,4.974808572080806);
\draw [line width=2pt] (22.40004633175385,4.974808572080806)-- (19.686289571089553,4.664875985748373);
\draw[line width=2pt, smooth,samples=100,domain=0:0.48431024875974704, variable=\t] plot({-0.5665600328886748*\t^3+3.7050869188508675*\t+16.397094021809764,-0.9957920885021287*\t^3-1.160458527781957*\t+6.556309717311731});
\draw[line width=2pt, smooth,samples=100,domain=0.48431024875974704:1, variable=\t] plot({0.5320850953615509*\t^3-1.5962552860846528*\t^2+4.478169713538588*\t+16.272290048274066,0.9351985625767647*\t^3-2.8055956877302943*\t^2+0.19832021766197552*\t+6.3369528932399275});
\draw[line width=2pt,dashed, smooth,samples=100,domain=0:0.36240483745319696, variable=\t] plot({0.5433632675051677*\t^3+1.5863971802034658*\t+22.40004633175385,-0.6714927094233843*\t^3+1.832326223098362*\t+4.974808572080806});
\draw[line width=2pt,dashed, smooth,samples=100,domain=0.36240483745319696:1, variable=\t] plot({-0.3088440568646777*\t^3+0.9265321705940331*\t^2+1.2506174395241774*\t+22.440609065867502,0.381671976991709*\t^3-1.145015930975127*\t^2+2.247285535444724*\t+4.924680818033948});
\draw (16.84015597525919,6.15) node[anchor=north west] {$\restr{\psi}{[0,t_1]}$};
\draw (19.5170595506974,4.3) node[anchor=north west] {$\restr{\psi}{[t_1,1]}$};
\draw (21.523953768867464,4.45) node[anchor=north west] {$\left. \phi\right|_{[0,t_2]}$};
\draw (23.317299770924688,6.1) node[anchor=north west] {$\left. \phi\right|_{[t_2,1]}$};
\draw (20.751210630934206,5.3989433893884) node[anchor=north west] {$\chi$};
\draw (20.9,3.6502691142) node[anchor=north west] {$p$};
\begin{scriptsize}
\draw [fill=black] (21.16636570058514,3.7411279409123743) circle (2.5pt);
\end{scriptsize}
\end{tikzpicture}
\end{figure}

%There exists a tubular neighbourhood $\kt$ of $\restr{\psi}{[0,t_3]}$ for some $t_3\in(t_1,1)$, so we have the sections $x\mapsto(x,x(x-t_3))$ and $x\mapsto(x,x(t_3-x))$ which only intersect $\psi$ at their end points. One of these sections $\sigma$ intersects $\chi$ once (since $\sigma\cap\kn_p$ is a subpath of a path homotopic to a radii of $\kn_p$) at the point $\chi(t_4)=\sigma(t_5)$.
%\begin{center}
%\includegraphics[]{images/simcc2.PNG}    \end{center}
%The path $\restr{\sigma}{[0,t_5]}\restr{\chi}{[t_4,1]}\restr{\phi}{[t_2,1]}$ only intersects $\psi$ at common end points, since $\psi$ only intersects $\sigma$ at $\sigma(0)$, $\chi$ at $\chi(0)$ and possibly $\restr{\phi}{[t_2,1]}$ at $\phi(1)$. As $\psi\simeq\sigma\restr{\phi}{[t_3,1]}$ and $\restr{\chi}{[t_4,1]}\simeq\restr{\sigma}{[t_5,1]}\restr{\psi}{[t_3,1]}\restr{\phi}{[0,t_2]}$, then $\psi\phi\simeq\restr{\sigma}{[0,t_5]}\restr{\chi}{[t_4,1]}\restr{\phi}{[t_2,1]}$. Hence, $\psi\phi$ only intersects $\phi$ at common end points when in minimal position.

By definition~\ref{dsimcc}, $\chi$ only intersects $\psi$ at $\psi(t_1)$. There exists a tubular neighbourhood $\kt$ of $\restr{\psi}{[0,t_3]}$ for some $t_1<t_3<1$ such that $\kt\cap\restr{\psi}{[0,t_1]}\chi$ is either a subset of $\kt^{+}$ or $\kt^{-}$. By Lemma~\ref{tbnhnsip} there exist a path $\sigma\simeq\restr{\psi}{[0,t_3]}$ that only intersects $\restr{\psi}{[0,t_1]}\chi$ at the point $\psi(0)$. We see this in the figure below.

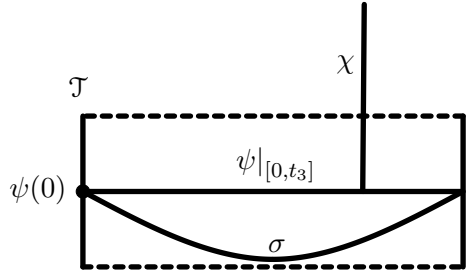
\begin{figure}[ht]
    \centering
\begin{tikzpicture}[line cap=round,line join=round,>=triangle 45,x=1cm,y=1cm]
\draw [line width=2pt] (0,2)-- (0,0);
\draw [line width=2pt] (5,2)-- (5,0);
\draw [line width=2pt,dashed] (5,2)-- (0,2);
\draw [line width=2pt,dashed] (0,0)-- (5,0);
\draw [fill=black] (0,1) circle (2.5pt);
%\draw[line width=2pt, smooth,samples=100,domain=0:0.5, variable=\t] plot({3.5527136788005006E-16*\t^3+5*\t,-3.599999999999998*\t^3+2.6999999999999975*\t+1.0000000000000009});
%\draw[line width=2pt, smooth,samples=100,domain=0.5:1, variable=\t] plot({-3.5527136788005006E-16*\t^3+1.0658141036401502E-15*\t^2+4.999999999999999*\t,3.599999999999998*\t^3-10.799999999999994*\t^2+8.099999999999994*\t+0.10000000000000127});
\draw [line width=2pt] (0,1)-- (5,1);
\draw[line width=2pt, smooth,samples=100,domain=0:0.5, variable=\t] plot({3.5527136788005006E-16*\t^3+5*\t,3.5999999999999983*\t^3-2.6999999999999984*\t+0.9999999999999994});
\draw[line width=2pt, smooth,samples=100,domain=0.5:1, variable=\t] plot({-3.5527136788005006E-16*\t^3+1.0658141036401502E-15*\t^2+4.999999999999999*\t,-3.5999999999999983*\t^3+10.799999999999995*\t^2-8.099999999999996*\t+1.8999999999999988});
\draw [line width=2pt] (3.68,1)-- (3.7,3.48);
\draw (1.84,1.725) node[anchor=north west] {$\left. \psi\right|_{[0,t_3]}$};
\draw (2.3,0.5) node[anchor=north west] {$\sigma$};
\draw (3.2,2.94) node[anchor=north west] {$\chi$};
\draw (-0.32,2.66) node[anchor=north west] {$\mathcal{T}$};
\draw (-1.1,1.35) node[anchor=north west] {$\psi(0)$};
\end{tikzpicture}\end{figure}

Note that, by our assumptions, (when in minimal position) the paths $\restr{\psi}{[0,t_1]}\chi$, $\restr{\psi}{[t_3,1]}$ and $\restr{\phi}{[t_2,1]}$ only possibly intersect $\psi$ at an end point of $\psi$. As $\kn_p$ is a (half-)disc then $\chi\simeq\restr{\psi}{[t_1,1]}\restr{\phi}{[0,t_2]}$ and $\psi\phi\simeq\restr{\psi}{[0,t_1]}\chi\restr{\phi}{[t_2,1]}$. By \cite[Corollary 4]{Thurston} there exist a path homotopic to $\psi$ which only intersects $\restr{\psi}{[0,t_1]}\chi\restr{\phi}{[t_2,1]}$ at their common end point. This proves the first part this result.

The paths $\restr{\psi}{[0,t_1]}$, $\restr{\phi}{[t_2,1]}$ and $\chi$ do not intersect except at common end points and $\restr{\phi}{[t_2,1]}$ is the only one of them to potentially have internal self-intersection points. This means the path $\restr{\psi}{[0,t_1]}\chi\restr{\phi}{[t_2,1]}$ only has internal intersection points that correspond to internal self-intersection points of $\restr{\phi}{[t_2,1]}$. This implies, when in minimal position, $\psi\phi$ has at most the same number of self-intersection points as $\restr{\phi}{[t_2,1]}$ since $\chi\simeq\restr{\psi}{[t_1,1]}\restr{\phi}{[0,t_2]}$ and $\psi\phi\simeq\restr{\psi}{[0,t_1]}\chi\restr{\phi}{[t_2,1]}$. These self-intersection points all correspond to self-intersection points of $\phi$ and the result follows. \qed

\subsection{Concatenations of arcs in an arc-collection}
In Section 6, we show a correspondence between the thick subcategories in $\tds$ and arc-collections. Let us define these now.

\defn A collection of arcs in $\Sigma$ with no internal intersection points is called an arc-collection. We say an arc-collection is connected if the union of the images of the arcs is path connected.

We denote the set of arcs corresponding to a set of string objects $\kS$ by $\ka_{\kS}$. Now we see some useful tools we apply throughout the remainder of the paper. The first helps us deal with concatenations of arcs in an arc-collection.

\defn \label{acscc} Let $\phi$ and $\psi$ be two distinct arcs in an arc-collection $\ka$ and $\phi\psi$ be a concatenation at a marked point $p=\phi(1)=\psi(0)$. Consider subpaths $\restr{\phi}{(t_1,1]}$ and $\restr{\psi}{[0,t_2)}$ in $\kn_p$ with $\phi(t_1)$ and $\psi(t_2)$ on the boundary of $\overline{\kn_p}$. The concatenation $\phi\psi$ is \textit{simple with respect to $\ka$} if it is a simple concatenation and when moved into minimal position the path $\restr{\phi}{[t_1,1]}\restr{\psi}{[0,t_2]}$ does not intersect any arc in $\ka\setminus\{\phi,\psi\}$.

We produce a version of Lemma~\ref{simcc} for simple concatenations with respect to an arc-collection.

\lem \label{lemacscc} If $\phi$ and $\psi$ are distinct arcs in an arc-collection $\ka$ and $\phi\psi$ is a simple concatenation with respect to $\ka$ then $\phi\psi$ does not internally intersect any arc in $\ka$.

\pf Let $\phi$ and $\psi$ be distinct arcs in an arc-collection $\ka$ and $\phi\psi$ be a simple concatenation with respect to $\ka$. Let $t_1$ and $t_2$ be as in Definition~\ref{acscc}. There exists $\beta$ which is fixed end point homotopic to $\restr{\phi}{[t_1,1]}\restr{\psi}{[0,t_2]}$ in $\kn_p$ such that $\restr{\phi}{[0,t_1]}\beta\restr{\psi}{[t_2,1]}$ is homotopic to $\phi\psi$ and $\restr{\phi}{[0,t_1]}\beta\restr{\psi}{[t_2,1]}$ only possibly intersects any arc in $\ka\setminus\{\phi, \psi\}$ at the marked points $\phi(0)$ and $\psi(1)$. By Lemma~\ref{simcc}, $\phi\psi$ only intersects either $\phi$ or $\psi$ at marked points when in minimal position. Therefore, by \cite[Corollary 4]{Thurston}, $\phi\psi$ is homotopic to a path which does not internally intersect any arc in $\ka$.\qed 

We now have a direct consequence of the above.
\cor Let $\ka$ be an arc-collection and suppose $\gamma_0,\gamma_1,...,\gamma_k\in\ka$. If we have an arc-collection $\ka^\prime=\ka\cup\{\gamma_0^\prime\gamma_1^\prime...\gamma_{k-1}^\prime\}$ and $\gamma_0^\prime\gamma_1^\prime...\gamma_k^\prime$ is a simple concatenation with respect to $\ka^\prime$ then $\ka\cup\{\gamma^\prime_0\gamma^\prime_1...\gamma_k^\prime\}$ is an arc-collection.

We would like to generalise Definition~\ref{acscc}.

\defn \label{ssc} Let $\ka$ be an arc-collection and suppose $\gamma_{S_0},\gamma_{S_1},...,\gamma_{S_n} \in\ka$. The concatenation $\gamma^\prime_{S_0}\gamma^\prime_{S_1}...\gamma^\prime_{S_n}$ is a \textit{sequence of simple concatenations} with respect to $\ka$ if and only if $\gamma^\prime_{S_0}\gamma^\prime_{S_1}...\gamma^\prime_{S_{i+1}}$ is a simple concatenation with respect to the arc-collection $\ka\cup\{\gamma^\prime_{S_0}\gamma^\prime_{S_1}...\gamma^\prime_{S_i}\}$ for all $i\in\{0,1,...,n-1\}$.

When discussing arc-collections we will use some terminology relating to graphs.
\defn For an arc-collection $\ka_{\overline{\kS}}$, there exists a connected graph $\Gamma(\overline{\kS})$ such that \begin{align*}
    &\Gamma(\overline{\kS})_0=\{\gamma_S(0),\gamma_S(1)\mid S\in\overline{\kS}\} \text{ and}\\
    &\Gamma(\overline{\kS})_1=\{\gamma_S\mid S\in\overline{\kS}\},
\end{align*} which has a ribbon graph structure induced by its embedding in the surface $\Sigma$. 

For a vertex $v$ of $\Gamma(\overline{\kS})$ we denote the degree of $v$ in $\Gamma(\overline{\kS})$ by $\deg_{\overline{\kS}}(v)$ and the distance given by the minimal path length between two vertices $v$ and $v^\prime$ on the graph $\Gamma(\overline{\kS})$ by $d_{\overline{\kS}}(v, v^\prime)$.

Our work will involve a particular type of arc-collection which we define now.

\begin{defn} An arc-collection $\ka_{\kS_1}$ is a \pac~at $v$ if: \begin{enumerate}[i)]
    \item $d_{\kS_1}(v,u)=1$ for all $u\not=v$ in $\Gamma(\kS_1)_0$ and
    \item $\deg_{\kS_1}(u)=1$ for all $u\not=v$ in $\Gamma(\kS_1)_0$.\end{enumerate}
\end{defn}

To conclude this section we talk about concatenations between arcs in a \pac. We now provide a setup and some notation which will use in the remainder of this section and we will refer to later.

\begin{setup}\label{supac}
    Let $\ka_{\kS_1}$ be a \pac~at $v$. When the arcs in $\ka_{\kS_1}$ are in minimal position, then we can consider the ribbon graph $\Gamma(\kS_1)$ and its half-edges incident to $v$. We have Figure~\ref{figsu},
    \begin{figure}[ht]
    \centering
    \begin{tikzpicture}[line cap=round,line join=round,>=triangle 45,x=1cm,y=1cm]
\draw [shift={(-5.93,-2)},line width=2pt,dashed]  plot[domain=0:3.141592653589793,variable=\t]({1*1.91*cos(\t r)+0*1.91*sin(\t r)},{0*1.91*cos(\t r)+1*1.91*sin(\t r)});
\draw [line width=2pt] (-7.84,-2)-- (-4.02,-2);
\draw [line width=2pt] (-5.86,-2)-- (-4.200213359278925,-1.1900998965409997);
\draw [line width=2pt] (-5.86,-2)-- (-4.8913269057935,-0.39710942252083337);
\draw [line width=2pt] (-5.86,-2)-- (-6.912686892866576,-0.36218851188904067);
\draw [line width=2pt] (-5.86,-2)-- (-7.636378686821466,-1.1418789262805344);
\draw (-6.25,-0.57) node[anchor=north west] {\textbf{...}};
\draw (-5.1,0.1) node[anchor=north west] {$2$};
\draw (-4.3,-0.85) node[anchor=north west] {$1$};
\draw (-8.1,-0.8) node[anchor=north west] {$b$};
\draw (-7.6,0.35) node[anchor=north west] {$b-1$};
\draw (-6.15,-2.05) node[anchor=north west] {$v$};
\begin{scriptsize}
\draw [fill=black] (-5.86,-2) circle (2.5pt);
\end{scriptsize}
\end{tikzpicture}
\caption{}
\label{figsu}
\end{figure}
where the half-edges attached to $v$ are labelled by their pre-image under the following map \[\kh:\{1,2,...,b\}\rightarrow\{\text{half-edges attached to } v\}
\] such that $i$ maps to the $i$th half-edge that the boundary of $\overline{\kn_v}$ crosses in the anti-clockwise direction.

Now we define regions of $\kn_v$ separated by the half-edges in $\Gamma(\kS_1)$ attached to $v$. Let \begin{itemize}
    \item $R_{0}$ be the region bounded by the boundary of $\Sigma$ and the half-edge $\kh(1)$;
    \item $R_{j}$ be the region bounded by the half-edges $\kh(j)$ and $\kh(j+1)$ for $1\leq j<b$ and
    \item $R_{b}$ be the region bounded by the boundary of $\Sigma$ and the half-edge $\kh(b)$.
\end{itemize}
Note that, the interior of these regions are disjoint.

We define $R_{b^\prime}$ to be a terminal region if either $\kh(b^\prime+1)$ does not exist, i.e. $b^\prime=b$, or $\kh(b^\prime+1)$ corresponds to an exceptional object. 

For a non-terminal region $R_{k}$, $\kh(k+1)$ corresponds to a spherelike object $S_{k+1}$ and there exists another half-edge which is incident to $v$ and corresponds to $S_{k+1}$. We define a map \[\tau:\left\{k^\prime\in\{0,1,...,b-1\}\mid R_{k^\prime} \text{ is a non-terminal region}\right\}\rightarrow\{1,2,...,b\}\] such that we have $\tau(k^\prime)=k^{\dprime}$ where $\kh(k^{\dprime})$ is the half-edge corresponding to $S_{k^\prime+1}$ that is not equal to $\kh(k^\prime+1)$.
\end{setup}

We now prove a technical lemma.

\lem \label{exccc}Let $Y$ be a string object such that $\ka_{\kS_1}\cup\{\gamma_Y\}$ is a \pac. Suppose there is a exceptional object $E\in\kS_1$ such that $\gamma_{Z}\simeq\gamma_Y^\prime\gamma_E^\prime$ is a simple concatenation with respect to $\ka_{\kS_1}\cup\{\gamma_Y\}$ at $v$. Then $\gamma_Z$ has a common end point $u\not=v$ with $\gamma_E$ and this is not an oriented intersection point from $\gamma_Z$ to $\gamma_E$.

\pf Let $f$ be a morphism from $Y$ to a shift of an exceptional object $E\in\kS_1$, which does not factor through any object in $\kS_1$. By \cite[Theorem 4.1]{Opper}, the cone $Z$ of $f$ corresponds to $\gamma_Z\simeq\gamma_Y^\prime\gamma_E^\prime$. As $\Sigma$ is an oriented surface, we have the diagram below:

\begin{figure}[ht]
    \centering    
\begin{tikzpicture}[line cap=round,line join=round,>=triangle 45,x=1cm,y=1cm]
\draw [line width=2pt] (3.7,-4.02)-- (2.58,-1.16);
\draw [line width=2pt] (3.7,-4.02)-- (4.62,-2.54);
\draw[dashed,line width=2pt, smooth,samples=100,domain=0:0.27808556644420207, variable=\t] plot({1.5603877301156541*\t^3-2.637878083842081*\t+4.62,3.6783193909157315*\t^3+0.14707159829326458*\t-2.54});
\draw[dashed,line width=2pt, smooth,samples=100,domain=0.27808556644420207:0.5679408281215428, variable=\t] plot({-1.219355819465929*\t^3+2.319019678665034*\t^2-3.282763984778898*\t+4.6797778203512985,-3.6051602913909724*\t^3+6.076291719417289*\t^2-1.5426574263811066*\t-2.3833702490120725});
\draw[dashed, line width=2pt, smooth,samples=100,domain=0.5679408281215428:1, variable=\t] plot({-0.1862814422450226*\t^3+0.5588443267350678*\t^2-2.283088537764665*\t+4.490525653274619,0.05112065592397254*\t^3-0.1533619677719176*\t^2+1.9954172476315535*\t-3.0531759357836084});
\draw (4.1,-3.1) node[anchor=north west] {$\gamma_Y$};
\draw (2.5,-2.6) node[anchor=north west] {$\gamma_E$};
\draw (3.4,-1.75) node[anchor=north west] {$\gamma_Z$};
\draw (3.44,-4.025) node[anchor=north west] {$v$};
\draw (2.25,-0.67) node[anchor=north west] {$u$};
\begin{scriptsize}
\draw [fill=black] (3.7,-4.02) circle (2.5pt);
\draw [fill=black] (2.58,-1.16) circle (2.5pt);
\draw [fill=black] (2.58,-1.16) circle (2.5pt);
\end{scriptsize}
\end{tikzpicture}
    \caption{}
   \label{nalg2}
\end{figure}

The path $\gamma_Z$ has one end point at $u$ and it is not oriented from $\gamma_Z$ to $\gamma_E$.\qed

Finally, we conclude the section by providing some information about sequences of simple concatenations with respect to a \pac, which we use later in this work. 
\lem\label{lsscreg} Let $\gamma_{S_1}^\prime$ and $\gamma_{S_2}^\prime$ be two closed arcs in $\ka_{\kS_1}$ such that there exists $j$ with $\kh(j)$ corresponding to $\gamma_{S_1}$ and $\kh(j-1)$ corresponding to $\gamma_{S_2}$. Then $\phi\simeq\gamma_{S_1}^\prime\gamma_{S_2}^\prime$ is simple concatenation with respect $\ka_{\kS_1}$ and $\phi\cap\kn_v\subset R_{k_1}\cup R_{k_2}$, where $\kh(k_1)$ corresponds to $\gamma_{S_1}$ and $\kh(k_2+1)$ corresponds to $\gamma_{S_2}$. 

\pf First, we consider the following. For a closed arc $\psi\in\ka_{\kS_1}$, there exists a tubular neighbourhood $\kt_\psi$ of $\restr{\psi}{[t_1,t_2]}$ such that $0<t_1<t_2<1$ and $\psi(t_1),\psi(t_2)\in\kn_v$. As $\ka_{\kS_1}$ is an arc-collection we can take such tubular neighbourhoods for all closed arcs in $\ka_{\kS_1}$ so the intersection between any two such tubular neighbourhoods is empty. We have the following diagram:

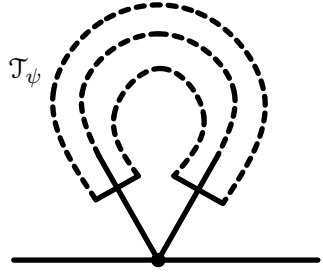
\begin{figure}[ht]
    \centering
    \definecolor{ududff}{rgb}{0,0,0}
\begin{tikzpicture}[line cap=round,line join=round,>=triangle 45,x=1cm,y=1cm]
\draw [line width=2pt] (0,0)-- (4,0);
\draw [line width=2pt] (1.9,0)-- (1.12,1.38);
\draw [line width=2pt] (1.9,0)-- (2.68,1.38);
\draw[line width=2pt,dashed, smooth,samples=100,domain=0:0.20320800823074292, variable=\t] plot({21.908573810187757*\t^3-1.9873160726202899*\t+1.12,-1.1068328479682712*\t^3+4.179400352607718*\t+1.3799999999999994});
\draw[line width=2pt,dashed, smooth,samples=100,domain=0.20320800823074292:0.5, variable=\t] plot({-15.000396811938606*\t^3+22.500595217907904*\t^2-6.559617210857534*\t+1.4297094024441168,-13.781170979924248*\t^3+7.726581022313185*\t^2+2.6092972126299983*\t+1.4863525105972348});
\draw[line width=2pt,dashed, smooth,samples=100,domain=0.5:0.7967919917692571, variable=\t] plot({-15.000396811938593*\t^3+22.500595217907883*\t^2-6.559617210857524*\t+1.4297094024441106,13.781170979924282*\t^3-33.61693191745961*\t^2+23.281053682516394*\t-1.9589402343838296});
\draw[line width=2pt,dashed, smooth,samples=100,domain=0.7967919917692571:1, variable=\t] plot({21.908573810187733*\t^3-65.7257214305632*\t^2+63.738405357942916*\t-17.241257737567445,1.1068328479682337*\t^3-3.320498543904701*\t^2-0.8589018087030128*\t+4.452567504639481});
\draw[line width=2pt, dashed, smooth,samples=100,domain=0:0.22214772677389394, variable=\t] plot({34.22604499794128*\t^3-4.209887120099562*\t+1.0799999999999998,1.2072439191564*\t^3+6.332565815341812*\t+0.7999999999999998});
\draw[line width=2pt, dashed, smooth,samples=100,domain=0.22214772677389394:0.3974178163091544, variable=\t] plot({-28.036288254829046*\t^3+41.49430738722466*\t^2-13.427753180228715*\t+1.7625759969879735,-34.53630529981966*\t^3+23.821044617478968*\t^2+1.0407749041893573*\t+1.1918531071584237});
\draw[line width=2pt, dashed, smooth,samples=100,domain=0.3974178163091544:0.5724460538049162, variable=\t] plot({-26.93828146308761*\t^3+40.18520500282515*\t^2-12.907492569295556*\t+1.6936557183517336,-3.78068831446974*\t^3-12.847445807196545*\t^2+15.613486296117037*\t-0.7386319392027232});
\draw[line width=2pt, dashed, smooth,samples=100,domain=0.5724460538049162:0.7797385260258203, variable=\t] plot({-25.895649871159353*\t^3+38.39465398171014*\t^2-11.88249870312191*\t+1.498071153729952,33.00573255256486*\t^3-76.02217018401888*\t^2+51.77760796584223*\t-7.6393015189210685});
\draw[line width=2pt, dashed, smooth,samples=100,domain=0.7797385260258203:1, variable=\t] plot({33.56745841073348*\t^3-100.70237523220045*\t^2+96.57681383070323*\t-26.691897009236246,-1.7938582228534286*\t^3+5.381574668560286*\t^2-11.696028058489805*\t+8.858311612782947});
\draw [line width=2pt] (1.08,0.8)-- (1.6375358166189113,1.1151289398280801);
\draw [line width=2pt] (2.1041547277936967,1.115042979942693)-- (2.75,0.75);
\draw[line width=2pt, dashed, smooth,samples=100,domain=0:0.24923573566606141, variable=\t] plot({15.093577207246362*\t^3-2.211627724160939*\t+1.6375358166189111,2.4113417606867196*\t^3+3.079567751621067*\t+1.1151289398280804});
\draw[line width=2pt, dashed, smooth,samples=100,domain=0.24923573566606141:0.5018039092849544, variable=\t] plot({-15.625101623802697*\t^3+22.968577551437942*\t^2-7.936218047396558*\t+2.113126643485062,-23.977266644473378*\t^3+19.730952687191067*\t^2-1.8380907567632496*\t+1.523681018525289});
\draw[line width=2pt, dashed, smooth,samples=100,domain=0.5018039092849544:0.7467672776563844, variable=\t] plot({-14.831796000032488*\t^3+21.77432596164105*\t^2-7.336937930966705*\t+2.012886275091314,25.747186823666738*\t^3-55.124822724920506*\t^2+35.72482997759091*\t-4.759392471027976});
\draw[line width=2pt, dashed, smooth,samples=100,domain=0.7467672776563844:1, variable=\t] plot({15.0762161875681*\t^3-45.22864856270431*\t^2+42.69869094945876*\t-10.442103846528852,-3.3653458977419044*\t^3+10.096037693225716*\t^2-12.979974403275182*\t+7.364325587734064});
\draw (-0.2,2.82) node[anchor=north west] {$\kt_\psi$};
\begin{scriptsize}
\draw [fill=ududff] (1.9,0) circle (2.5pt);
\end{scriptsize}
\end{tikzpicture}
\end{figure}

If $\psi$ corresponds the half-edges $\kh(a_1)$ and $\kh(a_2)$ then we can take $\kt_\psi^+$ and $\kt_\psi^-$ to be such that $\kt_\psi^+\cap\kn_v\subset R_{a_1}\cup R_{a_2-1}$ and $\kt_\psi^-\cap\kn_v\subset R_{a_1-1}\cup R_{a_2}$, as $\Sigma$ is an oriented surface.

Now we assume there exists $j$ such that $\kh(j)$ corresponding to $\gamma_{S_1}$ and $\kh(j-1)$ corresponding to $\gamma_{S_2}$. There exist a chord of $\kn_v$ which has end points on $\gamma_{S_1}$ and $\gamma_{S_2}$ and these are the only points where it intersects any arc in $\ka_{\kS_1}$. This implies the chord has end points on two half-edges of $\Gamma(\kS_1)$, one that corresponds to $S_1$ and the other that corresponds to $S_2$, which bound some region $R_{a^\prime}$. Therefore, $\gamma_{S_1}^\prime\gamma_{S_2}^\prime$ an simple concatenation with respect to $\ka_{\kS_1}$.

We have the following diagram of $R_{a^\prime}$.

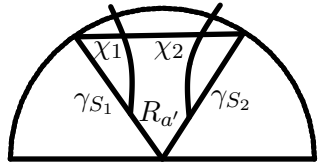
\begin{figure}[ht]
    \centering
    \begin{tikzpicture}[line cap=round,line join=round,>=triangle 45,x=1cm,y=1cm]
\draw [line width=2pt] (0,0)-- (4,0);
\draw [line width=2pt] (0.86,1.62)-- (2,0);
\draw [line width=2pt] (3.06,1.66)-- (2,0);
\draw (1.55,0.9) node[anchor=north west] {$R_{a^\prime}$};
\draw[line width=2pt, smooth,samples=100,domain=0:0.47830801385669, variable=\t] plot({-0.522997976145041*\t^3+0.12158060152918314*\t+1.5790769230769233,-0.09118303444667376*\t^3+1.4882125830900754*\t+0.5981538461538461});
\draw[line width=2pt, smooth,samples=100,domain=0.47830801385669:1, variable=\t] plot({0.4795053975628543*\t^3-1.4385161926885632*\t^2+0.8096344245547372*\t+1.4693763705709715,0.083600241640733*\t^3-0.25080072492219907*\t^2+1.6081725797014306*\t+0.5790279035800352});
\draw[line width=2pt, smooth,samples=100,domain=0:0.4469485666817181, variable=\t] plot({0.7897011885194979*\t^3-0.14040468716707832*\t+2.3522461538461537,-0.19526176207999926*\t^3+1.6239070806186102*\t+0.5516307692307689});
\draw[line width=2pt, smooth,samples=100,domain=0.4469485666817181:1, variable=\t] plot({-0.6381970881043038*\t^3+1.9145912643129113*\t^2-0.9961285085330727*\t+2.479734332324465,0.1578008109766117*\t^3-0.4734024329298351*\t^2+1.8354936194802383*\t+0.5201080024729855});
\draw [line width=2pt] (0.86,1.62)-- (3.06,1.66);
\draw [shift={(2,0)},line width=2pt,dash pattern=on 1pt off 1pt]  plot[domain=0:3.141592653589793,variable=\t]({1*2*cos(\t r)+0*2*sin(\t r)},{0*2*cos(\t r)+1*2*sin(\t r)});
\draw (0.95,1.7) node[anchor=north west] {$\chi_1$};
\draw (1.75,1.7) node[anchor=north west] {$\chi_2$};
\draw (0.7,1) node[anchor=north west] {$\gamma_{S_1}$};
\draw (2.5,1.05) node[anchor=north west] {$\gamma_{S_2}$};
\end{tikzpicture}
\end{figure}

In this diagram, $\chi_1$ and $\chi_2$ are sections in the tubular neighbourhoods of $\restr{\gamma_{S_1}}{[t_3,t_4]}$ and $\restr{\gamma_{S_2}}{[t_5,t_6]}$ respectively. We have the chord $\chi_3$ with end points $\chi_1(t_7)$ and $\chi_2(t_8)$ for some $t_3\leq t_7\leq t_4$ and $t_5\leq t_8\leq t_6$. We can orientate $\gamma_{S_1}$ and $\gamma_{S_2}$ such that $\gamma_{S_1}^\prime \gamma_{S_2}^\prime\simeq\restr{\gamma_{S_1}}{[0,t_3]}\restr{\chi_1}{[t_3,t_7]}\chi_3\restr{\chi_2}{[t_8,t_6]}\restr{\gamma_{S_2}}{[t_6,1]}$. 

Now if $\kh(a_3)\not=\kh(a_1)$ is the other half-edge that corresponds to $S_1$, we consider the diagram of $R_{a_3}$ (Figure~\ref{ra3}).
\begin{figure}[ht]
    \centering
    \begin{tikzpicture}[line cap=round,line join=round,>=triangle 45,x=1cm,y=1cm]
\draw [line width=2pt] (0,0)-- (4,0);
\draw [line width=2pt] (0.86,1.62)-- (2,0);
\draw [line width=2pt] (3.06,1.66)-- (2,0);
\draw (0.96,1.9) node[anchor=north west] {$R_{a_3}$};
\draw[line width=2pt, smooth,samples=100,domain=0:0.5651700159225248, variable=\t] plot({-0.13428979842206168*\t^3+0.5166525169092405*\t+2.3522461538461537,0.03358516835073817*\t^3+1.5611342491588138*\t+0.5516307692307689});
\draw[line width=2pt, smooth,samples=100,domain=0.5651700159225248:1, variable=\t] plot({0.17454308647424482*\t^3-0.5236292594227345*\t^2+0.8125920737946872*\t+2.2964940991538025,-0.04365230280017996*\t^3+0.13095690840053986*\t^2+1.4871213311529161*\t+0.5655740632467242});
\draw [shift={(2,0)},line width=2pt,dash pattern=on 1pt off 1pt]  plot[domain=0:3.141592653589793,variable=\t]({1*2*cos(\t r)+0*2*sin(\t r)},{0*2*cos(\t r)+1*2*sin(\t r)});
\draw (2.65,2.4) node[anchor=north west] {$\chi_1$};
\draw (2.7,1.34) node[anchor=north west] {$\gamma_{S_1}$};
\draw (1.8,1.58) node[anchor=north west] {$\chi_4$};
\draw [line width=2pt] (2,0)-- (2.7037515881699083,1.8320875890030819);
\end{tikzpicture}
\caption{}
\label{ra3}
\end{figure}
We know from above that $\chi_1\cap R_{a_3}$ is non-empty. There exists a path $\chi_4$ which is a radius of $\kn_v$ in $R_{a_3}$ and intersects $\chi_1$ at the point $\chi_1(t_9)$ on the boundary of $\kn_v$ and $\restr{\gamma_{S_1}}{[0,t_3]}\restr{\chi_1}{[t_3,t_9]}\simeq\chi_4$. 

Similarly, we have the radius $\chi_5$ of $\kn_v$ in $R_{a_4}$, (where $\kh(a_4+1)\not=\kh(a_2)$ corresponds to $\gamma_{S_2}$) such that $\gamma_{S_1}^\prime\gamma_{S_2}^\prime\simeq\chi_4\restr{\chi_1}{[t_9, t_7]}\chi_3\restr{\chi_2}{[t_8,t_{10}]}\chi_5$. As $\chi_4$ and $\chi_5$ are radii of $\kn_v$, the tubular neighbourhoods of $\restr{\gamma_{S_1}}{[t_3,t_4]}$ and $\restr{\gamma_{S_2}}{[t_5,t_6]}$ do not contain subpaths of any other arc in $\ka_{\kS_1}$ and $\chi_3$ is a contained in $R_{a^\prime}$ then $\rho=\chi_4\restr{\chi_1}{[t_9, t_7]}\chi_3\restr{\chi_2}{[t_8,t_{10}]}\chi_5$ only intersects any arc in $\ka_{\kS_1}$ at $v$. Since $\chi_4$ and $\chi_5$ are radii of $\kn_v$ and $\chi_3$ is a chord then we can take the radius of $\kn_v$ to be small enough such that $\kn_v\cap\rho\subset R_{a_3}\cup R_{a_4}$. Therefore, when $\gamma_{S_1}^\prime\gamma_{S_2}^\prime$ is in minimal position $\gamma_{S_1}^\prime\gamma_{S_2}^\prime\cap\kn_v\subset R_{a_3}\cup R_{a_4}$ and the result follows. \qed

\cor \label{sscreg} Let $\{\gamma_{S_1},\gamma_{S_2},...,\gamma_{S_n}\}$ be a set of closed arcs in $\ka_{\kS_1}$ such that there exists a sequence of integers $\{a_1,a_2,...,a_{n+1}\}$ with $\kh(a_i)$ corresponding to $S_{i}$ for all $1\leq i\le n$ and $\kh(a_j+1)$ being the other half-edge corresponding to $S_{j-1}$ for all $1\leq j\leq n+1$. Then the path $\phi$ homotopic to $\gamma^\prime_{S_1}\gamma^\prime_{S_2}...\gamma^\prime_{S_n}$ is a sequence of simple concatenations with respect to $\ka_{\kS_1}$ and $\kn_v\cup\phi\subset R_{a_1}\cup R_{a_{n+1}}$.

\pf We prove the result by induction, the base case comes from Lemma~\ref{lsscreg} and we assume it is true for $n=k$. By Definition~\ref{ssc}, $\ka_{\kS_1}\cup\{\gamma^\prime_{S_1}\gamma^\prime_{S_2}...\gamma^\prime_{S_{k}}\}$ is an arc-collection and applying Lemma~\ref{lsscreg} to two arcs in this arc-collection, the result follows for $n=k+1$.  \qed
\section{Thick subcategories generated by exceptional and spherelike objects}
We show that all thick subcategories in $\tds$ are generated by exceptional and spherelike objects. First, we consider the thick subcategories generated by one string object. By looking at the number of internal intersection points of arcs corresponding to objects generated in $\sT$, we prove the following:
\lem \label{strgen} Let $\sT$ be the thick subcategory generated by a string object $B\in\sD$. Then $\sT$ is generated by a finite set of exceptional and spherelike objects $\kS=\{A_1, A_2,...,A_n\}$. Furthermore, $\gamma_B$ is homotopic to a concatenation of the form $\gamma^\prime_{A_1}\gamma^\prime_{A_2}...\gamma^\prime_{A_n}$.

\pf Since $\sD$ is a bounded derived category of a finite dimensional algebra then $\dhomd(B,B)$ is finite and $\gamma_B$ (in minimal position) has a finite number of self-intersection points. If $B$ was exceptional or spherelike the proof would be trivial, therefore we assume it is not. 
 
By Lemma~\ref{excsph}, $\gamma_B$ has an internal intersection point. Therefore, there exists $t_1\in(0,1)$ which is minimal such that $p=\gamma_B(t_1)=\gamma_B(t_2)$ for some $t_2\in(0,t_1)$. By \cite[Theorem 4.1]{Opper} there exists two objects $M_1$ and $N_1$ corresponding to the arcs $\restr{\gamma_B}{[0,t_1]}(\restr{\gamma_B}{[0,t_2]})\inv$ and $\restr{\gamma_B}{[0,t_2]}\restr{\gamma_B}{[t_1,1]}$ respectively which are direct summands of cones of morphisms corresponding to $p$. We see these paths in a neighbourhood of $p$.

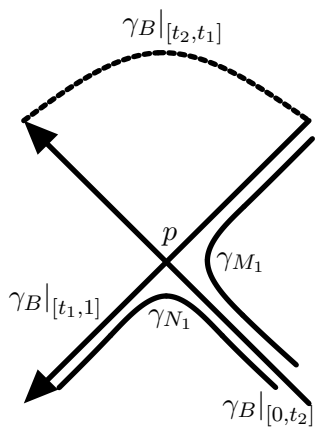
\begin{figure}[ht]
    \centering
    \begin{tikzpicture}[line cap=round,line join=round,>=triangle 45,x=.75cm,y=.75cm]
\draw [->,line width=2pt] (-1,-6) -- (-6,-1);
\draw [->,line width=2pt] (-1,-1) -- (-6,-6);
\draw[line width=2pt, smooth,samples=100,domain=0:0.3213187959380341, variable=\t] plot({-2.1416468011261296*\t^3-3.5134928208247893*\t-0.9399999999999995,-1.6903506553817538*\t^3-3.560087238268793*\t-1.3199999999999996});
\draw[line width=2pt, smooth,samples=100,domain=0.3213187959380341:0.5366330268884635, variable=\t] plot({35.56814967038461*\t^3-36.350599191882466*\t^2+8.16663794313696*\t-2.1910151844916563,1.4835516030183598*\t^3-3.0595033562823946*\t^2-2.5770113036597597*\t-1.4252935918747445});
\draw[line width=2pt, smooth,samples=100,domain=0.5366330268884635:0.7562008707794534, variable=\t] plot({-34.19495282743511*\t^3+75.96095538372289*\t^2-52.10345154331899*\t+8.589958332826816,5.45616189102634*\t^3-9.455005006788314*\t^2+0.8550261055213951*\t-2.0392084663025187});
\draw[line width=2pt, smooth,samples=100,domain=0.7562008707794534:1, variable=\t] plot({2.2064666303265232*\t^3-6.61939989097957*\t^2+10.34388502468764*\t-7.150951764034594,-3.9962654527380215*\t^3+11.988796358214065*\t^2-15.360795159515034*\t+2.048264254038989});
\draw[line width=2pt, smooth,samples=100,domain=0:0.2568395313984434, variable=\t] plot({2.267608961917259*\t^3+3.3545471478822066*\t-5.380000000000001,3.52221945872979*\t^3+3.349654461031808*\t-5.720000000000001});
\draw[line width=2pt, smooth,samples=100,domain=0.2568395313984434:0.49717792469183864, variable=\t] plot({-2.5215250610468227*\t^3+3.6901168147873102*\t^2+2.4067792743667176*\t-5.298858581163929,-30.926672946592507*\t^3+26.54351214773512*\t^2-3.4677687606613703*\t-5.136338738131821});
\draw[line width=2pt, smooth,samples=100,domain=0.49717792469183864:0.7375163179852338, variable=\t] plot({-1.9321645766784767*\t^3+2.811065747246417*\t^2+2.843824059824844*\t-5.371288254307748,31.203101676788204*\t^3-66.12514507873713*\t^2+42.60504192317546*\t-12.771800204968475});
\draw[line width=2pt, smooth,samples=100,domain=0.7375163179852338:1, variable=\t] plot({1.8590907626674875*\t^3-5.5772722880024626*\t^2+9.030360241597087*\t-6.892178716262112,-3.699588578932568*\t^3+11.098765736797704*\t^2-14.34885244191785*\t+1.2296752840527148});
\draw (-4.503333333333333,1.15) node[anchor=north west] {$\left.\gamma_B\right|_{[t_2,t_1]}$};
\draw (-6.5,-3.7) node[anchor=north west] {$\left.\gamma_B\right|_{[t_1,1]}$};
\draw (-2.75,-5.6) node[anchor=north west] {$\left.\gamma_B\right|_{[0,t_2]}$};
\draw (-2.8,-3.15) node[anchor=north west] {$\gamma_{M_1}$};
\draw (-4,-4.15) node[anchor=north west] {$\gamma_{N_1}$};
\draw[line width=2pt,dotted, smooth,samples=100,domain=0:0.47671301884549466, variable=\t] plot({0.11309487940970797*\t^3+4.9178738082542575*\t-6,-5.045419661981077*\t^3+3.663836106111361*\t-0.9999999999999996});
\draw[line width=2pt,dotted, smooth,samples=100,domain=0.47671301884549466:1, variable=\t] plot({-0.10302912803299898*\t^3+0.30908738409899694*\t^2+4.770527828293368*\t-5.976586084359366,4.596363611223211*\t^3-13.789090833669633*\t^2+10.237275224564751*\t-2.0445480021183275});
\draw (-3.75,-2.7) node[anchor=north west] {$p$};
\end{tikzpicture}
\end{figure}

First we show that $\gamma_{M_1}$ has no internal self-intersection points. Since $t_1$ is minimal, there exists $\restr{\gamma_B}{[0,t_3]}$ for some $t_3\in(t_2,t_1)$ with no self-intersection points. Let $\kt$ be a tubular neighbourhood of $\restr{\gamma_B}{[0,t_3]}$. Since the end points of $\restr{\gamma_B}{[t_3,t_1]}$ are the only points where it intersects $\restr{\gamma_B}{[0,t_3]}$ and $\restr{\gamma_B}{[0,t_2]}$ is a subpath of the zero section of $\kt$ then by making the tubular neighbourhood smaller if necessary, we may assume that $\kt\cap\restr{\gamma_B}{[t_3,t_1]}(\restr{\gamma_B}{[0,t_2]})\inv$ is a subspace of either $\kt^+$ or $\kt^-$. Using Lemma~\ref{tbnhnsip} there exists a path $\psi$ homotopic to $\restr{\gamma_B}{[0,t_3]}$ that only intersects $\restr{\gamma_B}{[t_3,t_1]}(\restr{\gamma_B}{[0,t_2]})\inv$ at its end points $\gamma_B(0)$ and $\gamma_B(t_3)$, which are also the end points of $\restr{\gamma_B}{[t_3,t_1]}(\restr{\gamma_B}{[0,t_2]})\inv$. We see this below.

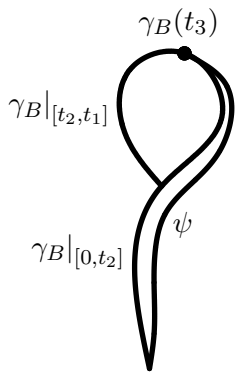
\begin{figure}[ht]
    \centering
    \begin{tikzpicture}[line cap=round,line join=round,>=triangle 45,x=1cm,y=1cm]
\draw[line width=2pt, smooth,samples=100,domain=0:0.3431550347195055, variable=\t] plot({18.177455700061426*\t^3-1.674231609648244*\t+2.5799999999999983,-6.841200575113684*\t^3+7.8577941739074175*\t-3.919999999999998});
\draw[line width=2pt, smooth,samples=100,domain=0.3431550347195055:0.5221505444546778, variable=\t] plot({-113.21734891708056*\t^3+135.26636622107398*\t^2-48.09156620662223*\t+7.889447355070496,40.40630561478685*\t^3-48.63965888101611*\t^2+24.548738005967405*\t-5.829193803397271});
\draw[line width=2pt, smooth,samples=100,domain=0.5221505444546778:0.653241497116008, variable=\t] plot({66.2481515128927*\t^3-145.85766005995123*\t^2+98.69749717530621*\t-17.659215766551522,-161.4569742646436*\t^3+267.56940560233875*\t^2-140.55999717552064*\t+22.908011519681786});
\draw[line width=2pt, smooth,samples=100,domain=0.653241497116008:0.7597643491972297, variable=\t] plot({176.46200952324776*\t^3-361.84645688879755*\t^2+239.79034217606704*\t-48.381782866769534,71.87002119188355*\t^3-189.68722188646643*\t^2+158.13900673148322*\t-42.132849980075385});
\draw[line width=2pt, smooth,samples=100,domain=0.7597643491972297:0.8626155270402535, variable=\t] plot({-62.5447759757892*\t^3+182.9200476263957*\t^2-174.10382659136846*\t+56.43889505660346,122.68624605405317*\t^3-305.51228991966565*\t^2+246.1387641664517*\t-64.4192094591064});
\draw[line width=2pt, smooth,samples=100,domain=0.8626155270402535:1, variable=\t] plot({-51.106607359376014*\t^3+153.31982207812803*\t^2-148.57021242953917*\t+49.096997710787136,-29.069011521087607*\t^3+87.20703456326281*\t^2-92.62702290128192*\t+32.9889998591067});
\draw[line width=2pt, smooth,samples=100,domain=0:0.1500983135278404, variable=\t] plot({-38.45825584676493*\t^3+4.597333439891937*\t+3.04,-23.30464402490219*\t^3-2.4063699526104907*\t+0.2600000000000002});
\draw[line width=2pt, smooth,samples=100,domain=0.1500983135278404:0.3325723079988216, variable=\t] plot({18.439738029891675*\t^3-25.62087877201064*\t^2+8.442984134671981*\t+2.8475914387654493,35.54013191735694*\t^3-26.49750488657019*\t^2+1.5708608435594054*\t+0.06100812166130316});
\draw[line width=2pt, smooth,samples=100,domain=0.3325723079988216:0.43727522178766215, variable=\t] plot({16.427351449100097*\t^3-23.613086622731505*\t^2+7.775248065604305*\t+2.9216149473067485,-27.95631840286325*\t^3+36.85397821161423*\t^2-19.49808810555213*\t+2.396657781366427});
\draw[line width=2pt, smooth,samples=100,domain=0.43727522178766215:0.5657518241313076, variable=\t] plot({36.59576665340038*\t^3-50.070531317429705*\t^2+19.344433062413266*\t+1.2353089688460408,-18.15004594575648*\t^3+23.989858320839417*\t^2-13.872927227210495*\t+1.57674329114372});
\draw[line width=2pt, smooth,samples=100,domain=0.5657518241313076:0.6618727305798515, variable=\t] plot({-25.492149423131135*\t^3+55.30852401299805*\t^2-40.27395971601023*\t+12.478380457569331,25.325130418381377*\t^3-49.7986226764845*\t^2+27.873040496903844*\t-6.2958758388710745});
\draw[line width=2pt, smooth,samples=100,domain=0.6618727305798515:0.7589134752732688, variable=\t] plot({-10.543922896157714*\t^3+25.627053486795845*\t^2-20.628603771207395*\t+8.144138663436074,-2.4290721686946997*\t^3+5.310626877638763*\t^2-8.602268985690186*\t+1.7514617231254799});
\draw[line width=2pt, smooth,samples=100,domain=0.7589134752732688:0.9022296153736542, variable=\t] plot({-13.772352851465271*\t^3+32.97735047797219*\t^2-26.206843205071685*\t+9.555272354989523,1.2397255092664219*\t^3-3.04227310992916*\t^2-2.2631406275149706*\t+0.147845078956786});
\draw[line width=2pt, smooth,samples=100,domain=0.9022296153736542:1, variable=\t] plot({14.660620003917373*\t^3-43.981860011752126*\t^2+43.228035674532364*\t-11.32679566669761,-1.0680742757675425*\t^3+3.204222827302628*\t^2-7.8989142543967*\t+1.8427657028616151});
\draw (2.3,1) node[anchor=north west] {$\gamma_B(t_3)$};
\draw (2.75,-1.7) node[anchor=north west] {$\psi$};
\draw (0.85,-2) node[anchor=north west] {$\left.\gamma_B\right|_{[0,t_2]}$};
\draw (0.55,-0.1) node[anchor=north west] {$\left.\gamma_B\right|_{[t_2,t_1]}$};
\begin{scriptsize}
\draw [fill=black] (3.04,0.26) circle (2.5pt);
\draw [fill=black] (3.04,0.26) circle (2.5pt);
\draw [fill=black] (3.04,0.26) circle (2.5pt);
\end{scriptsize}
\end{tikzpicture}
\end{figure}

Therefore, $\gamma_{M_1}\simeq\psi\restr{\gamma_B}{[t_3,t_1]}(\restr{\gamma_B}{[0,t_2]})\inv$ with the only self-intersection point being its end point. By Lemma~\ref{cosm} $\gamma_{M_1}$ is homotopic to a smooth path, which is a closed arc that has no internal self-intersection points, hence $M_1$ is spherelike by Lemma~\ref{excsph}.

Now consider $N_1$, the arc $\gamma_{N_1}$ is the concatenation at $\gamma_B(t_1)$ of two paths which only have one end point at $\gamma_B(t_1)$. This is the concatenation of two open paths $\restr{\gamma_B}{[0,t_2]}$ and $\restr{\gamma_B}{[t_1,1]}$, which is always simple. By Lemma~\ref{simcc} the number of internal self-intersection points of $\restr{\gamma_B}{[t_1,1]}$ gives an upper bound of the number of internal self-intersection points of $\gamma_{N_1}$. Each internal self-intersection point of $\restr{\gamma_B}{[t_1,1]}$ corresponds to a unique internal self-intersection point of $\gamma_B$. As $\gamma_B(t_1)$ is an end point of $\restr{\gamma_B}{[t_1,1]}$ it is an internal self-intersection point of $\gamma_B$ which does not correspond to an self-intersection point of $\restr{\gamma_B}{[t_1,1]}$. Hence, $\gamma_{N_1}$ has fewer internal self-intersection points than $\gamma_B$.

Finally, $\gamma_{M_1}$ and $\gamma_{N_1}$ have a common end point and  there exists a morphism corresponding to this point with the cone corresponding to  $\gamma_{M_1}\gamma_{N_1}\simeq\restr{\gamma_B}{[0,t_1]}(\restr{\gamma_B}{[0,t_2]})\inv\restr{\gamma_B}{[0,t_2]}\restr{\gamma_B}{[t_1,1]}\simeq\gamma_B$ by \cite[Theorems 3.7 and 4.1]{Opper}. Therefore $\td(B)=\td(M_1,N_1)$ where $M_1$ is a spherelike object, $\gamma_{N_1}$ has fewer internal self-intersection points than $\gamma_B$ and $\gamma_B\simeq\gamma_{M_1}\gamma_{N_1}$.

Working inductively, we can write $\td(N_i)=\td(M_{i+1}, N_{i+1})$ where $M_{i+1}$ is spherelike, $\gamma_{N_{i+1}}$ has fewer self-intersections than $\gamma_{N_i}$ and $\gamma_{N_i}\simeq\gamma_{M_{i+1}}\gamma_{N_{i+1}}$. Since $\gamma_B$ has a finite number of self-intersection points, there exists a finite $k\in\IN$ such that $\gamma_{N_k}$ has no internal self-intersection points. By Lemma~\ref{excsph}, $N_k$ is an exceptional or spherelike object, $\td(B)=\td(M_1,...,M_k, N_k)$ and $\gamma_B\simeq\gamma_{M_1}\gamma_{M_2}...\gamma_{M_k}\gamma_{N_k}$ as required.  \qed

We apply the above to more general situations and see what it tells us about all thick subcategories in $\tds$.

\thm\label{ctwsgen} Let $\sT$ be a Ext-connected thick subcategory containing a string object. Then $\sT$ is generated by exceptional and spherelike string objects. Furthermore, if $\sT\in\tds$ then it is generated by finitely many exceptional and spherelike string objects.

\pf By Theorem~\ref{cstobj}, $\sT$ is generated by string objects. It follows from the definition of thick subcategories and Lemma~\ref{strgen} each string object can be replaced by a finite number of exceptional and spherelike string objects. The result follows.\qed

So now we just have to consider thick subcategories generated by finitely many exceptional and spherelike string objects.

\section{Thick subcategories corresponding to arc-collections}
In this section we first show that every thick subcategory in $\tds$ corresponds to an arc-collection.

\defn Let $\sT$ be a thick subcategory generated by a finite set $\kS$ of string objects. Suppose the objects in $\kS$ correspond to a set of arcs which form a (pointed) arc-collection. Then we say $\sT$ corresponds to a finite (pointed) arc-collection.

Our next result says if we would like to classify the Ext-connected thick subcategories generated by finitely many exceptional and spherelike string objects then we only need to consider thick subcategories that correspond to finite arc-collections.

\thm\label{strESc} $\sT\in \tds$ if and only if $\sT$ corresponds to an connected finite arc-collection.

\pf ($\impliedby$) If $\sT$ corresponds to a connected finite arc-collection $\ka$, then $\sT$ is generated by a finite set of string objects by definition. It remains to show that $\sT$ is Ext-connected. %Let the objects in $\sT_{\leq 0}$ be all of the shifts of the objects corresponding to arcs in $\ka$ and $\sT_{\leq v}$ be the set of objects in $\sT$ generated in at most $v$ steps by taking direct summands and cones between objects generated in fewer steps. The collection of objects $T_{\leq v}$ is Ext-connected if for any two direct summands $Y$ and $Y^\prime$ of objects in $\sT_{\leq v}$ there exists a finite sequence of indecomposable objects in $\sT_{\leq v}$\[Y=Y_0,Y_1,...,Y_a=Y^\prime\] with morphisms between shifts of consecutive objects in the sequence. 
As the union of images of arcs in $\ka$ is path connected then for arcs $\beta_1, \beta_2\in\ka$ there exists a path which is a (possibly trivial) finite sequence of arcs in $\ka$, $\delta_1\delta_2...\delta_n$ with $\delta_1$ and $\delta_n$ having a common end point with $\beta_1$ and $\beta_2$ respectively. As $\delta_k$ has a common end point with $\delta_{k+1}$ there exists a morphism between shifts of their corresponding objects (using the correspondence between morphisms and intersection points of the geometric model \cite[Theorem 3.3]{Opper}). We have the sequence of objects corresponding to arcs in $\ka$\[X_0,X_1,...,X_{n+1}\] where there exists morphisms between consecutive objects, $X_0$ corresponds to $\beta_1$, $X_{k^\prime}$ corresponds to $\delta_{k^\prime}$ for $1\leq k^{\prime}\leq n$ and $X_{n+1}$ corresponds to $\beta_2$. Hence, $\sT_{\leq0}$ is Ext-connected. It follows from Lemma~\ref{tsecgsec} that $\sT$ is Ext-connected.

($\implies$) Assume $\sT\in\tds$ i.e. $\sT$ is Ext-connected and finitely generated by string objects. Applying Theorem~\ref{ctwsgen} we deduce that $\sT$ is generated by a set $\kS=\{D_1,...,D_m\}$ of finitely many exceptional or spherelike string objects. Let $\ku_1=\{\alpha_1,...,\alpha_m\}$ be the set of arcs in the geometric model, such that $\alpha_i$ corresponds to $D_i$ for all $1\leq i\leq m$.

Since $\sD$ is hom-finite there exist a finite number of intersection points between any two arcs in $\ku_1$. Since there are a finite number of objects in $\kS$ this implies that there is a finite number $x$ of internal intersection points between all of the arcs in $\ku_1$. Take any $ D_i$, such that $\alpha_i$ internally intersects some other arc in $\ku_1$. Let $\alpha_j\in\ku_1$ be an arc and $p_1=\alpha_i(t_1)=\alpha_j(t_2)$ be the internal intersection point with $\alpha_i$ such that $t_1\in(0,1)$ is minimized for any $\alpha_k\in \ku_1$. We define subpaths $\phi_1=\restr{\alpha_i}{[0,t_1]}$, $\phi_2=\restr{\alpha_i}{[t_1,1]}$ of $\alpha_i$ and $\psi_1=\restr{\alpha_j}{[0,t_2]}$, $\psi_2=\restr{\alpha_j}{[t_2,1]}$ of $\alpha_j$. The point $p_1$ corresponds to two morphisms between $D_i$ and shifts of $D_j$ and by \cite[Theorem 4.1]{Opper} the cones of these are direct sums of objects corresponding to the arcs $\phi_1\psi_2$, $\phi_1\psi_1^{-1}$, $\phi\inv_2\psi_2$ and $\psi_1\phi_2$. In Figure~\ref{figtsfsc}, we see these in a neighbourhood of $p_1$.

\begin{figure}[ht]
    \centering
\begin{tikzpicture}[line cap=round,line join=round,>=triangle 45,x=1cm,y=.95cm]
\draw[<-,line width=2pt, smooth,samples=100,domain=0:0.3213187959380341, variable=\t] plot({-2.1416468011261296*\t^3-3.5134928208247893*\t-0.9399999999999995,-1.6903506553817538*\t^3-3.560087238268793*\t-1.3199999999999996});
\draw[line width=2pt, smooth,samples=100,domain=0.3213187959380341:0.5366330268884635, variable=\t] plot({35.56814967038461*\t^3-36.350599191882466*\t^2+8.16663794313696*\t-2.1910151844916563,1.4835516030183598*\t^3-3.0595033562823946*\t^2-2.5770113036597597*\t-1.4252935918747445});
\draw[line width=2pt, smooth,samples=100,domain=0.5366330268884635:0.7562008707794534, variable=\t] plot({-34.19495282743511*\t^3+75.96095538372289*\t^2-52.10345154331899*\t+8.589958332826816,5.45616189102634*\t^3-9.455005006788314*\t^2+0.8550261055213951*\t-2.0392084663025187});
\draw[line width=2pt, smooth,samples=100,domain=0.7562008707794534:1, variable=\t] plot({2.2064666303265232*\t^3-6.61939989097957*\t^2+10.34388502468764*\t-7.150951764034594,-3.9962654527380215*\t^3+11.988796358214065*\t^2-15.360795159515034*\t+2.048264254038989});
\draw[<-,line width=2pt, smooth,samples=100,domain=0:0.2568395313984434, variable=\t] plot({2.267608961917259*\t^3+3.3545471478822066*\t-5.380000000000001,3.52221945872979*\t^3+3.349654461031808*\t-5.720000000000001});
\draw[line width=2pt, smooth,samples=100,domain=0.2568395313984434:0.49717792469183864, variable=\t] plot({-2.5215250610468227*\t^3+3.6901168147873102*\t^2+2.4067792743667176*\t-5.298858581163929,-30.926672946592507*\t^3+26.54351214773512*\t^2-3.4677687606613703*\t-5.136338738131821});
\draw[line width=2pt, smooth,samples=100,domain=0.49717792469183864:0.7375163179852338, variable=\t] plot({-1.9321645766784767*\t^3+2.811065747246417*\t^2+2.843824059824844*\t-5.371288254307748,31.203101676788204*\t^3-66.12514507873713*\t^2+42.60504192317546*\t-12.771800204968475});
\draw[line width=2pt, smooth,samples=100,domain=0.7375163179852338:1, variable=\t] plot({1.8590907626674875*\t^3-5.5772722880024626*\t^2+9.030360241597087*\t-6.892178716262112,-3.699588578932568*\t^3+11.098765736797704*\t^2-14.34885244191785*\t+1.2296752840527148});
\draw (-5.35,-3.1) node[anchor=north west] {$\psi_1\phi_2$};
\draw (-6.6,-5.4) node[anchor=north west] {$\psi_1$};
\draw (-1.2,-5.4) node[anchor=north west] {$\phi_1$};
\draw (-2.76,-3.1) node[anchor=north west] {$\phi_1\psi_2$};
\draw (-4.25,-4.25) node[anchor=north west] {$\phi_1\psi_1^{-1}$};
\draw [line width=2pt] (-1,-6)-- (-6,-1);
\draw [line width=2pt] (-1,-1)-- (-6,-6);
\draw (-1.52,-0.4) node[anchor=north west] {$\psi_2$};
\draw (-4.15,-1.9) node[anchor=north west] {$\phi_2^{-1}\psi_2$};
\draw (-6.2,-0.4) node[anchor=north west] {$\phi_2$};
\draw[<-,line width=2pt, smooth,samples=100,domain=-0:0.32131879593803414, variable=\t] plot({2.1416468011260017*\t^3+3.513492820824797*\t-6.06,-1.6903506553817798*\t^3-3.560087238268787*\t-1.3200000000000005});
\draw[line width=2pt, smooth,samples=100,domain=0.32131879593803414:0.5366330268884635, variable=\t] plot({-35.56814967038393*\t^3+36.350599191881706*\t^2-8.16663794313671*\t-4.808984815508368,1.4835516030184803*\t^3-3.059503356282537*\t^2-2.577011303659708*\t-1.4252935918747491});
\draw[line width=2pt, smooth,samples=100,domain=0.5366330268884635:0.7562008707794535, variable=\t] plot({34.19495282743431*\t^3-75.96095538372127*\t^2+52.10345154331797*\t-15.589958332826612,5.4561618910262215*\t^3-9.45500500678807*\t^2+0.8550261055212396*\t-2.0392084663024876});
\draw[line width=2pt, smooth,samples=100,domain=0.7562008707794535:1, variable=\t] plot({-2.206466630326236*\t^3+6.619399890978709*\t^2-10.343885024686797*\t+0.1509517640343227,-3.9962654527379824*\t^3+11.98879635821395*\t^2-15.360795159514922*\t+2.048264254038953});
\draw[line width=2pt, smooth,samples=100,domain=-0.08:0.2568395313984434, variable=\t] plot({2.267608961917259*\t^3+3.3545471478822066*\t-5.380000000000001,-3.5222194587297766*\t^3-3.349654461031806*\t-1.2800000000000002});
\draw[line width=2pt, smooth,samples=100,domain=0.2568395313984434:0.49717792469183864, variable=\t] plot({-2.5215250610468227*\t^3+3.6901168147873102*\t^2+2.4067792743667176*\t-5.298858581163929,30.926672946592475*\t^3-26.543512147735086*\t^2+3.4677687606613636*\t-1.863661261868179});
\draw[line width=2pt, smooth,samples=100,domain=0.49717792469183864:0.7375163179852338, variable=\t] plot({-1.9321645766784767*\t^3+2.811065747246417*\t^2+2.843824059824844*\t-5.371288254307748,-31.203101676788204*\t^3+66.12514507873712*\t^2-42.60504192317545*\t+5.771800204968474});
\draw[->,line width=2pt, smooth,samples=100,domain=0.7375163179852338:1.08, variable=\t] plot({1.8590907626674875*\t^3-5.5772722880024626*\t^2+9.030360241597087*\t-6.892178716262112,3.6995885789325733*\t^3-11.09876573679772*\t^2+14.348852441917865*\t-8.229675284052718});
\end{tikzpicture}
    \caption{}
\label{figtsfsc}
\end{figure}
The cone of the morphism corresponding to the common end point of $\phi_1\psi_2$ and $\phi_1\psi^{-1}_1$ is $\psi_1\phi_1^{-1}\phi_1\psi_2\simeq\alpha_j$. Let $\ku_2=\ku_1\backslash\{\alpha_j\}\cup\{ \phi_1\psi_2, \phi_1\psi^{-1}_1\}$. By the definition of thick subcategories $\sT$ is generated by the objects corresponding to the arcs in $\ku_2$.% Note that, both $\phi_1\psi_2$ and $\phi_1\psi^{-1}_1$ have the end point $\phi_1(0)=\alpha_i(0)$, so they have a common end point with an arc in $\ku_1\backslash\{\alpha_j\}$.

%Since $D_i$ and $D_j$ are exceptional or spherelike the arcs $\alpha_i$ and $\alpha_j$ have no self-intersection points. The diagram above shows how the arcs $\phi_1\psi_2$ and $\phi_1\psi^{-1}_1$ differ from subpaths of $\alpha_i$ and $\alpha_j$. The arcs $\phi_1\psi_2$ and $\phi_1\psi^{-1}_1$ have no self-intersection points around $p_1$. As $t_1$ is minimal then $\phi_1$ does not intersect $\psi_1$ or $\psi_2$ apart from at common end points, so the arcs $\phi_1\psi_2$ and $\phi_1\psi^{-1}_1$ have no internal self-intersection points and their corresponding objects are exceptional and spherelike.
%Globally, the subpaths of $\phi_1\psi_2$ and $\phi_1\psi^{-1}_1$ outside of the neighbourhood of $p_1$, . So, they have no self-intersection points and their corresponding objects are exceptional or spherelike. 

%Since $t_1$ is minimal, the subpath $\phi_1$ has no intersection with any $\alpha_k$. Any intersection that $\psi_1$ and $\psi_2$ have with any $\alpha_k$ corresponds to an intersection point between $\alpha_j$ and $\alpha_k$. Since $p_1$ is a simple intersection point, then any intersection point between $\alpha_k$ and the arcs $\phi_1\psi_2$ or $\phi_1\psi^{-1}_1$ would correspond to an intersection between $\alpha_k$ and $\psi_2$ or $\psi_1$ respectively. Hence it corresponds to an intersection point between $\alpha_k$ and $\alpha_j$. The point $p_1$ does not correspond to an intersection point between the arcs in $V_2$ as we can deduce from the diagram that $\alpha_i$ does not intersect $\phi_1\psi_2$ or $\phi_1\psi^{-1}_1$ at $p_1$.

Now we need to show the arcs in $\ku_2$ have no internal self-intersection points and at most $x-1$ internal intersection points with each other. The arcs in $\ku_1\setminus\{\alpha_j\}$ correspond to spherelike and exceptional objects, so these arcs have no internal self-intersection points. Considering $\phi_1\psi_1\inv$ and $\phi_1\psi_2$, since the paths $\psi_1$ and $\psi_2$ only intersect $\phi_1$ at common end points then $\phi_1\psi_1\inv$ and $\phi_1\psi_2$ have no internal self intersection points. Therefore, an object corresponding to an arc in $\ku_2$ is exceptional or spherelike by Lemma~\ref{excsph}.

The internal intersection points between two arcs in $\ku_1\setminus\{\alpha_j\}$ correspond to distinct internal intersection point between the same two arcs in $\ku_1$. Since $\ku_1\setminus\{\alpha_j\}$ does not include $\alpha_j$, none of these points corresponds to $\alpha_i(t_1)$. Since the arcs in $\ku_1$ are assumed to be in minimal position and $t_1$ is minimal, any internal intersection point $p_2$ between an arc in $\ku_1\setminus\{\alpha_i, \alpha_j\}$ and another in $\{\phi_1\psi_1\inv, \phi_1\psi_2\}$ corresponds to an internal intersection point between the same arc in $\ku_1$ and $\alpha_j$. The point $p_2$ does not correspond to an internal intersection point between two arcs in $\ku_1\setminus\{\alpha_j\}$ or an internal intersection point between $\alpha_i$ and $\alpha_j$. 

The only internal intersection points between arcs in $\ku_2$ we have left to consider are those between $\phi_1\psi_1\inv$, $\phi_1\psi_2$ and $\alpha_i$. There exists a subpath $\restr{\alpha_j}{[s_1,s_2]}$ for some $0<s_1<t_2<s_2<1$ such that the only internal intersection point between $\alpha_j$ and any other arc in $\ku_1$ on $\restr{\alpha_j}{[s_1,s_2]}$ is $p_1$. Any internal intersection point between $\alpha_i$ and $\alpha_j$ which is not $p_1$ is an internal intersection point between $\phi_2$ and $\restr{\alpha_j}{[0,s_1]}$ or $\restr{\alpha_j}{[s_2,1]}$, so it corresponds to a distinct internal intersection point between $\alpha_i$ and $\alpha_j$ which is not $p_1$. When in minimal position, $\phi_1(\restr{\alpha_j}{[s_1,t_2]})\inv$ and $\phi_1\restr{\alpha_j}{[t_2,s_2]}$ only intersect $\phi_1$ at $\phi_1(0)$ by Lemma~\ref{simcc} as $\restr{\alpha_j}{[s_1,t_2]}$ and $\restr{\alpha_j}{[t_2,s_2]}$ only intersect at the point $p_1$. Furthermore, $\restr{\alpha_j}{[0,s_1]}$ and $\restr{\alpha_j}{[s_2, 1]}$ only possibly intersect at $\alpha_j(0)$ as $\alpha_j$ has no internal intersection points. Hence, any internal intersection points between $\phi_1\psi_1\inv$, $\phi_1\psi_2$ and $\alpha_i$ corresponds to a distinct internal intersection point between $\alpha_i$ and $\alpha_j$ which is not $p_1$.

%So when in minimal position $\phi_1\psi_1\inv$ and $\phi_1\psi_2$ only intersect at common end points and they only internally intersect $\alpha_i$ at a point which corresponds to an internal intersection point between $\phi_2$ and $\restr{\alpha_j}{[0,s_1]}$ or $\restr{\alpha_j}{[s_2,1]}$. So any internal intersection point between $\alpha_i$ and $\phi_1\psi_1\inv$ and $\phi_1\psi_2$ corresponds to a distinct internal intersection point between $\alpha_i$ and $\alpha_j$ which is not $p_1$. 

We have shown that any internal intersection point between arcs in $\ku_2$ when in minimal position corresponds to an distinct internal intersection point between the arcs in $\ku_1$ which is not $p_1$, Therefore, the arcs in $\ku_2$ have at most $x-1$ intersection points when in minimal postition. We proceed iteratively until we have constructed a set of arcs $\ku^\prime$ which have no internal intersection points and such that $\sT$ is generated by the objects corresponding to the arcs in $\ku^\prime$.  

If $\ku^\prime$ was an arc-collection which is not connected, then the set of objects corresponding to $\ku^\prime$ is a union of two sets $\kv_1$ and $\kv_2$ where no object in $\kv_1$ has a morphism to or from any object in $\kv_2$. %This implies $\kv_1\cup\kv_2$ is not Ext-connected. However, Lemma~\ref{tsecgsec} would imply the contradiction that $\sT$ is not Ext-connected. Therefore, $\ku^\prime$ must be a connected arc collection
Applying Lemma~\ref{zhom}, $\dhomd(\td(\kv_1),\td(\kv_2))=0$ and by Remark~\ref{rzhom} $\td(\kv_1)\cap\td(\kv_2)=0$. %This implies that $\td(\sV_1)\oplus\td(\sV_2)=\td(\sV_1,\sV_2)$. 
%Now consider $\td(\sV_1)\oplus\td(\sV_2)$. Let $V_1\oplus V_2$ be an object in $\td(\sV_1)\oplus\td(\sV_2)$ such that $V_1\in\td(\sV_1)$ and $V_2\in\td(\sV_2)$. Any shift of $V_1\oplus V_2$ is the direct sum of a shift of $V_1$ and a shift of $V_2$, hence it is in $\td(\sV_1)\oplus\td(\sV_2)$. A direct summand of $V_1\oplus V_2$ is a direct sum of a (possibly trivial) direct summand of $V_1$ and a (possibly trivial) direct summand of $V_2$, hence it is in $\td(\sV_1)\oplus\td(V_2)$. Let $V^*_1\oplus V^*_2$ be a object in $\td(\sV_1)\oplus\td(\sV_2)$ such that $V^*_1\in\td(\sV_1)$ and $V^*_2\in\td(\sV_2)$. As there are no morphisms between an object in $\td(\sV_1)$ and an object in $\td(\sV_2)$, any morphism from $V_1\oplus V_2$ to $V^*_1\oplus V^*_2$ has the form \[\begin{pmatrix}    f&0\\0&g\end{pmatrix}\] where $f:V_1\rightarrow V_1^*$ and $g:V_2\rightarrow V_2^*$. The cone of this morphism is $C(f)\oplus C(g)$, by the definition of thick subcategories $C(f)\in\td(\sV_1)$ and $C(g)\in\td(\sV_2)$. Hence, $C(f)\oplus C(g)\in\td(sV_1)\oplus\td(\sV_2)$. Therefore, $\td(\sV_1)\oplus\td(\sV_2)$ is closed under taking shifts, direct summands and cones of morphisms between objects. By definition $\td(\sV_1)\oplus\td(\sV_2)$ is a thick subcategory, in fact it is the smallest thick subcategory containing $\sV_1$ and $\sV_2$ as $\td(\sV_1)\cap\td(\sV_2)=0$. 
This implies $\sT=\td(\kv_1, \kv_2)=\td(\kv_1)\oplus\td(\kv_2)$ and contradicts the fact that $\sT$ is Ext-connected. We conclude that $\ku^\prime$ must be a connected arc collection and the result follows.\qed

We assume, for the remainder of this paper, that we have a thick subcategory $\sT=\td(\kS)$ where $\ka_\kS$ is an arc-collection.

We consider the ribbon graph $\Gamma(\overline{\kS})$ corresponding to an arc-collection $\ka_{\overline{\kS}}$. %Let us introduce some terminology for talking about a half-edge in relation to another. We now define a permutation of half-edges around each vertex of the ribbon graph. 
Let $v$ be a vertex in $\Gamma(\overline{\kS})$ and $H_{v}$ be the set of half-edges adjacent to $v$. As $\Gamma(\overline{\kS})$ is an ribbon graph then we take the radius of the half-disc $\kn_{v}$ to be small enough such that the boundary of $\overline{\kn_{v}}$ intersects any half-edge in $H_v$ exactly once. % There exists a cyclic permutation  \[\kf_{\overline{\kS},v}:H_{v}\rightarrow H_{v}\] such that if the half-edge $b$ is not the last half-edge that the boundary of $\kn_v$ in the anti-clockwise direction then it maps to the half-edge that the boundary of $\overline{\kn_{v}}$ intersects in the anticlockwise direction after $b$. Otherwise, $b$ maps to the first edge that the boundary of $\overline{\kn_v}$ intersects in the anti-clockwise direction. We see this now.
We give the boundary of $\overline{\kn_v}$ an anti-clockwise orientation. Let $h$ be a half-edge incident to the vertex $v$. Following the orientation of $\overline{\kn_v}$, let $h_1$ be the half-edge in $H_v$ intersected by the boundary of $\overline{\kn_v}$ directly before $h$ and $h_2$ be the half-edge in $H_v$ intersected by the boundary of $\overline{\kn_v}$ directly after $h$. We say that $h_1$ is directly clockwise of $h$ at $v$ and $h_2$ is directly anti-clockwise of $h$ at $v$.
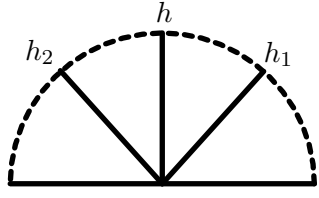
\begin{figure}[ht]
    \centering
\begin{tikzpicture}[line cap=round,line join=round,>=triangle 45,x=2cm,y=2cm]
\draw [shift={(0,0)},dashed,line width=2pt]  plot[domain=0:3.141592653589793,variable=\t]({1*1*cos(\t r)+0*1*sin(\t r)},{0*1*cos(\t r)+1*1*sin(\t r)});
\draw [line width=2pt] (-1,0)-- (1,0);
\draw [line width=2pt] (0,0)-- (-0.67,0.75);
\draw [line width=2pt] (0,0)-- (0,1);
\draw [line width=2pt] (0,0)-- (0.6667948594698256,0.7452413135250994);
\draw (-0.11,1.275) node[anchor=north west] {$h$};
\draw (0.6,1.01) node[anchor=north west] {$h_1$};
\draw (-0.97,1.025) node[anchor=north west] {$h_2$};
\end{tikzpicture}
\end{figure}

\subsection{Thick subcategories corresponding to \pac}
We can now simplify the arc collections we have to deal with and show that any thick subcategory in $\tds$ corresponds to a \pac.

\lem \label{accep} Let $v$ be the end point of some arc in a connected arc-collection $\ka_\kS$. Then there exists a collection $\ka_{\kS_1}$ of arcs such that: \begin{enumerate}[i)]
    \item $\td(\kS)=\td(\kS_1)$;
    \item $\ka_{\kS_1}$ is a \pac~at $v$ and
    \item every arc in $\ka_{\kS_1}$ can be written as a concatenation of arcs $\ka_{\kS}$.
\end{enumerate}

\pf Consider the ribbon graph $\Gamma(\kS)$. If $v$ is the only vertex in $\Gamma(\kS)_0$ then the result is trivial, so we assume otherwise. %For any vertex $p\in\Gamma(\kS)_0$ let $n_\kS(p)=d_\kS(q,p)$. 

Let $u\in\Gamma(\kS)_0$ be any vertex such that $d_\kS(v,u)=1$, such a vertex exists as $\ka_\kS$ is an connected arc-collection. Then $u\not=v$ and there exists an arc $\gamma_X\in\ka_{\kS}=\Gamma(\kS)_1$ incident to $v$ and $u$.

Let $h$ be the half-edge of $\gamma_X$ incident to $u$. If $\deg_\kS(u)>1$, then we consider a half-edge $h_1\in H_u$ which is directly clockwise or anti-clockwise of $h$ at $u$. We treat the anti-clockwise case but the clockwise case is similar.

Suppose $h_1$ is a half-edge of $\gamma_Y$. Then $\dhomd(X,Y)\not=0$ and there exist a morphism from $X$ to a shift of $Y$ corresponding to the vertex $u$ which does not factor through any object in $\kS$ \cite[Theorem 3.3]{Opper}. By \cite[Theorem 4.1]{Opper} the cone $\widehat{Y}$ of this morphism corresponds to the arc $\gamma_{\widehat{Y}}\simeq\gamma_X^\prime\gamma_Y^\prime$. Since $h_1$ is directly anti-clockwise of $h$, then $\gamma_{\widehat{Y}}$ is a simple concatenation with respect to $\ka_\kS$ and by Lemma~\ref{lemacscc} it does not internally intersect any arc in $\ka_\kS$. Therefore, by setting $\kS^\prime=\kS\setminus\{Y\}\cup\{\widehat{Y}\}$, it follows that $\ka_{\kS^\prime}$ is an connected arc-collection, every arc in $\ka_{\kS^\prime}$ can be written as a concatenation of arcs in $\ka_\kS$ and $\td(\kS)=\td(\kS^\prime)$. Note that, \[
\deg_{\kS^\prime}(w)=\begin{cases}
    \deg_\kS(w)& w\notin\{u,v\}\\
    \deg_\kS(w)-1& w=u\\
    \deg_\kS(w)+1& w=v
\end{cases}
\]
and $d_{\kS^\prime}(v,w)\leq d_\kS(v,w)$ for all $w\in\Gamma(\kS_1)_0=\Gamma(\kS)_0$. 

Working inductively, we can construct the collection $\widetilde{\kS}$ such that every object in $\widetilde{\kS}$ corresponds to an arc which can be written as a concatenation of arcs in $\ka_{\kS}$, $\deg_{\widetilde{\kS}}(u)=1$, $\deg_{\widetilde{\kS}}(w)=\deg_\kS(w)$ for all $w\notin\{u,v\}$, $d_{\widetilde{\kS}}(v,w)\leq d_\kS(v,w)$ for all $w\in\Gamma(\widetilde{\kS})_0=\Gamma(\kS)_0$ and $\td(\widetilde{\kS})=\td(\kS)$.

We can repeat the above process for any vertex $\widetilde{u}\in\Gamma(\widetilde{\kS})_0$ such that $d_{\widetilde{\kS}}(v,\widetilde{u})=1$. Again, working inductively, we construct a collection $\widehat{\kS}$ such that $\td(\widehat{\kS})=\td(\kS)$, every object in $\widehat{\kS}$ corresponds to an arc which can be written as a concatenation of arcs in $\ka_{\kS}$ and every vertex $x\in\Gamma(\widehat{\kS})_0$ with $d_{\widehat{\kS}}(v,x)=1$ has degree one.

If there existed a vertex $y\in\Gamma(\widehat{\kS})_0$ such that $d_{\widehat{\kS}}(v,y)=n>1$ then there would exist a shortest path $p\simeq\gamma_{Z_1}^\prime\gamma_{Z_2}^\prime...\gamma_{Z_n}^\prime$ from $v$ to $y$ on $\Gamma(\widehat{\kS})$ such that $Z_i\in\widehat{\kS}$ for all $1\leq i\leq n$. The path $p$ is the concatenation of $n$ distinct open paths and it passes through any vertex at most once (otherwise it could be shortened). Therefore, $\gamma_{Z_1}^\prime\gamma_{Z_2}^\prime$ is a concatenation of two distinct arcs at a vertex $v_1\not=v$. This implies $v_1$ is a vertex such that $d_{\widehat{\kS}}(v,v_1)=1$ but $\deg_{\widehat{\kS}}(v_1)>1$. This is a contradiction, so the shortest path from $v$ to every other vertex in $\Gamma(\widehat{\kS})_0$ has length one and the result follows. \qed

The next result follows from the previous one and Theorem~\ref{strESc}.
\cor $\sT\in\tds$ if and only if $\sT$ corresponds to an finite \pac.

\subsection{$\ka_{\kS}$-generated objects} 

In this section, we assume $\sT=\td(\kS)$ corresponds to the finite arc-collection $\ka_\kS=(\gamma_S)_{S\in\kS}$. We would like to identify which string objects in $\sD$ are in $\sT$. We now define the objects that, we will go on to prove, are the only string objects in $\sT$.

\defn A string object $X\in\sD$ is an $\ka_\kS$-generated object if \[\gamma_X\simeq\gamma_{S_1}^\prime\gamma_{S_2}^\prime...\gamma_{S_n}^\prime\] where $S_i\in\kS$ and $\gamma_{S_j}^\prime=\gamma_{S_j}^{\pm 1}$.

\lem \label{aginT} If $X\in \sD$ is an $\ka_\kS$-generated object then $X\in \sT$.

\pf If $X\in \sD$ and $\gamma_X\simeq\gamma_{S_1}^\prime\gamma_{S_2}^\prime...\gamma_{S_l}^\prime$ then we say that $\gamma_X$ is the concatenation of $l$ arcs in $\ka_\kS$. We now complete an induction argument on the number of arcs in $\ka_\kS$ that form the concatenation $\gamma_X$.

For $X\in \sD$ if $\gamma_X\simeq\alpha$ for some $\alpha\in\ka_\kS$ then some shift of $X$ is an object in $\kS$ and trivially $X\in \sT$. Now let $n\geq1$ and assume that for any $m\leq n$, if $\gamma_X\simeq\alpha_1^\prime\alpha_2^\prime...\alpha_m^\prime$ where %$\alpha_i^\prime\simeq\alpha_i^{\pm 1}$ and 
$\alpha_i\in\ka_\kS$ then $X\in \sT$. 

Consider $\gamma_X \simeq\alpha_1^\prime\alpha_2^\prime...\alpha_{n+1}^\prime$ where $\alpha_j\in\ka_\kS$ for all $1\leq j\leq n+1$. The arc $\gamma_X$ is the concatenation of $\gamma_Y\simeq\alpha_1^\prime\alpha_2^\prime...\alpha_{n}^\prime$ and $\gamma_Z\simeq\alpha_{n+1}^\prime$ at a common end point $u$. By \cite[Theorem 3.7]{Opper} there is a morphism between $Y$ and a shift of $Z$ corresponding to $u$ with the cone corresponding to $\gamma_X$ by \cite[Theorem 4.1]{Opper}. Since $Y$ and $Z$ are in $\sT$ by our assumption and $\sT$ is closed under taking cones of morphisms then $X\in \sT$. The result follows by induction.\qed

We have shown that any $\ka_\kS$-generated string object is in $\td(\kS)$. It remains to show any non-$\ka_\kS$-generated string object in $\sD$ is not in $\td(\kS)$. We can simplify the arguments by only considering exceptional and spherelike objects whose arcs form an arc-collection with the arcs in $\ka_\kS$. The following result allows us to make this simplification without loss of generality.

\lem \label{objac} %Let $\sT\subseteq\sT^\prime\subseteq\sD$ be a thick subcategories. If $X\in \sT^\prime$ is a string object then $\gamma_X\simeq\gamma_{Y_1}\gamma_{Y_2}...\gamma_{Y_n}$, where each $Y_i$ is an exceptional or spherelike string object in $\sT^\prime$ and each $\gamma_{Y_i}$ has no internal intersection points with any arc in $\ka_\kS$, if it is not homotopic to an arc in $\ka_\kS$.
Let $\sT\subseteq\sT^\prime\subseteq\sD$ be thick subcategories in $\tds$. Suppose $\sT$ corresponds to an arc collection $\ka_\kS$. Then, this can be extended to an arc-collection $\ka_{\kS^\prime}$ corresponding to $\sT^\prime$ such that $\ka_\kS\subseteq\ka_{\kS^\prime}$.

\pf Let $\sT=\td(\kS)$ where $\ka_\kS$ is an arc-collection. Since $\sT\subseteq\sT^\prime$, then $\sT^\prime=\td(\kS,\kS^\prime)$. Let $X$ be a string object in the $\kS^\prime$. By Lemma~\ref{strgen}, $\td(X)=\td(Z_1, Z_2,...,Z_m)$, where each $Z_j$ is exceptional or spherelike and $\gamma_X\simeq\gamma_{Z_1}\gamma_{Z_2}...\gamma_{Z_m}$.

Now we need to show that each $\gamma_{Z_k}\simeq\gamma_{V_1}\gamma_{V_2}...\gamma_{V_l}$, where each $\gamma_{V_r}$ has no internal intersection points with itself, any other $\gamma_{V_{r^\prime}}$ or any arc in $\ka_\kS$. If $\gamma_{Z_k}$ does not internally intersect any arc in $\ka_\kS$ then the result is trivial. Otherwise, let $W=Z_k$ and $\gamma_S$ be an arc in $\ka_\kS$ which internally intersects $\gamma_W$. There exists an internal intersection point $\gamma_S(t_1)=\gamma_W(t_2)$ between $\gamma_S$ and $\gamma_W$ such that \[t_1=\min\{t\in(0,1)\mid \gamma_S(t)=\gamma_W(t^\prime)\text{ for some } t^\prime\in(0,1)\}.\] We define the subpaths $\beta_1=\restr{\gamma_W}{[0,t_2]}$, $\beta_2=\restr{\gamma_W}{[t_2,1]}$, $\beta_3=\restr{\gamma_S}{[0,t_1]}$ and $\beta_4=\restr{\gamma_S}{t_1, 1]}$. The point $\gamma_S(t_1)$ corresponds to a morphism in each direction between $S$ and shifts of $W$. The cones of these morphisms are direct sums of objects corresponding to the arcs $\beta_1\beta_3\inv$, $\beta_1\beta_4$, $\beta_3\beta_2$, $\beta_4\inv\beta_2$ by \cite[Theorem 4.1]{Opper}. We see these paths in a neighbourhood of $\gamma_S(t_1)$.

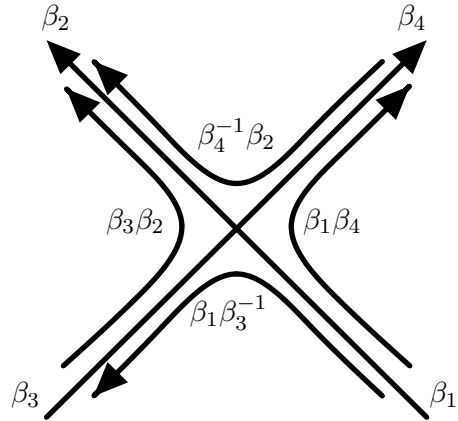
\begin{figure}[ht]
    \centering\begin{tikzpicture}[line cap=round,line join=round,>=triangle 45,x=1cm,y=1cm]
\draw[<-,line width=2pt, smooth,samples=100,domain=0.08:0.3213187959380341, variable=\t] plot({-2.1416468011261296*\t^3-3.5134928208247893*\t-0.9399999999999995,-1.6903506553817538*\t^3-3.560087238268793*\t-1.3199999999999996});
\draw[line width=2pt, smooth,samples=100,domain=0.3213187959380341:0.5366330268884635, variable=\t] plot({35.56814967038461*\t^3-36.350599191882466*\t^2+8.16663794313696*\t-2.1910151844916563,1.4835516030183598*\t^3-3.0595033562823946*\t^2-2.5770113036597597*\t-1.4252935918747445});
\draw[line width=2pt, smooth,samples=100,domain=0.5366330268884635:0.7562008707794534, variable=\t] plot({-34.19495282743511*\t^3+75.96095538372289*\t^2-52.10345154331899*\t+8.589958332826816,5.45616189102634*\t^3-9.455005006788314*\t^2+0.8550261055213951*\t-2.0392084663025187});
\draw[line width=2pt, smooth,samples=100,domain=0.7562008707794534:1, variable=\t] plot({2.2064666303265232*\t^3-6.61939989097957*\t^2+10.34388502468764*\t-7.150951764034594,-3.9962654527380215*\t^3+11.988796358214065*\t^2-15.360795159515034*\t+2.048264254038989});
\draw[<-,line width=2pt, smooth,samples=100,domain=0:0.2568395313984434, variable=\t] plot({2.267608961917259*\t^3+3.3545471478822066*\t-5.380000000000001,3.52221945872979*\t^3+3.349654461031808*\t-5.720000000000001});
\draw[line width=2pt, smooth,samples=100,domain=0.2568395313984434:0.49717792469183864, variable=\t] plot({-2.5215250610468227*\t^3+3.6901168147873102*\t^2+2.4067792743667176*\t-5.298858581163929,-30.926672946592507*\t^3+26.54351214773512*\t^2-3.4677687606613703*\t-5.136338738131821});
\draw[line width=2pt, smooth,samples=100,domain=0.49717792469183864:0.7375163179852338, variable=\t] plot({-1.9321645766784767*\t^3+2.811065747246417*\t^2+2.843824059824844*\t-5.371288254307748,31.203101676788204*\t^3-66.12514507873713*\t^2+42.60504192317546*\t-12.771800204968475});
\draw[line width=2pt, smooth,samples=100,domain=0.7375163179852338:1, variable=\t] plot({1.8590907626674875*\t^3-5.5772722880024626*\t^2+9.030360241597087*\t-6.892178716262112,-3.699588578932568*\t^3+11.098765736797704*\t^2-14.34885244191785*\t+1.2296752840527148});
\draw (-5.35,-3.1) node[anchor=north west] {$\beta_3\beta_2$};
\draw (-6.6,-5.4) node[anchor=north west] {$\beta_3$};
\draw (-1.1,-5.4) node[anchor=north west] {$\beta_1$};
\draw (-2.76,-3.1) node[anchor=north west] {$\beta_1\beta_4$};
\draw (-4.25,-4.25) node[anchor=north west] {$\beta_1\beta_3^{-1}$};
\draw [->,line width=2pt] (-1,-6)-- (-6,-1);
\draw [<-,line width=2pt] (-1,-1)-- (-6,-6);
\draw (-1.52,-0.4) node[anchor=north west] {$\beta_4$};
\draw (-4.15,-1.9) node[anchor=north west] {$\beta_4^{-1}\beta_2$};
\draw (-6.2,-0.4) node[anchor=north west] {$\beta_2$};
\draw[<-,line width=2pt, smooth,samples=100,domain=0.09:0.32131879593803414, variable=\t] plot({2.1416468011260017*\t^3+3.513492820824797*\t-6.06,-1.6903506553817798*\t^3-3.560087238268787*\t-1.3200000000000005});
\draw[line width=2pt, smooth,samples=100,domain=0.32131879593803414:0.5366330268884635, variable=\t] plot({-35.56814967038393*\t^3+36.350599191881706*\t^2-8.16663794313671*\t-4.808984815508368,1.4835516030184803*\t^3-3.059503356282537*\t^2-2.577011303659708*\t-1.4252935918747491});
\draw[line width=2pt, smooth,samples=100,domain=0.5366330268884635:0.7562008707794535, variable=\t] plot({34.19495282743431*\t^3-75.96095538372127*\t^2+52.10345154331797*\t-15.589958332826612,5.4561618910262215*\t^3-9.45500500678807*\t^2+0.8550261055212396*\t-2.0392084663024876});
\draw[line width=2pt, smooth,samples=100,domain=0.7562008707794535:1, variable=\t] plot({-2.206466630326236*\t^3+6.619399890978709*\t^2-10.343885024686797*\t+0.1509517640343227,-3.9962654527379824*\t^3+11.98879635821395*\t^2-15.360795159514922*\t+2.048264254038953});
\draw[<-,line width=2pt, smooth,samples=100,domain=0:0.2568395313984434, variable=\t] plot({2.267608961917259*\t^3+3.3545471478822066*\t-5.380000000000001,-3.5222194587297766*\t^3-3.349654461031806*\t-1.2800000000000002});
\draw[line width=2pt, smooth,samples=100,domain=0.2568395313984434:0.49717792469183864, variable=\t] plot({-2.5215250610468227*\t^3+3.6901168147873102*\t^2+2.4067792743667176*\t-5.298858581163929,30.926672946592475*\t^3-26.543512147735086*\t^2+3.4677687606613636*\t-1.863661261868179});
\draw[line width=2pt, smooth,samples=100,domain=0.49717792469183864:0.7375163179852338, variable=\t] plot({-1.9321645766784767*\t^3+2.811065747246417*\t^2+2.843824059824844*\t-5.371288254307748,-31.203101676788204*\t^3+66.12514507873712*\t^2-42.60504192317545*\t+5.771800204968474});
\draw[line width=2pt, smooth,samples=100,domain=0.7375163179852338:1, variable=\t] plot({1.8590907626674875*\t^3-5.5772722880024626*\t^2+9.030360241597087*\t-6.892178716262112,3.6995885789325733*\t^3-11.09876573679772*\t^2+14.348852441917865*\t-8.229675284052718});
\end{tikzpicture}
\end{figure}

Since $W, S \in \sT^\prime$ and a thick subcategory is closed under taking cones and direct summands then the objects $U_1$ and $U_2$ corresponding to $\beta_1\beta_3\inv$ and $\beta_3\beta_2$ respectively are in $\sT^\prime$. The arc homotopic to $\gamma_{U_1}\gamma_{U_2}\simeq\beta_1\beta_3\inv\beta_3\beta_2\simeq\beta_1\beta_2$ is homotopic to $\gamma_W$. Since $W$ and $S$ are exceptional or spherelike and $t_1$ is minimal then the paths $\beta_1$, $\beta_2$ and $\beta_3$ only intersect at common end points and they are paths with no self-intersection points. This implies the paths $\beta_1\beta_3\inv$ and $\beta_3\beta_2$ have no internal intersection points, hence the objects $U_1$ and $U_2$ are exceptional or spherelike. We would now like to show $\beta_1\beta_3\inv$ and $\beta_3\beta_2$ collectively have less intersections with arcs in $\ka_\kS$ than $\gamma_W$.

We can assume $\gamma_W$ and $\gamma_S$ are in minimal position and intersection points between them are transversal. There exists $\delta_1=\restr{\gamma_S}{[t_1, t_3]}$ for some $t_1<t_3<1$ such that there are no intersection points between $\gamma_S$ and $\gamma_W$ on $\delta_1$ apart from $\gamma_S(t_1)$. 

There exists a tubular neighbourhood $\kt$ of $\beta_3\delta_1$ since it is a subpath of $\gamma_S$, which has no self-intersection points. There exist a subpath $\chi=\restr{\gamma_W}{[t_4,t_5]}$ for some $0<t_4<t_2<t_5<1$ of $\gamma_W$ contained in $\kt$ which only intersects $\beta_3\delta_1$ at $\gamma_S(t_1)$. Now we see the tubular neighbourhood $\kt$. 
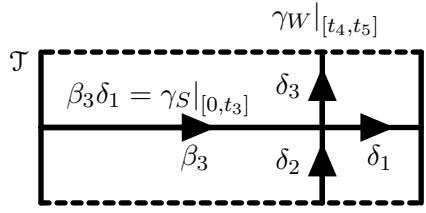
\begin{figure}[ht]
    \centering\begin{tikzpicture}[line cap=round,line join=round,>=triangle 45,x=1cm,y=1cm]
\draw [line width=2pt] (0,2)-- (0,0);
\draw [line width=2pt] (5,2)-- (5,0);
\draw [line width=2pt,dashed] (5,2)-- (0,2);
\draw [line width=2pt,dashed] (0,0)-- (5,0);
\draw [line width=2pt] (0,1)-- (5,1);
\draw [->,line width=2pt] (0,1)-- (2.3,1);
\draw [->,line width=2pt] (0,1)-- (4.65,1);
\draw [line width=2pt] (3.7,2)-- (3.7,0);
\draw [->,line width=2pt] (3.7,0)-- (3.7,0.8);
\draw [->,line width=2pt] (3.7,1)-- (3.7,1.8);
\draw (0.2,1.75) node[anchor=north west] {$\beta_3\delta_1=\restr{\gamma_S}{[0,t_3]}$};
\draw (1.7,.9) node[anchor=north west] {$\beta_3$};
\draw (4.15,.9) node[anchor=north west] {$\delta_1$};
\draw (2.95,.85) node[anchor=north west] {$\delta_2$};
\draw (2.95,1.85) node[anchor=north west] {$\delta_3$};
\draw (2.9,2.8) node[anchor=north west] {$\restr{\gamma_W}{[t_4,t_5]}$};
\draw (-0.55,2.15) node[anchor=north west] {$\mathcal{T}$};
\end{tikzpicture}
\end{figure}

As $\chi$ only intersects $\beta_3\delta_1$ at $\gamma_S(t_1)$ and $\beta_3$ is a subpath of the zero section, if $\delta_2=\restr{\gamma_W}{[t_4,t_2]}$ and $\delta_3=\restr{\gamma_W}{[t_2, t_5]}$ then $\delta_2\beta_3\inv$ and $\beta_3\delta_3$ are contained in $\kt^{+}$ or $\kt^{-}$. Applying Lemma~\ref{tbnhnsip} there exist a paths homotopic to $\beta_3\delta_1$, one of which does not internally intersect $\delta_2\beta_3\inv$ and another that does not internally intersect $\beta_3\delta_3$. Note that neither $\delta_2\beta_3\inv$ nor $\beta_3\delta_3$ internally intersect any arc in $\ka_\kS$ other than $\gamma_S$ as $\beta_3$ is a subpath of an arc in $\ka_\kS$ and by the constructions of $\delta_2$ and $\delta_3$.

The path $\beta_1\beta_3\inv$ is homotopic to $\phi_1\phi_2$ where $\phi_1=\restr{\gamma_W}{[0,t_4]}$, $\phi_2$ is a path homotopic to $\delta_2\beta_3\inv$ which intersects no arc in $\ka_\kS$ away from a marked point. Similarly, $\beta_3\beta_2$ is homotopic to $\psi_1\psi_2$ where $\psi_1\simeq\beta_3\delta_3$ intersects no arc in $\ka_\kS$ away from a marked point and $\psi_2=\restr{\gamma_W}{[t_5,1]}$. As $\phi_1$ and $\psi_2$ are subpaths of $\gamma_W$ with no intersection then each intersection point between an arc in $\ka_\kS$ and $\phi_1$ or $\psi_2$ corresponds to a distinct intersection point between $\gamma_W$ and the same arc in $\ka_\kS$. No intersection point between arcs in $\ka_\kS$ corresponds to $\gamma_S(t_1)$, therefore $\beta_1\beta_3\inv$ and $\beta_3\beta_2$ collectively have less internal intersection points with arcs in $\ka_\kS$ than $\gamma_W$. 

Working inductively, $\gamma_W\simeq\gamma_{V_1}\gamma_{V_2}...\gamma_{V_l}$, where each $V_r\in\sT^\prime$ and each $\gamma_{V_r}$ has no internal intersection points with itself or any arc in $\ka_\kS$. We can complete the argument for any exceptional or spherelike object. Hence, there exists a set $\{W_1,W_2,...,W_n\}$ of exceptional or spherelike objects such that $\sT^\prime=\td(\kS, W_1,W_2,...,W_n, \kS^\prime\setminus\{X\})$, where $\ka_\kS\cup\{\gamma_{W_1}, \gamma_{W_2},...,\gamma_{W_n}\}$ is an arc-collection. Using the above argument, $\sT^\prime=\td(\kS_1)$ where $\ka_{\kS_1}$ is an arc-collection and $\ka_{\kS_1}=\ka_\kS\cup\ka_{\kS^*}$ for some set of objects $\kS^*\subset\sT^\prime$. \qed

For the remainder of this chapter we refer to Setup~\ref{supac}. Recall that we defined terminal regions and the map $\tau$ in Setup~\ref{supac}. We prove a property of the non-terminal regions.
\lem \label{taufinite} For a non-terminal region $R_a$ there exists a minimal finite integer $n$ such that either: \begin{enumerate}[i)]
    \item $\tau^n(a)=a$ and $\tau$ is a cyclic permutation on $\{a, \tau(a), \tau^2(a),...,\tau^{n-1}(a)\}$ or
    \item $\tau^n(a)=c$ where $R_c$ is a terminal region.
\end{enumerate} Furthermore, $\tau^m(a)\not=\tau^l(a)$ for any $0\leq m,l<n$.

\pf As there are a finite number of half-edges in $\Gamma(\kS_1)$ then for a non-terminal region $R_a$ there exists a finite sequence of integers $\{a=a_0,a_1,...,a_n\}$ such that $\tau(a_l)=a_{l+1}$ and $R_{a_n}$ is either a terminal region or it is equal to some $R_{a_{k}}$ for some $0\leq k<n$.

We can take $n$ to be minimal with respect to the condition above, so $\tau^m(a)\not=\tau^l(a)$ for any $0\leq m,l<n$. 

If $R_{a_n}$ is a terminal region then the result follows.

Otherwise, $a_n=a_k$ for some $0\leq k<n$. By its definition, $\tau$ is an injective map. This implies that if $k\not=0$ then $a_{n-1}=a_{k-1}$ but this is a contradiction that $n$ is minimal. Therefore, $a_n=a_0$, by definition $\tau$ is a permutation on the set $\{a=a_0,a_1,...,a_{n-1}\}$ and the result follows. \qed

We distinguish two types of non-terminal region.
\defn Let $R_a$ be a non-terminal region and $\{a=a_0,a_1,...,a_n\}$ be a sequence satisfying the conditions of Lemma~\ref{taufinite}. If $a_n=a$ then is a cyclic region. Otherwise, $R_a$ is a terminating region.

We now prove a technical result, which constructs a path for each non-terminal region.

\lem \label{trpth} Let $\ka_{\kS_1}$ be an \pac~at $v$. For a non-terminal region $R_a$, there exists a closed path $\psi_{R_a}$ satisfying the following conditions: \begin{enumerate}[i)]
    \item $\psi_{R_a}$ has no internal intersection points with itself or any arc in $\ka_{\kS_1}$;
    \item the intersection of $\psi_{R_a}$ and $\kn_v$ is two radii of $\kn_v$ with one of these being in $R_a$ and the other being in either $R_a$ or a terminal region. 
\end{enumerate}

\pf For the non-terminal region $R_a$, there exists a finite sequence $\{a=~a_0,a_1,...,a_n\}$ of integers such that $\tau(a_i)=a_{i+1}$ for $0\leq i<n$, satisfying the conditions in Lemma~\ref{taufinite}.

We let $\psi_{R_a}\simeq\gamma^\prime_{S_{a_1}}\gamma^\prime_{S_{a_2}}...\gamma^\prime_{S_{a_{n-1}}}$. By the construction of the sequence $\{a_0,a_1,...,a_n\}$, $\psi_{R_a}$ is a sequence of simple concatenations with respect to $\ka_{\kS_1}$. %We see this as we start with a simple concatenation with respect to $\ka_{\kS_1}$ then continue taking simple concatenations, of arcs corresponding to $\ka_{\kS_1}$-generated objects, with respect to the an arc-collection that includes $\ka_{\kS_1}$. 

By Definition~\ref{ssc}, $\psi_{R_a}$ does not internally intersect itself or any arc in $\ka_{\kS_1}$ when in minimal position, so condition (i) is satisfied. By Lemma~\ref{taufinite}, $R_{a_n}$ is either a terminal region or equal to $R_a$, and the result follows from Lemma~\ref{sscreg}. \qed

Now we show that for any string object $X\in\sD$, if is not $\ka_\kS$-generated then $X\notin\td(\kS)$. We do this in the next two results. In the first, we prove the statement under additional assumptions and then show the general case reduces to this.

\lem\label{acnag} Suppose $\ka_{\kS_1}$ is a \pac~and $Y\not=0$ is a string object in $\sD$ such that either: \begin{enumerate}[a)]
    \item $\gamma_Y$ does not intersect any arc in $\ka_{\kS_1}$, or
    \item $\ka_{\kS_1}\cup\{\gamma_Y\}$ is a \pac.
\end{enumerate}
If $Y$ is not $\ka_{\kS_1}$-generated then $Y\notin\td(\kS_1)$.

\pf Suppose $Y\not=0$ is not $\ka_{\kS_1}$-generated.

Assume statement a) is true. By \cite[Theorem 3.3]{Opper}, $Y$ has no morphisms to the shift of any object in $\kS_1$. It follows from Remark~\ref{rzhom} that $Y\notin\td(\kS_1)$.

Now assume statement b) is true, i.e. $\ka_{\kS_1}\cup\{\gamma_Y\}$ is a pointed arc-collection at a marked point $v$. First, we consider the case when $Y$ is a spherelike object. Up to homotopy, we can assume that the intersection $\kn_v\cap\gamma_Y$ is two radii of $\kn_v$ and is contained in the (possibly non-distinct) regions $R_{c_1}$ and $R_{c_2}$ for some $0\leq c_1,c_2\leq b$. We consider the following two cases:

\underline{Case A:} $R_{c_1}$ or $R_{c_2}$ is a cyclic region. We treat the case where $R_{c_1}$ is a cyclic region and note that there is an analogous argument for when $R_{c_2}$ is cyclic.

There exist the sequence $\{c_1=d_0,d_1,...,d_n\}$ satisfying condition i) of Lemma~\ref{taufinite}. The path $\psi_{R_{c_1}}$ can be written as $h_{c_1+1}^\prime\phi_{S_{d_0+1}}^\prime\rho^\prime_1\phi_{S_{d_1+1}}^\prime\rho^\prime_2...\phi_{S_{d_{n-1}+1}}^\prime h^\prime_{c_1}$, where $h_{c_1+1}=\kh(c_1+1)$, $h_{c_1}=\kh(c_1)$, $S_i$ is the spherelike object that corresponds to $\kh(i)$, $\phi_{S_j}$ is the subarc of $\gamma_{S_j}$ outside of $\kn_v$ and $\rho_k$ is the chord of $\kn_v$ with end points on the half edges $\kh(d_k)$ and $\kh(d_k+1)$. 

Since $\gamma_Y$ does not internally intersect any arc in $\ka_{\kS_1}$, it follows that $\psi_{R_{c_1}}$ can only internally intersect $\gamma_Y$ in $R_{c_2}$ and this would be along some $\rho_{k^\prime}$ for $0<k^\prime<n-1$. By Lemma~\ref{trpth}, this happens if and only if $d_l=c_2$ for some $0<l\leq n-1$. Therefore, $\psi_{R_{c_1}}$ and $\gamma_Y$ internally intersect at most once, when in minimal position. 

Let $\chi$ be a closed curve homotopic to $\psi_{R_{c_1}}$ as in Section 2.1. The path $\psi_{R_{c_1}}$ is a sequence of simple concatenations with respect to $\ka_{\kS_1}$ so it does not internally intersect any arc in $\ka_{\kS_1}$. By Corollary~\ref{sscreg}, $\kn_v\cap\psi_{R_{c_1}}\subset R_{c_1}$ and so by construction, $\chi$ does not intersect any arc in $\ka_{\kS_1}$ when in minimal position. 

We now show $\chi$ intersects $\gamma_Y$ when in minimal position. Let $d_l=c_2$ for some $0<l\leq n-1$. Let $\chi$, $\psi_{R_{c_1}}$ and $\gamma_Y$ be in minimal position. As $\psi_{R_{c_1}}$ and $\gamma_Y$ internally intersect then $\gamma_Y$ has a subpath $\eta$ on the cylinder $\kC$ bounded by $\chi$ and $\psi_{R_{c_1}}$ with one end point being $\psi_{R_{c_1}}(t_1)$ for $0<t_1<1$. If $\eta$ had another end point on $\psi_{R_{c_1}}$ then, as $\gamma_Y$ has no internal self-intersection points and it only internally intersects $\psi_{R_{c_1}}$ once, there would be a bigon between $\gamma_Y$ and $\psi_{R_{c_1}}$. However, this would be a contradiction that $\gamma_Y$ and $\psi_{R_{c_1}}$ are in minimal position. Therefore, the other end point of $\eta$ must be on $\chi$. This implies $\chi$ and $\gamma_Y$ intersect, when in minimal position.

Otherwise, $d_l\not=c_2$ for any $0<l\leq n-1$. We let $\psi$ be a subpath of $\psi_{R_{c_1}}$ such that $\chi\simeq\psi^\prime\rho^\prime$, where $\rho$ is the chord of $\kn_v$ with end points on the half-edges $\kh(c_1)$ and $\kh(c_1+1)$. The path $\rho^\prime\psi^\prime$ intersects $\gamma_Y$ at most twice and these intersections are on $\rho$. As $\gamma_Y\cap\kn_v$ contains at least one radius of $\kn_v$ in $R_{c_1}$ then $\rho$ intersects $\gamma_Y$ at least once. If $\rho^\prime\psi^\prime$ and $\gamma_Y$ intersect once then the closed curve homotopic to $\phi_{R_{c_1}}$ and $\gamma_Y$ intersect in minimal position.

Otherwise, $\gamma_Y$ and $\rho^\prime\psi^\prime$ intersect twice and we have the following diagram.

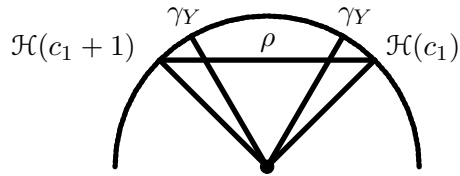
\begin{figure}[ht]
    \centering
    \begin{tikzpicture}[line cap=round,line join=round,>=triangle 45,x=1cm,y=1cm]
\draw [line width=2pt] (0,0)-- (1.4142135623730951,1.4142135623730951);
\draw [line width=2pt] (0,0)-- (-1.4142135623730951,1.4142135623730951);
\draw [line width=2pt] (0,0)-- (1,1.7320508075688772);
\draw [line width=2pt] (0,0)-- (-1,1.7320508075688772);
\draw [shift={(0,0)},line width=2pt,dash pattern=on 1pt off 1pt]  plot[domain=0:3.141592653589793,variable=\t]({1*2*cos(\t r)+0*2*sin(\t r)},{0*2*cos(\t r)+1*2*sin(\t r)});
\draw [line width=2pt] (1.4142135623730951,1.4142135623730951)-- (-1.4142135623730951,1.4142135623730951);
\draw (1.43,1.9) node[anchor=north west] {$\mathcal{H}(c_1)$};
\draw (-3.5,1.9) node[anchor=north west] {$\mathcal{H}(c_1+1)$};
\draw (-1.45,2.25) node[anchor=north west] {$\gamma_Y$};
\draw (0.8,2.25) node[anchor=north west] {$\gamma_Y$};
\draw (-0.225,1.9) node[anchor=north west] {$\rho$};
\begin{scriptsize}
\draw [fill=black] (0,0) circle (2.5pt);
\end{scriptsize}
\end{tikzpicture}
\end{figure}

If these two intersections caused $\gamma_Y$ and $\rho^\prime\psi^\prime$ to form a bigon then we have two cases:

\textit{Case 1:} The subpath of $\rho^\prime\psi^\prime$ between the intersection points which is not on $\rho$ is homotopic to a subpath of $\gamma_Y$. This implies that $\rho^\prime\psi^\prime$ is a closed curve homotopic to $\gamma_Y$ and therefore, since $\gamma_Y$ and $\phi_{R_{c_1}}$ have the same end point, $\gamma_Y\simeq\phi_{R_{c_1}}$. However, this is a contradiction as the object corresponding to $\phi_{R_{c_1}}$ is $\ka_{\kS_1}$-generated but $Y$ is not. Hence, we cannot have this case.

\textit{Case 2:} The subpath of $\gamma_Y$ between the two intersection points which includes the subpath outside of $\kn_v$ is homotopic to a subpath of $\rho$. This implies that $\gamma_Y$ is homotopic to a closed path contained in $\kn_v$ with end point $v$, so $\gamma_Y$ is homotopic to the trivial path at $v$. This is a contradiction, as $Y$ is non-zero and we cannot have this case either.

Therefore, we have shown, if we have Case A, that $\gamma_Y$ and $\chi$ always intersect in minimal position. There exists an object $P$ which is either an object in the unbounded homotopy category or a band object in $\sD$. By Lemma~\ref{hominf} or \cite[Theorem 3.3]{Opper}, $\hom_{\sK/\Sigma}(Y,P)\not=0$. As we have mentioned, $\chi$ does not intersect any arc in $\ka_{\kS_1}$, so $\hom_{\sK/\Sigma}(\kS_1,P)=0$. Therefore, it follows from Lemma~\ref{zhom} that $Y\notin\td(\kS_1)$.

\underline{Case B:} The region $R_{c_1}$ and $R_{c_2}$ are either terminal or terminating regions. First, we construct a spherelike object $Z$ such that $\ka_{\kS_1}\cup\{\gamma_Z\}$ is an arc-collection, if $Y\in\td(\kS_1)$ then $Z$ is in $\td(\kS_1)$ and either \begin{enumerate}[i)]
    \item $\kn_v\cap\gamma_Z\subset R_{c_2}\cup R_t$ such that $R_t$ is a terminal region or
    \item $\kn_v\cap\gamma_Z\subset R_{c_2}$
\end{enumerate}

If $R_{c_1}$ is a terminal region then we let $Z=Y$ and note that condition i) is satisfied.

Otherwise $R_{c_1}$ is a terminating region and there exist a sequence $\{c_1=f_0,f_1,...,f_m\}$ where $R_{f_m}$ is a terminal region. If $f_u=c_2$ for some $0<u\le n-1$ then the path $\gamma_{Y}^\prime\gamma_{S_{f_0+1}}^\prime\gamma_{S_{f_1+1}}^\prime...\gamma_{S_{f_{u-1}+1}}^\prime$, where $S_h$ is the object corresponding to $\kh(h)$ is a sequence of simple concatenations with respect to an arc collection $\ka_{\kS_1}\cup\{\gamma_Y\}$. By Definition~\ref{ssc} $\ka_{\kS_1}\cup\{\gamma_Y^\prime\gamma_{S_{f_0+1}}^\prime\gamma_{S_{f_1+1}}^\prime...\gamma_{S_{f_{u-1}+1}}^\prime\}$ is an arc-collection. If $Y\in\td(\kS_1)$ then $\td(\kS_1)=\td(\kS_1,Y)$ and by Lemma~\ref{aginT} the object corresponding $\gamma_Y^\prime\gamma_{S_{f_0+1}}^\prime\gamma_{S_{f_1+1}}^\prime...\gamma_{S_{f_{u-1}+1}}^\prime$ is in $\td(\kS_1)$. Finally, by Lemma~\ref{sscreg}, $\kn_v\cap\gamma_Y\subset R_{c_2}$. We let $Z$ correspond to $\gamma_Y^\prime\gamma_{S_{f_0+1}}^\prime\gamma_{S_{f_1+1}}^\prime...\gamma_{S_{f_{u-1}+1}}^\prime$ and condition ii) holds.

Otherwise, $R_{c_1}$ is a terminating region and $f_u\not=c_2$ for any $0<u\le n-1$. Then let $\gamma_Z\simeq\gamma_{Y}^\prime\psi^\prime_{R_{c_1}}$ and note by an similar argument to above that condition i) is satisfied.

Now we construct a spherelike object $W$ such that $\gamma_W\cap\kn_v\subset R_{t_1}\cup R_{t_2}$ where $R_{t_1}$ and $R_{t_2}$ are terminal regions, $\ka_{\kS_1}\cup\{\gamma_W\}$ is an arc-collection and if $Y\in\td(\kS_1)$ then $W\in\td(\kS_1)$.

If $R_{c_2}$ is a terminal then let $W=Z$ and note it satisfies the conditions we required for $W$.

Otherwise, $R_{c_2}$ is a terminating region and we have the path $\psi_{R_{c_2}}$. If $Z$ satisfies condition i) then by a similar argument we produced to show $Z$ satisfied a condition we required, we let $\gamma_W\simeq\gamma_Z^\prime\psi_{R_{c_2}}^\prime$ and say that $W$ satisfies the conditions we require. Similarly if $Z$ satisfies condition ii) we let $\gamma_W\simeq\psi_{R_{c_2}}^\prime\gamma_Z^\prime\psi_{R_{c_2}}^\prime$.

Now we show there is an object $V$ such that \begin{itemize}[(C1)]
    \item $\gamma_V$ has no oriented intersection points to any arc in $\ka_{\kS_1}$ and if $Y\in\td(\kS_1)$ then $V\in\td(\kS_1)$.
\end{itemize} 

If $t_1=b=t_2$ then let $V=W$. The only intersection point between $\gamma_V$ to an arc in $\ka_{\kS_1}$ is $v$ and by \cite[Theorem 3.7]{Opper}, $v$ is not an oriented intersection point from $\gamma_V$ to any arc in $\ka_{\kS_1}$.

If $t_1=b\not=t_2$ then there exists an exceptional object $E_1$ such that $\gamma_{E_1}^\prime\gamma_W^\prime$ is a simple concatenation with respect to $\ka_{\kS_1}\cup\{\gamma_W\}$ and by Lemma~\ref{exccc} this path has no oriented intersection points to any arc in $\ka_{\kS_1}$. We let $\gamma_V\simeq\gamma_{E_1}^\prime\gamma_W^\prime$. A similar result holds in the two other cases where $t_1\not=b$.

Therefore, if we have Case B then we can find an object $V$ satisfying condition (C1). By Remark~\ref{rzhom} $V\notin\td(\kS_1)$, so $Y\notin\td(\kS_1)$.

We have proven that there are no spherelike objects in $\td(\kS_1)$ which are not $\ka_{\kS_1}$-generated and satisfy the conditions of this lemma. We can write a similar arguement for exceptional objects and conclude there are no non-$\ka_{\kS_1}$-generated objects in $\td(\kS_1)$ which satisfy the conditions of this lemma. \qed

We now show that we can relax the assumptions in Lemma~\ref{acnag}
\lem \label{objacnip} Let $\ka_\kS$ be an arc-collection and suppose there exists a string object $X$ in $\td(\kS)$ which is not $\ka_\kS$-generated. Then there exists a pointed arc-collection $\ka_{\kS_1}$ and string object $Y$ in $\td(\kS_1)$ that is not $\ka_{\kS_1}$-generated and such that either \begin{enumerate}[a)]
    \item $\gamma_Y$ does not intersect any arc in $\ka_{\kS_1}$, or
    \item $\ka_{\kS_1}\cup\{\gamma_Y\}$ is a \pac.
\end{enumerate}

\pf Let $\sT=\td(\kS)$ and $v$ be the end point of some arc in $\ka_{\kS}$. By Lemma~\ref{accep}, there exists a \pac~at $v$ such that $\td(\kS_1)=\sT$. Furthermore, each object in $\kS_1$ is $\ka_\kS$-generated, so any string object in $\sD$ that is $\ka_{\kS_1}$-generated is also $\ka_\kS$-generated. 

Suppose $X$ is a string object in $\sD$ which is not $\ka_\kS$-generated, so $X$ is not $\ka_{\kS_1}$-generated. We now construct the desired $Y$ in a number of steps.

\textit{Step 1:} We replace $X$ with an object which is exceptional or spherelike. By Lemma~\ref{strgen}, there exist exceptional or spherelike objects $Z_1,Z_2,...,Z_n$ for a finite $n\in\IN$ such that \begin{equation}
    \label{objacnip1} \td(X)=\td(Z_1,Z_2,...,Z_n)
\end{equation}and \begin{equation} \label{objacnip2}
    \gamma_X\simeq\gamma_{Z_1}^\prime\gamma_{Z_2}^\prime...\gamma_{Z_n}^\prime
\end{equation} Since $X$ is not $\ka_{\kS_1}$-generated, it follows from (\ref{objacnip2}) that there exists $Z:=Z_j$ which is not $\ka_{\kS_1}$-generated for some $j\in\{1,2,...,n\}$. It follows from (\ref{objacnip1}) that if $X\in\sT$, then $Z\in\sT$.

\textit{Step 2:} We now replace $Z$ with an string object corresponding to an arc which does not internally intersect any arc in $\ka_\kS$. By Lemma~\ref{objac}, if $Z\in\sT$ then there exists exceptional or spherelike objects $W_1, W_2,..., W_m$ such that \begin{equation}\label{objacnip3}
    \gamma_Z\simeq\gamma_{W_1}^\prime\gamma_{W_2}^\prime...\gamma_{W_m}^\prime
\end{equation} where each $\gamma_{W_i}$ is either in $\ka_{\kS_1}$ or has no internal intersection points with any arc in $\ka_{\kS_1}$. Since $Z$ is not $\ka_{\kS_1}$-generated, it follows from (\ref{objacnip3}) that there exists $k\in\{1,2,...,m\}$ such that $W:=W_k$ is not $\ka_{\kS_1}$-generated. We note that if $X\in\sT$ then $Z\in\sT$ and so $W\in \sT$ by construction.

\textit{Step 3:} If $\gamma_W$ satisfies either condition a) or b), then taking $Y:=W$ the result follows. Otherwise, $\gamma_W$ has a common end point $v_1\not=v$ with some arc in $\ka_{\kS_1}$. There exists an unique object $S_1\in\kS_1$ such that $\gamma_{S_1}$ has the end points $v$ and $v_1$. This implies there is some morphism between $W$ and $S_1$ upto shift corresponding to $v_1$ which does not factor through a shift of $W$ or any object in $\kS_1$ \cite[Theorem 3.7]{Opper}. We can choose such a morphism $f$. By \cite[Theorem 4.1]{Opper}, the cone $C_1$ of $f$ corresponds to the arc $\gamma_{C_1}\simeq\gamma_W^\prime\gamma_{S_1}^\prime$. As $f$ does not factor through a shift of any object in $\kS_1\cup \{W\}$, then $\gamma_{C_1}$ is a simple concatenation with respect to the arc-collection $\ka_{\kS_1}\cup\{\gamma_W\}$, so $\ka_{\kS_1}\cup\{\gamma_{C_1}\}$ is an arc-collection. We note two things, if $X\in\sT$ then $W$ is in $\sT$ and so is $C_1$. Also, $\gamma_{C_1}$ has an end point at $v$. We have the following three possible cases:

\textit{Case 1:} The arc $\gamma_{C_1}$ is a closed arc at $v$ and the result follows if we take $Y:=C_1$.

\textit{Case 2:} The arc $\gamma_{C_1}$ is open and it has an end point which is not an end point of any arc in $\ka_{\kS_1}$. Again, by taking $Y:=C_1$, the result follows.

\textit{Case 3:} The arc $\gamma_{C_1}$ is open and it has an end point $v_2$ not equal to $v$, in common with some arc in $\ka_{\kS_1}$. There exists a unique object $S_2\in\kS_1$ such that $\gamma_{S_2}$ has the end points $v$ and $v_2$. By \cite[Theorem 3.7]{Opper}, there exist a morphism between $C_1$ and $S_2$ upto shift corresponding to $v_2$. The cone $C_2$ of this morphism corresponds to the arc $\gamma_{C_2}\simeq\gamma_{S_2}^\prime\gamma_{C_1}^\prime$. The arc $\gamma_{C_2}$ is a simple concatenation at $v_2$ with respect to the arc-collection $\ka_{\kS_1}\cup\{\gamma_{C_1}\}$ and $\ka_{\kS_1}\cup\{\gamma_{C_2}\}$ is an arc-collection, since $\gamma_{S_2}$ and $\gamma_{C_1}$ are open arcs and the only arcs in $\ka_{\kS_1}\cup\{C_1\}$ with an end point at $v_2$. As the other end point of both $\gamma_{S_2}$ and $\gamma_{C_1}$ is $v$ then $\gamma_{C_2}$ is a closed arc with end point $v$. By setting $Y:=C_2$ the result follows. \qed

\cor \label{nalg}Let $\ka_\kS$ be an arc-collection and $X$ be a string object in $\sD$. If $X$ is not $\ka_\kS$-generated then $X\notin\td(\kS)$. 

\pf Consider the string object $X$ which is not $\ka_\kS$-generated and suppose, to the contrary, that $X\in\td(\kS)$. Applying Lemma~\ref{objacnip}, there exists an arc-collection $\ka_{\kS_1}$ and the object $Y\in\td(\kS_1)$, which is not $\ka_{\kS_1}$-generated, satisfying the conditions of Lemma~\ref{acnag}. This is a contradiction and $Y$ cannot be in $\td(\kS_1)$. Therefore, $X\notin\td(\kS)$.\qed 

\subsection{Partial order on $\tds$}
We denote by $\acs$ the set of arc-collections on a geometric model $\Sigma$. We can define a binary relation $\leqgen$ on $\acs$ by saying $(\alpha_i\mid i\in I)\leqgen(\beta_j\mid j\in J)=\beta$  if and only if every arc in $(\alpha_i\mid i\in I)$ corresponds to a $\beta$-generated object.

\lem The binary operation $\leqgen$ is a well-defined pre-order on $\acs$.

\pf The operation is clearly reflexive. 

Suppose $(\alpha_i\mid i\in I)\leqgen(\beta_j\mid j\in J)$ and $(\beta_j\mid j\in J)\leqgen(\delta_k\mid k\in K)$. Each $\alpha_{i^\prime}\simeq\beta_{i_0}^\prime\beta_{i_1}^\prime...\beta_{i_n}^\prime$ and each $\beta_{i_m}\simeq\delta_{m_0}\delta_{m_1}...\delta_{m_t}$. Every $\alpha_i\simeq\delta_{0_0}^\prime\delta_{0_1}^\prime...\delta_{0_{r_0}}^\prime\delta_{1_0}^\prime\delta_{1_1}^\prime...\delta_{1_{r_1}}^\prime...\delta_{n_0}^\prime\delta_{n_1}^\prime...\delta_{n_{r_n}}^\prime$, hence we have $(\alpha_i\mid i\in I)\leqgen(\delta_k\mid k\in K)$. The operation is transitive and we have a well-defined pre-order on $\acs$.\qed

\defn We say $(\alpha_i\mid i\in I)\simgen(\beta_j\mid j\in J)$ if and only if $(\alpha_i\mid i\in I)\leqgen(\beta_j\mid j\in J)$ and $(\beta_j\mid j\in J)\leqgen(\alpha_i\mid i\in I)$. 

Hence, on the set of equivalence classes $\acs/\simgen$ the antisymmetric property holds for the binary operation $\leqgen$ and we have the following:

\cor The binary operation $\leqgen$ is a well-defined partial order on $\acs/\simgen$.

Now we see main result of the paper and show $\leqgen$ induces a partial order on the thick subcategories generated by string objects of the derived category of a gentle algebra with the next couple of results.

\lem \label{pots} $\td(\kS)\subseteq\td(\kS^\prime)$ if and only if $(\gamma_A)_{A\in\kS}\leqgen(\gamma_B)_{B\in\kS^\prime}$.

\pf ($\impliedby$) Let $(\gamma_A)_{A\in\kS}\leqgen(\gamma_B)_{B\in\kS^\prime}$, so any arc in $(\gamma_A)_{A\in\kS}$ can be written as a concatenation of arcs in $(\gamma_B)_{B\in \kS^\prime}$. By Lemma~\ref{aginT}, every object in $\kS$ is in $\td(\kS^\prime)$, hence $\td(\kS)\subseteq\td(\kS^\prime)$.

($\implies$) Assume $\td(\kS)\subseteq\td(\kS^\prime)$, by Lemmas~\ref{aginT} and \ref{nalg}, every object in $\kS$ corresponds to an arc which is a concatenation of the arcs in $(\gamma_B)_{B\in\kS^\prime}$, hence $(\gamma_A)_{A\in\kS}\leqgen (\gamma_B)_{B\in\kS^\prime}$.\qed 

We denote the the equivalence class of $\acs/\simgen$ containing the arc-collection $\ka$ by $[\ka]$. We define a map $\varphi:\acs/\simgen\rightarrow\tds$ such that $\varphi([(\gamma_A)_{A\in\kS}])=\td(\kS)$.

\lem \label{iacts} The map $\varphi$ is a well-defined bijection.

\pf First show $\varphi$ is a well-defined map. For two arc collections $(\gamma_A)_{A\in\kS}$ and $(\gamma_B)_{B\in\kS^\prime}$ in $\acs$, if $[(\gamma_A)_{A\in\kS}]=[(\gamma_B)_{B\in\kS^\prime}]$ then by the antisymmetric condition of our binary operation $(\gamma_A)_{A\in\kS}\leqgen(\gamma_B)_{B\in\kS^\prime}$ and $(\gamma_B)_{B\in\kS^\prime}\leqgen(\gamma_A)_{A\in\kS}$. By Lemma~\ref{pots} $\td(\kS)\subseteq\td(\kS^\prime)$ and $\td(\kS^\prime)\subseteq\td(\kS)$, hence $\td(\kS)=\td(\kS^\prime)$.

Next we see $\varphi$ is injective, if $\td(\kS)=\td(\kS^\prime)$ then by Lemma~\ref{pots} $(\gamma_A)_{A\in\kS}\leqgen(\gamma_B)_{B\in\kS^\prime}$ and $(\gamma_B)_{B\in\kS^\prime}\leqgen(\gamma_A)_{A\in\kS}$. By our defined antisymmetric property $(\gamma_A)_{A\in\kS}\simgen(\gamma_B)_{B\in\kS^\prime}$ so they are in the same equivalence class in $\acs/\simgen$.

Finally $\varphi$ is surjective since every thick subcategory in $\tds$ corresponds to an arc-collection from Theorem~\ref{strESc}.\qed

Lemmas~\ref{pots} and \ref{iacts} give us the following theorem.
\thm\label{mt} There is an isomorphism of posets:\[\faktor{\acs}{\simgen}\xleftrightarrow{1-1}\tds.\]

\section{Examples}

\exm \label{exm1} Remark~\ref{posnlat} says that the set $\tds$ is not a sublattice of the lattice of all thick subcategories. Consider the example where $\sD=\sD^b(\text{mod } \sk Q)$ and $Q$ is given by 
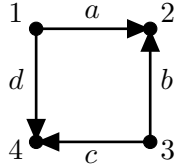
\begin{figure}[ht]
    \centering
\definecolor{ududff}{rgb}{0,0,0}
\begin{tikzpicture}[line cap=round,line join=round,>=triangle 45,x=1cm,y=1cm]
\draw [->,line width=1pt] (-1,4) -- (0.5,4);
\draw [->,line width=1pt] (0.5,2.5) -- (0.5,4);
\draw [->,line width=1pt] (0.5,2.5) -- (-1,2.5);
\draw [->,line width=1pt] (-1,4) -- (-1,2.5);
\draw (-1.5,4.475) node[anchor=north west] {$1$};
\draw (0.5,4.475) node[anchor=north west] {$2$};
\draw (0.5,2.65) node[anchor=north west] {$3$};
\draw (-1.5,2.65) node[anchor=north west] {$4$};
\draw (-0.5,4.45) node[anchor=north west] {$a$};
\draw (0.5,3.6) node[anchor=north west] {$b$};
\draw (-0.5,2.525) node[anchor=north west] {$c$};
\draw (-1.5,3.6) node[anchor=north west] {$d$};
%\draw (-2.3,3.6) node[anchor=north west] {$Q:$};
\begin{scriptsize}
\draw [fill=ududff] (-1,4) circle (2.5pt);
\draw [fill=ududff] (0.5,4) circle (2.5pt);
\draw [fill=ududff] (-1,2.5) circle (2.5pt);
\draw [fill=ududff] (0.5,2.5) circle (2.5pt);
\end{scriptsize}
\end{tikzpicture} 
\end{figure}

We have the exceptional objects $A, B, C$ and $D$, where $A$ is the projective module corresponding to the vertex $2$, $B$ is the string object corresponding to $\boldsymbol{d}\boldsymbol{\overline{c}}$, C is the projective module corresponding to the vertex 4 and $D$ is the string object corresponding to $\boldsymbol{b\overline{a}}$. We see the arcs corresponding to these objects on the geometric model of $\sD$.
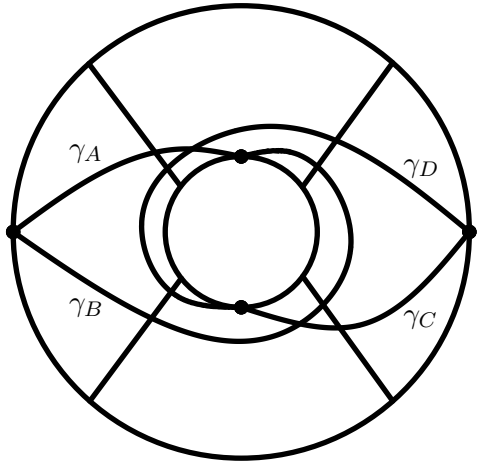
\begin{figure}[ht]
    \centering
\begin{tikzpicture}[line cap=round,line join=round,>=triangle 45,x=1cm,y=1cm]
\draw [line width=2pt] (0,0) circle (3cm);
\draw [line width=2pt] (0,0) circle (1cm);
\draw[line width=2pt, smooth,samples=100,domain=0:0.794806422143927, variable=\t] plot({0.1713384710219546*\t^3+2.835875633250096*\t-3,-1.177223425692723*\t^3+2.1276571345160034*\t});
\draw[line width=2pt, smooth,samples=100,domain=0.794806422143927:1, variable=\t] plot({-0.663670464502017*\t^3+1.9910113935060512*\t^2+1.2534069911297574*\t-2.5807479201337915,4.5599123949929155*\t^3-13.679737184978746*\t^2+13.000400102378197*\t-2.880575312392365});
\draw[line width=2pt, smooth,samples=100,domain=0:0.5453219084985116, variable=\t] plot({-1.8519220336833437*\t^3+7.519078284238384*\t-3.0000000000000044,10.349679882793524*\t^3-5.278281337251183*\t});
\draw[line width=2pt, smooth,samples=100,domain=0.5453219084985116:0.7077928295334395, variable=\t] plot({-54.38255255338643*\t^3+85.93831106890394*\t^2-39.344965520995075*\t+5.518663269275919,-20.87079506780528*\t^3+51.07562695287154*\t^2-33.13093970494911*\t+5.062888272610029});
\draw[line width=2pt, smooth,samples=100,domain=0.7077928295334395:0.9120354159661325, variable=\t] plot({36.912486266009005*\t^3-107.91561047613132*\t^2+97.86345012551881*\t-26.85304764613952,-72.28479483008542*\t^3+160.24700806099906*\t^2-110.40166044354416*\t+23.29344229649615});
\draw[line width=2pt, smooth,samples=100,domain=0.9120354159661325:1, variable=\t] plot({26.219363373090907*\t^3-78.65809011927273*\t^2+71.1795553767137*\t-18.740828630531876,142.22341530541865*\t^3-426.67024591625596*\t^2+424.88766142530187*\t-139.44083081446448});
\draw[line width=2pt, smooth,samples=100,domain=0:0.5513695644544098, variable=\t] plot({-0.7116077650226995*\t^3-2.431617008731506*\t+3,2.5245379620799304*\t^3-3.0164260551767392*\t});
\draw[line width=2pt, smooth,samples=100,domain=0.5513695644544098:1, variable=\t] plot({0.8745703197461061*\t^3-2.6237109592383185*\t^2-0.9849826398820126*\t+2.734123279374225,-3.10267268182072*\t^3+9.30801804546216*\t^2-8.148583910836997*\t+0.9432385471955578});
\draw[line width=2pt, smooth,samples=100,domain=0:0.47728429654334176, variable=\t] plot({-1.281951657139702*\t^3-6.1611479341635595*\t+2.9999999999999982,-10.166026417086243*\t^3+5.1233748563694625*\t});
\draw[line width=2pt, smooth,samples=100,domain=0.47728429654334176:0.712616395462456, variable=\t] plot({41.26020120201809*\t^3-60.914104502467325*\t^2+22.912197582864163*\t-1.6254170877520322,8.855331428798074*\t^3-27.235786196316216*\t^2+18.122587911883105*\t-2.068106752939283});
\draw[line width=2pt, smooth,samples=100,domain=0.712616395462456:0.8756998712268171, variable=\t] plot({-26.58801354712853*\t^3+84.13514619683123*\t^2-80.45227661500012*\t+22.927655919499323,77.39909220336175*\t^3-173.77200940014757*\t^2+122.54670309607931*\t-26.872885606912043});
\draw[line width=2pt, smooth,samples=100,domain=0.8756998712268171:1, variable=\t] plot({-38.309925660407785*\t^3+114.92977698122337*\t^2-107.4191308273697*\t+30.799279506554125},{-79.27885545602061*\t^3+237.83656636806185*\t^2-237.8988736999954*\t+78.34116278795415});
\draw [line width=2pt] (-0.7966506514071133,0.6044400215179518)-- (-2.004988051082459,2.231596494668461);
\draw [line width=2pt] (0.8038687976796928,0.5948066543987257)-- (1.9908956305686534,2.244177931489088);
\draw [line width=2pt] (0.8076897832811484,-0.589607678022642)-- (1.992418475572746,-2.2428260338680692);
\draw [line width=2pt] (-1.9916583707595037,-2.2435010439452885)-- (-0.7888548296293945,-0.6145795780616038);
\begin{scriptsize}
\draw [fill=black] (0,1) circle (2.5pt);
\draw [fill=black] (0,-1) circle (2.5pt);
\draw [fill=black] (-3,0) circle (2.5pt);
\draw [fill=black] (3,0) circle (2.5pt);
\draw [fill=black] (-3,0) circle (2.5pt);
\draw [fill=black] (0,1) circle (2.5pt);
\draw [fill=black] (-3,0) circle (2.5pt);
\draw [fill=black] (0,1) circle (2.5pt);
\draw [fill=black] (3,0) circle (2.5pt);
\draw [fill=black] (0,-1) circle (2.5pt);
\draw [fill=black] (3,0) circle (2.5pt);
\draw [fill=black] (0,-1) circle (2.5pt);
\end{scriptsize}
\draw (-2.4,1.3) node[anchor=north west] {$\gamma_A$};
\draw (-2.4,-0.75) node[anchor=north west] {$\gamma_B$};
\draw (2,-0.9) node[anchor=north west] {$\gamma_C$};
\draw (2,1.1) node[anchor=north west] {$\gamma_D$};
\end{tikzpicture}
\end{figure}

Consider the thick subcategories $\td(A,B)$ and $\td(C,D)$. By Lemma~\ref{pots} there is no non-zero thick subcategory in $\tds$ that is contained both of these. However, both thick subcategories contain $\td(X)$ where $X$ is the band object corresponding to $\boldsymbol{a\overline{b}c\overline{d}}$.

\begin{exm} \label{exm2}
Remark~\ref{nsphexcbds} says not all thick subcategories are not generated by exceptional and spherelike objects we will see an example of such a thick subcategory. Consider the derived category $\sD$ corresponding to the gentle algebra $\sk Q/I$, where

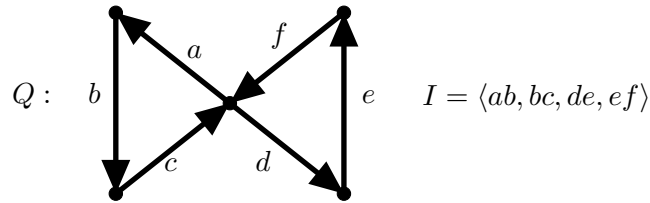
\begin{figure}[ht]
    \centering   
    \definecolor{ududff}{rgb}{0,0,0}
\begin{tikzpicture}[line cap=round,line join=round,>=triangle 45,x=2cm,y=2cm]
\draw [->,line width=2pt] (0,0) -- (-0.75,0.6);
\draw [->,line width=2pt] (-0.75,0.6) -- (-0.75,-0.6);
\draw [->,line width=2pt] (-0.75,-0.6) -- (0,0);
\draw [->,line width=2pt] (0,0) -- (0.75,-0.6);
\draw [->,line width=2pt] (0.75,-0.6) -- (0.75,0.6);
\draw [->,line width=2pt] (0.75,0.6) -- (0,0);
\draw (-0.35,0.45) node[anchor=north west] {$a$};
\draw (-1,0.2) node[anchor=north west] {$b$};
\draw (-1.5,0.2) node[anchor=north west] {$Q:$};
\draw (-0.5,-0.3) node[anchor=north west] {$c$};
\draw (0.1,-0.25) node[anchor=north west] {$d$};
\draw (0.8,0.15) node[anchor=north west] {$e$};
\draw (0.2,0.6) node[anchor=north west] {$f$};
\draw (1.2,0.2) node[anchor=north west] {$I=\langle ab, bc, de, ef\rangle$};
\begin{scriptsize}
\draw [fill=ududff] (0,0) circle (2.5pt);
\draw [fill=ududff] (0.75,0.6) circle (2.5pt);
\draw [fill=ududff] (0.75,-0.6) circle (2.5pt);
\draw [fill=ududff] (-0.75,0.6) circle (2.5pt);
\draw [fill=ududff] (-0.75,-0.6) circle (2.5pt);
\end{scriptsize}
\end{tikzpicture}
\end{figure}

In this example, any band object corresponds to a closed curve with at least one self-intersection point. It follows from \cite[Theorem 3.3]{Opper} that for any band object $B\in\sD$, $\dhomd(B,B)\geq4$. From Lemma~\ref{bdgnstr} we know a thick subcategory generated by band objects contains no string objects. This implies any thick subcategory of $\sD$ generated by band objects contains no exceptional or spherelike objects. Therefore, there exists thick subcategories generated by band objects which are not equivalent to any thick subcategory generated by exceptional or spherelike objects.
\end{exm}

\Addresses
\end{document}